\documentclass[11pt,reqno]{amsart}

\usepackage[T1]{fontenc}

\usepackage{amsmath}
\usepackage{mathtools}
\usepackage{amsthm}
\usepackage{enumitem}
\usepackage{xcolor}
\usepackage{cases}

\usepackage{amssymb}
\usepackage[a4paper,margin=1in]{geometry}

\usepackage[numbers,sort&compress]{natbib}
\usepackage{indentfirst}
\usepackage[colorlinks=true,bookmarksnumbered=true,pagebackref]{hyperref}
\hypersetup{linktoc=all}

\setenumerate[1]{itemsep=0pt,parsep=0pt,label={\textnormal{(\roman*)}}}
\allowdisplaybreaks

\newtheorem{theorem}{Theorem}[section]
\newtheorem{lemma}[theorem]{Lemma}

\newtheorem{proposition}[theorem]{Proposition}
\theoremstyle{definition}

\theoremstyle{remark}
\newtheorem{remark}[theorem]{Remark}
\numberwithin{equation}{section}

\DeclareMathOperator{\dist}{dist}

\usepackage{microtype}
\usepackage{anyfontsize}

\usepackage{nccmath}

\def\ds{\displaystyle}
\def\lag{\langle}
\def\ol{\overline}
\def\rag{\rangle}
\def\what{\widehat}
\def\wtilde{\widetilde}
\def\ve{\varepsilon}
\def\vp{\varphi}

\newcommand{\mbH}{\mathbb{H}}
\newcommand{\mbR}{\mathbb{R}}
\newcommand{\mbS}{\mathbb{S}}
\newcommand{\mbRn}{{\mathbb{R}^n}}

\newcommand{\mcE}{\mathcal{E}}
\newcommand{\mcI}{\mathcal{I}}
\newcommand{\mcM}{\mathcal{M}}
\newcommand{\mcT}{\mathcal{T}}
\newcommand{\mcU}{\mathcal{U}}
\newcommand{\mcV}{\mathcal{V}}
\newcommand{\mcZ}{\mathcal{Z}}
\newcommand{\mrc}{\mathrm{c}}
\newcommand{\mrd}{\mathrm{d}}
\newcommand{\mre}{\mathrm{e}}
\newcommand{\mrloc}{\mathrm{loc}}
\newcommand{\nit}{{-\infty}}
\newcommand{\pit}{{+\infty}}
\newcommand{\Ccit}{C_\mrc^\infty}
\newcommand{\limktoit}{\lim\limits_{k\to\infty}}
\newcommand{\norm}[1]{\lVert #1\rVert}
\newcommand{\lrnorm}[1]{\left\lVert #1\right\rVert}
\newcommand{\bignorm}[1]{\bigl\lVert #1\bigr\rVert}
\newcommand*{\diop}{\mathop{}\mspace{-3mu}\mathrm{d}}
\newcommand{\diopr}{\diop r}
\newcommand{\diopt}{\diop t}
\newcommand{\diopv}{\diop v}
\newcommand{\diopx}{\diop x}
\newcommand{\diopy}{\diop y}
\newcommand{\diopz}{\diop z}
\newcommand{\diopomega}{\diop\omega}
\newcommand{\dioptau}{\diop\tau}
\newcommand{\dioptheta}{\diop\theta}
\DeclareMathOperator{\supp}{supp}
\DeclareMathOperator{\MaOpdiv}{div}
\DeclareMathOperator{\MaOpspan}{span}

\newcommand{\mbBn}{{\mathbb{B}^n}}

\hypersetup{
  pdftitle={Critical Poincare--Sobolev stability inequalities on hyperbolic spaces: sharp constants and extremals},
  pdfauthor={Jungang Li, Yansheng Shen and Zeming Wang},
  pdfsubject={Quantitative stability of the critical Poincare--Sobolev inequality on hyperbolic space},
  pdfkeywords={Poincare--Sobolev inequality, Bianchi--Egnell stability, sharp stability constant, existence of extremals, hyperbolic space, critical exponent}
}

\makeatletter
\renewcommand{\l@section}{\@tocline{1}{4pt plus 1pt}{0pt}{}{}}
\renewcommand{\l@subsection}{\@tocline{2}{0pt}{2em}{}{}}
\renewcommand{\l@subsubsection}{\@tocline{3}{0pt}{4em}{}{}}
\makeatother

\title[Critical Poincar\'e--Sobolev stability on hyperbolic spaces]{Critical Poincar\'e--Sobolev stability inequalities on hyperbolic spaces: sharp constants and extremals}

\author{Jungang Li}
\address{School of Mathematical Sciences, University of Science and Technology of China, Hefei 230026, The People's Republic of China}
\email{jungangli@ustc.edu.cn}

\author{Yansheng Shen}
\address{School of Mathematical Sciences, Jiangsu University, Zhenjiang 212013, The People's Republic of China}
\email{ysshen@mail.bnu.edu.cn}

\author[Zeming Wang]{Zeming Wang$^{*}$}
\address{School of Mathematics and Statistics, Beijing Institute of Technology, Beijing 102488, The People's Republic of China}
\email{zmwangmath@bit.edu.cn}
\thanks{$^{*}$Corresponding author.}

\begin{document}

\begin{abstract}
We quantify Bianchi--Egnell-type stability for critical Poincar\'e--Sobolev inequalities on hyperbolic spaces $\mathbb{B}^n$ with $n\geqslant 4$, and show the existence of extremal functions.
  Precisely, let $c_{\mathrm{BE}}$ be the sharp constant in the stability inequality
  \[
    \|f\|_\lambda^2-S_\lambda\|f\|_{L^{2^\ast}(\mathbb{B}^n)}^2
    \geqslant c_{\mathrm{BE}}\operatorname{dist}(f,\mathcal{Z})^2\quad\text{for every }f\in H^1(\mathbb{B}^n),
  \]
  which was established by Bhakta--Ganguly--Karmakar--Mazumdar \cite{BGKM2025}. Here $\frac{n(n-2)}4<\lambda<\frac{(n-1)^2}{4}$, $2^\ast:=\frac{2n}{n-2}$, $S_\lambda$ is the sharp Poincar\'e--Sobolev constant, $\mathcal{Z}$ is the set of nonzero extremals of the Poincar\'e--Sobolev inequality,
  and the distance is measured with respect to the norm $\|f\|_\lambda^2:=\int_{\mathbb{B}^n}(|\nabla_{\mathbb{B}^n}f|^2-\lambda f^2)\,\mathrm{d}v$.
  We prove the explicit bounds
  \[
    0<\min\Bigl\{\sup_{0<s<1}s\nu(s),\,1-2^{-\frac{2}{n}}\Bigr\}
    \leqslant c_{\mathrm{BE}}
    <\min\Bigl\{1-\frac{2^\ast-1}{\mu_3},\,1-\frac{S_\lambda}{S_n},\,
               2-2^{\frac{n-2}{n}}\Bigr\},
  \]
  where $S_n$ is the sharp Euclidean Sobolev constant, $\mu_3$ is the third eigenvalue of the weighted eigenvalue problem for the linearized operator at the ground state, and $\nu$ is an explicit function of $s$ depending only on $n$ and $\mu_3$.
  Moreover, we prove that $c_{\mathrm{BE}}$ is attained by some $f\notin\mathcal{Z}\cup\{0\}$. This generalizes the Euclidean result of K\"onig \cite{Koenig2025}.

  \medskip
  \noindent\textit{2020 Mathematics Subject Classification.} Primary 46E35, 35A23; Secondary 26D10, 35B33, 35B35, 35J61, 49J40, 58J05.

  \noindent\textit{Key words and phrases.} Poincar\'e--Sobolev inequality, Bianchi--Egnell stability, sharp stability constant, existence of extremals, hyperbolic space, critical exponent.
\end{abstract}

\maketitle
\tableofcontents

\section{Introduction}\label{sec:intro}
The main purpose of the present article is to establish the quantitative stability of critical Poincar\'e--Sobolev inequalities on hyperbolic spaces, and to show that extremal functions exist. To better illustrate our results, we first briefly recall the history of the related works.

We begin with the Euclidean space $\mbRn$ with $n\geqslant 3$. For $1\leqslant p<\infty$, let $L^p(\mbRn)$ be the Lebesgue space with norm $\norm{f}_{L^p(\mbRn)}:=\bigl(\int_{\mbRn}|f|^p\diopx\bigr)^{\frac1p}$. Let $\nabla$ be the Euclidean gradient operator and let $\dot W^{1,p}(\mbRn)$ be the completion of the space $\Ccit(\mbRn)$ of smooth functions with compact support in $\mbRn$ with respect to the norm $\norm{\nabla f}_{L^p(\mbRn)}$, and write $\dot H^1(\mbRn):=\dot W^{1,2}(\mbRn)$. The sharp Sobolev inequality asserts that
\begin{equation}\label{equ:sobolev-Rn}
  \norm{\nabla f}_{L^2(\mbRn)}^2\geqslant S_n\norm{f}_{L^{2^{\ast}}(\mbRn)}^2\quad\text{for every }f\in\dot H^1(\mbRn),
\end{equation}
where $2^{\ast}:=\frac{2n}{n-2}$ and 
\begin{equation*}
  S_n:=\inf_{f\in\dot H^1(\mbRn)\setminus\{0\}}\frac{\norm{\nabla f}_{L^2(\mbRn)}^2}{\norm{f}_{L^{2^{\ast}}(\mbRn)}^2}.
\end{equation*}
Aubin \cite{Aubin1976} and Talenti \cite{Talenti1976} independently obtained the explicit $S_n$ and showed that the nonzero extremals of \eqref{equ:sobolev-Rn} form an $(n+2)$-dimensional manifold
\begin{equation*}
  \mcZ_{\mbRn}:=\Bigl\{c\Bigl(\frac{\lambda}{1+\lambda^2|x-x_0|^2}\Bigr)^{\frac{n-2}{2}}:c\in\mbR\setminus\{0\},\ \lambda>0,\ x_0\in\mbRn\Bigr\}.
\end{equation*}
Brezis and Lieb \cite{BrezisLieb1985} proposed the stability problem of \eqref{equ:sobolev-Rn}: find a suitable distance $\dist(f,\mcZ_{\mbRn})$ which measures how close $f$ is to the Aubin--Talenti bubbles as \eqref{equ:sobolev-Rn} becomes almost an equality.
In the celebrated work \cite{BianchiEgnell1991}, Bianchi and Egnell answered this question affirmatively by showing that, with the \emph{deficit} $\delta_{\mbRn}(f)^2:=\norm{\nabla f}_{L^2(\mbRn)}^2-S_n\norm{f}_{L^{2^{\ast}}(\mbRn)}^2$ of \eqref{equ:sobolev-Rn} and the distance $\dist(f,\mcZ_{\mbRn}):=\inf_{z\in\mcZ_{\mbRn}}\norm{\nabla(f-z)}_{L^2(\mbRn)}$, the constant
\begin{equation}\label{equ:cbe-Rn}
  c_{\mathrm{BE}}(\mbRn):=\inf_{f\in\dot H^1(\mbRn)\setminus(\mcZ_{\mbRn}\cup\{0\})}\frac{\delta_{\mbRn}(f)^2}{\dist(f,\mcZ_{\mbRn})^2}
\end{equation}
is positive, that is,
\begin{equation}\label{equ:BE-Rn}
  \delta_{\mbRn}(f)^2
  \geqslant c_{\mathrm{BE}}(\mbRn)\dist(f,\mcZ_{\mbRn})^2\quad\text{for every }f\in\dot H^1(\mbRn).
\end{equation}

After \cite{BianchiEgnell1991} there is a vast literature on stability estimates of this type for related inequalities; see \cite{Frank2024} for a survey. 
One direction is to extend \eqref{equ:BE-Rn} from $p=2$ to general $1<p<n$.
Precisely, consider the sharp $p$-Sobolev inequality 
\begin{equation*}
  \norm{\nabla f}_{L^p(\mbRn)}^2\geqslant S_{n,p}\norm{f}_{L^{\frac{np}{n-p}}(\mbRn)}^2 \quad \text{for every\ } f\in \dot W^{1,p}(\mbRn)
\end{equation*}
with sharp constant $S_{n,p}$, so that $S_{n,2}=S_n$. For $p\neq 2$ the Hilbertian structure used by Bianchi and Egnell is unavailable. Figalli and Zhang \cite{FigalliZhang2022} proved that the deficit, normalized by $\norm{\nabla f}_{L^p(\mbRn)}^2$, controls the $\max\{2,p\}$-th power of the relative distance to the extremals in $\dot W^{1,p}(\mbRn)$. We refer to \cite{CFMP2009,FigalliNeumayer2019,Neumayer2020} for earlier work. Fractional and weighted stability inequalities have also received much attention. Chen, Frank and Weth \cite{CFW2013} extended \eqref{equ:BE-Rn} to the sharp fractional Sobolev inequality in $\dot H^{\alpha}(\mbRn)$ for every $0<\alpha<\frac{n}{2}$, while Wei and Wu \cite{WeiWu2022} proved the analogue of \eqref{equ:BE-Rn} for the Caffarelli--Kohn--Nirenberg inequality in the parameter region where the extremal function is unique up to dilations and multiplication by constants.

Since the argument of Bianchi and Egnell is by contradiction and compactness, it gives no quantitative information on $c_{\mathrm{BE}}(\mbRn)$, and it is natural to ask for explicit bounds. More recently, Dolbeault, Esteban, Figalli, Frank and Loss \cite{DEFFL2025} proved that there exists an explicit constant $\beta>0$, independent of $n$, such that $c_{\mathrm{BE}}(\mbRn)\geqslant\frac{\beta}{n}$ for every $n\geqslant3$. They also deduced from this bound a stability estimate for the Gaussian logarithmic Sobolev inequality with an explicit dimension-free constant. 
Chen, Lu and Tang \cite{CLTa2024,CLTc2025,CLTb2026,CLTW2026} treated the mutually dual Hardy--Littlewood--Sobolev and Sobolev inequalities of order $\alpha$. They obtained explicit lower bounds for all $0<\alpha<\frac{n}{2}$ in \cite{CLTa2024}, and lower bounds with optimal dependence on the dimension for $0<\alpha<1$ in \cite{CLTc2025} and for $1\leqslant\alpha<\frac{n}{2}$ in \cite{CLTb2026}. They also extended these results from Euclidean space to the Heisenberg group in \cite{CLTW2026} and obtained there the optimal asymptotic lower bound for the Sobolev stability constant, with a CR Yamabe flow in place of the unavailable rearrangement techniques.

As for quantitative upper bounds, the quotient in \eqref{equ:cbe-Rn} tends to $\frac{4}{n+4}$ along perturbations of an element of $\mcZ_{\mbRn}$ in the direction of a spherical harmonic of degree $2$ \cite{CFW2013}. Together with the lower bound of \cite{DEFFL2025}, this shows that the order $\frac{1}{n}$ is optimal. The quotient also tends to $2-2^{\frac{n-2}{n}}$ along two bubbles, that is, two elements in $\mcZ_{\mbRn}$, whose scales separate. K\"onig \cite{Koenig2023a,Koenig2025} proved that both upper bounds are strict, that is,
\begin{equation*}
  c_{\mathrm{BE}}(\mbRn)<\min\Bigl\{\frac{4}{n+4},2-2^{\frac{n-2}{n}}\Bigr\}\quad\text{for every }n\geqslant3,
\end{equation*}
and deduced from these strict bounds that every minimizing sequence for \eqref{equ:cbe-Rn}, normalized in $\dot H^1(\mbRn)$, converges strongly up to a subsequence. Consequently, $c_{\mathrm{BE}}(\mbRn)$ is attained by some $f\in\dot H^1(\mbRn)\setminus(\mcZ_{\mbRn}\cup\{0\})$.

 Another important line of research concerns the Euler--Lagrange equation $-\Delta u=|u|^{q-1}u$ of \eqref{equ:sobolev-Rn}, where $q:=2^{\ast}-1=\frac{n+2}{n-2}$. Struwe \cite{Struwe1984} proved that, on a bounded domain, every Palais--Smale sequence of the Brezis--Nirenberg functional splits, up to a subsequence, into a solution of the Dirichlet problem and finitely many rescaled solutions of this equation on $\mbRn$. On $\mbRn$ this gives the following qualitative picture \cite{FigalliGlaudo2020}. If a sequence $\{u_k\}$ of nonnegative functions with bounded energy satisfies $\norm{\Delta u_k+u_k^{q}}_{H^{-1}(\mbRn)}\to0$, where $H^{-1}(\mbRn)$ is the dual of $\dot H^1(\mbRn)$, then, up to a subsequence, $u_k$ is close in $\dot H^1(\mbRn)$ to a sum of weakly interacting bubbles. Ciraolo, Figalli and Maggi \cite{CFM2018} obtained the first quantitative version of this picture, in the case of one bubble. As an application, they obtained a quantitative rate of convergence for a fast diffusion equation related to the Yamabe flow. For sums of several bubbles, Figalli and Glaudo \cite{FigalliGlaudo2020} discovered a striking dependence on the dimension. With $\Gamma(u):=\norm{\Delta u+|u|^{q-1}u}_{H^{-1}(\mbRn)}$, they proved that the distance from $u$ to sums of bubbles is bounded by a constant multiple of $\Gamma(u)$ if $3\leqslant n\leqslant5$, and that this linear bound fails for $n\geqslant6$. More recently, Deng, Sun and Wei \cite{DSW2025} settled the remaining dimensions.  They proved that, near a sum of weakly interacting bubbles, the distance is bounded by a constant multiple of $\Gamma(u)|\log\Gamma(u)|^{\frac12}$ if $n=6$ and of $\Gamma(u)^{\frac{n+2}{2(n-2)}}$ if $n\geqslant7$. Moreover, these bounds are sharp, so that, together with the result of Figalli and Glaudo, the optimal dependence on the dimension is now known in every dimension $n\geqslant3$.

Quantitative stability of Euler--Lagrange equations has also been studied in other settings. Wei and Wu \cite{WeiWu2022} obtained a stability estimate for the Euler--Lagrange equation of the Caffarelli--Kohn--Nirenberg inequality. Chen, Kim and Wei \cite{CKW2026} obtained sharp quantitative stability estimates for the Brezis--Nirenberg problem $-\Delta u-\lambda u=u^{q}$ in a smooth bounded domain of $\mbRn$. 

 The results recalled so far concern Euclidean space. Stability of Sobolev-type inequalities has also been studied on Riemannian manifolds, where new phenomena appear. On a closed Riemannian manifold, Engelstein, Neumayer and Spolaor \cite{ENS2022} proved that, if a function nearly minimizes the Yamabe energy, then the corresponding conformal metric is quantitatively close to a minimizing Yamabe metric in the conformal class. Generically this distance is controlled quadratically by the Yamabe energy deficit, but they also gave an example for which the quadratic estimate fails. Motivated by their work, Frank \cite{Frank2022} showed that on $\mathbb S^1\bigl(\frac{1}{\sqrt{n-2}}\bigr)\times\mathbb S^{n-1}(1)$ the conformally invariant Sobolev inequality holds with a remainder term given by the fourth power of the distance to the optimizers, and that the fourth power is optimal. Thus, in contrast with \eqref{equ:BE-Rn}, on a manifold the square of the distance is not always the right power. On Riemannian manifolds $(M,g_M)$ with Ricci curvature bounded below, Nobili and Violo \cite{NobiliViolo2024} studied the qualitative stability of the sharp Sobolev inequality. If $\mathrm{Ric}_{g_M}\geqslant(n-1)g_M$, they proved that an almost extremal function is close to an extremal function of the round sphere, and if the Ricci curvature is nonnegative and the volume growth is Euclidean, that it is close to an extremal function of the Euclidean Sobolev inequality. We refer to the survey \cite{Nobili2026} for further results in this direction.

 Hyperbolic space is a complete, simply connected Riemannian manifold of constant sectional curvature $-1$. Being noncompact with negative Ricci curvature, it is not covered by the results on manifolds recalled above. Let $\mbBn$ with $n\geqslant3$ be hyperbolic $n$-space in the Poincar\'e ball model, with metric $g:=\rho^2|\mrd x|^2$, where $\rho(x):=\frac{2}{1-|x|^2}$, volume element $\diopv:=\rho^n\diopx$ and hyperbolic gradient $\nabla_{\mbBn}$.  For $1\leqslant p\leqslant\infty$, let $L^p(\mbBn)$ be the Lebesgue space with respect to $\diopv$, and let $H^1(\mbBn)$ be the Sobolev space with respect to $\diopv$, that is, the completion of the space $\Ccit(\mbBn)$ of smooth functions with compact support in $\mbBn$ with respect to the norm $\norm{f}_{H^1(\mbBn)}:=\bigl(\int_{\mbBn}(|\nabla_{\mbBn}f|^2+f^2)\diopv\bigr)^{\frac12}$. By \cite{McKean1970}, the bottom of the $L^2$ spectrum of the Laplace--Beltrami operator $-\Delta_{\mbBn}$ is $\lambda_1:=\frac{(n-1)^2}{4}$. Hence, for $\lambda<\lambda_1$, the quadratic form
\begin{equation*}
  \norm{f}_\lambda^2:=\int_{\mbBn}\bigl(|\nabla_{\mbBn}f|^2-\lambda f^2\bigr)\diopv
\end{equation*}
is positive definite, and $\norm{\cdot}_\lambda$ is a norm on $H^1(\mbBn)$ equivalent to $\norm{\cdot}_{H^1(\mbBn)}$, with inner product
\begin{equation*}
  \lag f,h\rag_\lambda:=\int_{\mbBn}(\nabla_{\mbBn}f\cdot\nabla_{\mbBn}h-\lambda fh)\diopv.
\end{equation*}
The \emph{critical Poincar\'e--Sobolev inequality} asserts that, for every $\lambda<\lambda_1$,
\begin{equation}
  \norm{f}_\lambda^2\geqslant S_\lambda\norm{f}_{L^{2^{\ast}}(\mbBn)}^2\quad\text{for every }f\in H^1(\mbBn),\label{equ:PS}
\end{equation}
where
\begin{equation}
  S_\lambda:=\inf_{f\in H^1(\mbBn)\setminus\{0\}}
  \frac{\norm{f}_\lambda^2}{\norm{f}_{L^{2^{\ast}}(\mbBn)}^2}>0.\label{equ:Slambda}
\end{equation}
When $\lambda=\lambda_0:=\frac{n(n-2)}{4}$, the operator $-\Delta_{\mbBn}-\lambda_0$ is the conformal Laplacian, and \eqref{equ:PS} is equivalent to the Sobolev inequality \eqref{equ:sobolev-Rn}; in particular, $S_{\lambda_0}=S_n$.

 From now on, we assume that $n$ and $\lambda$ satisfy
\begin{equation*}
  n\geqslant 4\quad\text{and\quad} \lambda_0<\lambda<\lambda_1.
\end{equation*}
Among other results, Mancini and Sandeep \cite{ManciniSandeep2008} showed that $S_\lambda<S_n$ (see \cite[proof of Theorem~5.2, p.~662]{ManciniSandeep2008}) and that the infimum in \eqref{equ:Slambda} is attained by a positive function $\mcU\in H^1(\mbBn)$, radial about the origin $o$ of $\mbBn$, which satisfies
\begin{equation}\label{equ:EL}
  -\Delta_{\mbBn}\mcU-\lambda\mcU=\mcU^{q},
\end{equation}
where we recall that $q=2^\ast-1=\frac{n+2}{n-2}$. 
For $\xi\in\mbBn$, let $\mcU_\xi:=\mcU\circ\gamma_\xi^{-1}$, where $\gamma_\xi$ is the hyperbolic isometry sending $o$ to $\xi$. They also proved that the $(n+1)$-dimensional cone
\begin{equation}
  \mcZ=\mcZ_{\mbBn}:=\{c\mcU_\xi:c\in\mbR\setminus\{0\},\ \xi\in\mbBn\}\label{equ:Zo}
\end{equation}
is exactly the set of nonzero extremals of \eqref{equ:PS}, and that the positive finite-energy solutions of \eqref{equ:EL} are exactly the hyperbolic translates $\mcU_\xi$ of $\mcU$. Bhakta and Sandeep \cite{BhaktaSandeep2012} showed that Palais--Smale sequences of the associated energy functional decompose, up to a vanishing remainder, into hyperbolic bubbles and localized Aubin--Talenti bubbles, and Ganguly and Sandeep \cite{GangulySandeep2015} proved the nondegeneracy of $\mcU$.

The stability problem for \eqref{equ:PS} can be formulated as follows: define the \emph{deficit} of \eqref{equ:PS} by
\begin{equation*}
  \delta(f)^2:=\norm{f}_\lambda^2-S_\lambda\norm{f}_{L^{2^{\ast}}(\mbBn)}^2\geqslant0,
\end{equation*}
and measure the distance from $f$ to $\mcZ$ by
\begin{equation}\label{equ:deficit-sigma}
  \dist(f,\mcZ):=\inf_{z\in\mcZ}\norm{f-z}_\lambda,\qquad
  \sigma(f):=\frac{\dist(f,\mcZ)^2}{\norm{f}_\lambda^2}\in[0,1],
\end{equation}
where $\sigma(f)$ is called the \emph{relative distance} of $f$ to $\mcZ$. By \eqref{equ:Zo}, every nonzero $f\in H^1(\mbBn)$ with $\delta(f)=0$ lies in $\mcZ$. For $f\in H^1(\mbBn)\setminus(\mcZ\cup\{0\})$, we call $\frac{\delta(f)^2}{\dist(f,\mcZ)^2}$ the \emph{stability quotient} of $f$, and define the \emph{stability constant} of Bianchi--Egnell type as 
\begin{equation}
  c_{\mathrm{BE}}=c_{\mathrm{BE}}(\mbBn,\lambda):=
  \inf_{f\in H^1(\mbBn)\setminus(\mcZ\cup\{0\})}
  \frac{\delta(f)^2}{\dist(f,\mcZ)^2}\geqslant0.\label{equ:cbe-def}
\end{equation}

The first qualitative stability result for \eqref{equ:PS} is due to Bhakta, Ganguly, Karmakar and Mazumdar \cite{BGKM2025}. They proved that there exists a constant $C>0$ such that
\begin{equation*}
  \dist(f,\mcZ)\leqslant C\,\delta(f)\quad\text{for every }f\in H^1(\mbBn).
\end{equation*}
 They also proved a quantitative stability result for \eqref{equ:EL}, which is the hyperbolic counterpart of the one-bubble estimate of Ciraolo, Figalli and Maggi \cite{CFM2018}. As an application, they studied the extinction rate of the subcritical fast diffusion flow for radial data. In \cite{BGKM2025b}, the same authors proved the hyperbolic counterpart of the result of Figalli and Glaudo \cite{FigalliGlaudo2020}. For $3\leqslant n\leqslant5$ and power nonlinearities with exponent in $(2,q]$, under a further condition on $\lambda$, they controlled the $H^1(\mbBn)$-distance to sums of bubbles by the $H^{-1}(\mbBn)$-norm of the residual of the corresponding equation.

In this article, we study the quantification of $c_{\mathrm{BE}}$ and the existence of extremal functions. Our first result concerns functions near $\mcZ$. For such functions, the stability quotient is governed by the weighted eigenvalue problem obtained by linearizing \eqref{equ:EL} at $\mcU$,
\begin{equation}\label{equ:const-linearized}
  (-\Delta_{\mbBn}-\lambda)\Psi=\mu\mcU^{q-1}\Psi,\quad \Psi\in H^1(\mbBn).
\end{equation}
From \cite{BGKM2025,GangulySandeep2015}, the eigenvalues of \eqref{equ:const-linearized} form a discrete sequence tending to infinity, and the first two are $\mu_1=1$ and $\mu_2=q$ (see Section~\ref{subsec:linearized-spectrum} for more details). Let $\mu_3>q$ be the third eigenvalue, counted without multiplicity, and for $0<s<1$ let
\begin{equation*}
  \mcT_s:=\{f\in H^1(\mbBn)\setminus\{0\}:0<\sigma(f)\leqslant s\}
\end{equation*}
be the \emph{tubes} around $\mcZ$. Our first result reads as follows.
\begin{theorem}\label{thm:crit-local}
  The \emph{local stability constant} is
  \begin{equation*}
    c_{\mrloc}:=\lim_{s\downarrow0}
    \inf_{f\in\mcT_s}\frac{\delta(f)^2}{\dist(f,\mcZ)^2}
    =1-\frac{q}{\mu_3}.
  \end{equation*}
\end{theorem}

\begin{remark}
    A local stability estimate near $\mcZ$ was obtained in \cite{BGKM2025}. Theorem~\ref{thm:crit-local} identifies the sharp local constant, which plays a key role in the proof of Theorem~\ref{thm:crit-global-bounds}.
\end{remark}

As we will show in Section~\ref{sec:attainment}, a minimizing sequence $\{f_k\}$ for \eqref{equ:cbe-def}, normalized by $\norm{f_k}_{L^{2^{\ast}}(\mbBn)}=1$, can fail to converge, modulo isometries, to an element of $H^1(\mbBn)\setminus(\mcZ\cup\{0\})$ in only three ways: it can approach $\mcZ$, concentrate at a point, or split into two pieces that move apart. We call these the \emph{local channel}, the \emph{concentration channel} and the \emph{dichotomy channel}, respectively. Along each channel the stability quotient has an explicit limit, which we compute on a suitable family of test functions. In the local channel we take functions approaching $\mcZ$ along an eigenfunction associated with $\mu_3$, and the limit is $c_{\mrloc}$ by Theorem~\ref{thm:crit-local}. In the concentration channel we take Euclidean bubbles concentrating at a point. Their energy $\norm{f}_\lambda^2$ tends to $S_n$ and their relative distance $\sigma$ to $\mcZ$ tends to $1$, so the limit is $c_{\mathrm{conc}}:=1-\frac{S_\lambda}{S_n}$. In the dichotomy channel we take two ground states whose centers move apart. Their energy tends to $2^{\frac{2}{n}}S_\lambda$, and $\sigma$ tends to $\frac{1}{2}$ since only one of the two pieces is at positive distance from $\mcZ$; hence the limit is $c_{\mathrm{dich}}:=2-2^{\frac{n-2}{n}}$. Let
\begin{equation}\label{equ:const-nu-explicit}
  \nu(s):=\begin{cases}
    1-\frac{q}{\mu_3}-2\sqrt{\frac{s}{1-s}}-\frac{s}{2(1-s)},&n=4,\\
    1-\frac{q}{\mu_3}-\frac{28}{27}\bigl(\sqrt{\frac{s}{1-s}}+\bigl(\frac{s}{1-s}\bigr)^{\frac{2}{3}}\bigr),&n=5,\\
    1-\frac{q}{\mu_3}-\frac{2}{2^{\ast}}\bigl(\frac{s}{1-s}\bigr)^{\frac{q-1}{2}},&n\geqslant6,
  \end{cases}.
\end{equation}
Our second result is a global quantitative stability estimate:
\begin{theorem}\label{thm:crit-global-bounds}
  The stability constant $c_{\mathrm{BE}}$ defined in \eqref{equ:cbe-def} satisfies
  \begin{equation}\label{equ:intro-global-two-sided}
    0 < \min\Bigl\{\sup_{0<s<1}s\nu(s),\frac{1}{2}c_{\mathrm{dich}}\Bigr\}
    \leqslant c_{\mathrm{BE}}
    <\min\bigl\{c_{\mrloc},c_{\mathrm{conc}},c_{\mathrm{dich}}\bigr\}.
  \end{equation}
\end{theorem}

Since the three upper bounds in \eqref{equ:intro-global-two-sided} are strict, a normalized minimizing sequence cannot escape through any of the three channels. This gives our last main result, a hyperbolic version of K\"onig's result \cite{Koenig2025}.
\begin{theorem}\label{thm:attainment}
  There exists $f\in H^1(\mbBn)\setminus(\mcZ\cup\{0\})$ such that $\delta(f)^2=c_{\mathrm{BE}}\dist(f,\mcZ)^2$.
\end{theorem}

Hyperbolic geometry changes both the compactness properties and the local analysis of the stability problem. The absence of dilation symmetry makes point concentration an additional obstruction to attainment, since it cannot be removed by hyperbolic isometries. Moreover, a conformal change of variables leaves a potential term, so the explicit spectral calculation on the sphere cannot be used directly. For the strict local upper bound, we prove that the third-order term in the expansion of the deficit is nonzero. The main difficulty is the radial direction, which corresponds to dilations in the Euclidean case, where this term vanishes; here its nonvanishing is a purely hyperbolic effect, and we prove it by a radial comparison argument. The global argument combines the local estimate with a constrained energy flow and compactness below the concentration and dichotomy thresholds.

The paper is organized as follows. Section~\ref{sec:prelim} collects preliminary facts on the geometry and analysis of hyperbolic space, and Section~\ref{sec:crit-local} establishes explicit local stability. Section~\ref{sec:crit-global-bounds-proof} derives a preliminary global lower bound using the constrained flow. Section~\ref{sec:strict-upper} proves the strict upper bounds and completes the proof of Theorem~\ref{thm:crit-global-bounds}. Finally, Section~\ref{sec:attainment} proves the existence of extremals.

\section{Preliminaries}\label{sec:prelim}
Throughout the article, all function spaces consist of real-valued functions, and $|x|$ and $x\cdot y$ denote the Euclidean norm and inner product. For a Euclidean domain $\Omega$, $H^1(\Omega)$, $H_0^1(\Omega)$ and $H^{-1}(\Omega)$ are the usual Sobolev spaces with respect to $\diopx$. For $1<p<\infty$, $p':=\frac{p}{p-1}$ denotes the conjugate exponent of $p$, so that $\frac{1}{p}+\frac{1}{p'}=1$; in particular, $(2^{\ast})'=\frac{2n}{n+2}$. For positive functions $f$ and $g$ of one variable, we write $f\sim g$ as the variable tends to an endpoint if there are constants $0<c\leqslant C$ such that $cg\leqslant f\leqslant Cg$ near that endpoint. Inside proofs, an unnamed constant $C$ may change from line to line; constants named in a statement are fixed.

\subsection{The geometry and analysis of hyperbolic space}\label{subsec:geometry-analysis}
We first briefly recall the geometry of hyperbolic space used in this article. As a set, the Poincar\'e ball model is $\mbBn=\{x\in\mbRn:|x|<1\}$, with the metric $g$ of Section~\ref{sec:intro}. The distance between $x,y\in\mbBn$ is
\begin{equation*}
  d(x,y):=\operatorname{arcosh}\Bigl(1+\frac{2|x-y|^2}{(1-|x|^2)(1-|y|^2)}\Bigr),
\end{equation*}
where $\operatorname{arcosh}$ is the inverse function of the hyperbolic cosine function. 
In particular, the distance to the origin $o:=0\in \mbBn$ is $d(x,o)=\log\frac{1+|x|}{1-|x|}$, and we write $B_{\mbH}(x,r)$ for the geodesic ball with center $x$ and radius $r$. The isometries of $\mbBn$ are exactly the M\"obius transformations of $\mbRn\cup\{\infty\}$ that map $\mbBn$ onto itself. 
In particular, for $\xi\in\mbBn$, the hyperbolic translation $\gamma_\xi$ introduced in Section~\ref{sec:intro}, which is also the M\"obius translation, is given by
\begin{equation*}
  \gamma_\xi(x)=\frac{(1-|\xi|^2)x+(1+2\xi\cdot x+|x|^2)\xi}{1+2\xi\cdot x+|\xi|^2|x|^2}.
\end{equation*}
Moreover, every isometry of $\mbBn$ can be expressed as $\varphi:=\gamma_\xi\circ A$ for some $\xi\in\mbBn$ and $A\in O(n)$, where $O(n)$ is the orthogonal group of $\mbRn$, which acts on $\mbBn$ as the group of isometries fixing $o$. This fact, together with the fact that $\mcU$ is radial, yields that the cone \eqref{equ:Zo} can be rewritten as
\begin{equation*}
  \mcZ=\{c\,\mcU\circ\varphi^{-1}:\ c\in\mbR\setminus\{0\},\ \varphi\text{ is an isometry of }\mbBn\}.
\end{equation*}
We also write $B_r\subset\mbRn$ for the Euclidean ball of radius $r$ centered at $0$; thus $B_1$ is the same set as $\mbBn$, regarded with the Euclidean metric.

Besides geometry, we introduce the analysis of hyperbolic space. Here and below, quantities without the subscript $\mbBn$, such as the gradient $\nabla$ and the Laplacian $\Delta$, are Euclidean. The hyperbolic gradient and the Laplace--Beltrami operator are expressed through the Euclidean $\nabla$ and $\Delta$ by
\begin{equation*}
  \nabla_{\mbBn}=\rho^{-2}\nabla,\qquad
  \Delta_{\mbBn}=\rho^{-n}\MaOpdiv(\rho^{n-2}\nabla)=\rho^{-2}\bigl(\Delta+(n-2)\rho^{-1}\nabla\rho\cdot\nabla\bigr).
\end{equation*}
Writing $|\cdot|_g$ and $\lag\cdot,\cdot\rag_g$ for the length and the inner product of tangent vectors with respect to $g$, we have
\begin{equation*}
  |\nabla_{\mbBn}f|_g=\rho^{-1}|\nabla f|,\qquad
  \lag\nabla_{\mbBn}f,\nabla_{\mbBn}h\rag_g=\rho^{-2}\nabla f\cdot\nabla h,
\end{equation*}
and in particular, using $\diopv=\rho^n\diopx$,
\begin{equation*}
  \int_{\mbBn}|\nabla_{\mbBn}f|_g^2\diopv=\int_{\mbBn}\rho^{n-2}|\nabla f|^2\diopx.
\end{equation*}
For simplicity, when no confusion can arise, we drop the subscript $g$ and write $|\nabla_{\mbBn}f|$ and $\nabla_{\mbBn}f\cdot\nabla_{\mbBn}h$. In geodesic polar coordinates $(r,\theta)$ centered at $o$, where $r=d(x,o)$ and $\theta\in\mbS^{n-1}$, we have
\begin{align}
  &g=\diopr^2+\sinh^2r\,g_{\mbS^{n-1}},\nonumber\\*
  &\diopv=\sinh^{n-1}r\diopr\dioptheta,\label{equ:polar-volume}\\*
  &|\nabla_{\mbBn}f|^2=(\partial_rf)^2+\sinh^{-2}r\,|\nabla_{\mbS^{n-1}}f|^2,\label{equ:polar-gradient}\\*
  &\Delta_{\mbBn}f=\partial_r^2f+(n-1)\coth r\,\partial_rf+\sinh^{-2}r\,\Delta_{\mbS^{n-1}}f,\label{equ:polar-laplacian}
\end{align}
where $g_{\mbS^{n-1}}$, $\dioptheta$, $\nabla_{\mbS^{n-1}}$ and $\Delta_{\mbS^{n-1}}$ are the round metric, the surface measure, the gradient and the Laplace--Beltrami operator of the unit sphere $\mbS^{n-1}$, respectively. In particular, 
\begin{equation*}
  \Delta_{\mbBn}f=f''+(n-1)\coth r\,f'
\end{equation*}
for every radial function $f=f(r)$. Finally, since the isometries of $\mbBn$ preserve $g$, they preserve $\diopv$ and $|\nabla_{\mbBn}f|$, and commute with $\Delta_{\mbBn}$. Consequently, for every isometry $\varphi$, the map $f\mapsto f\circ\varphi$ preserves $\norm{\cdot}_\lambda$ and $\norm{\cdot}_{L^{2^{\ast}}(\mbBn)}$, and hence, by the description of $\mcZ$ above, also $\delta$, $\dist(\cdot,\mcZ)$ and $\sigma$.

\subsection{The ground state}\label{subsec:ground}
\begin{lemma}\label{lem:ground-state}
  The ground state $\mcU$ is strictly
  decreasing in $d(x,o)$, $\mcU\in C^2(\mbBn)\cap L^\infty(\mbBn)$, and for every $\ve>0$
  \begin{equation}\label{equ:ground-state-tail}
    \mcU(x)\leqslant C_\ve\mre^{-(\beta_+-\ve)d(x,o)},\qquad\beta_+:=\frac{n-1}{2}+\sqrt{\lambda_1-\lambda}>\frac{n-1}{2},
  \end{equation}
  for a constant $C_\ve$ depending only on $n$, $\lambda$ and $\ve$.
\end{lemma}
\begin{proof}
  The ground state is strictly decreasing in the radial variable \cite[Lemma~3.2]{ManciniSandeep2008}. Local integrability and elliptic regularity for critical potentials \cite{BrezisKato1979} (see also \cite[p.~257]{BhaktaSandeep2012}) give $\mcU\in C^2(\mbBn)\cap L^\infty(\mbBn)$. By \cite[Lemma~3.4]{ManciniSandeep2008}, $\frac{1}{d(x,o)}\log\mcU(x)\to-\beta_+$ as $d(x,o)\to\infty$. Combining this limit with the boundedness of $\mcU$ gives \eqref{equ:ground-state-tail} for every $\ve>0$. This finishes the proof of Lemma~\ref{lem:ground-state}.
\end{proof}
Here $\beta_+$ is the larger root of $\beta(n-1-\beta)=\lambda$, the indicial equation of $\Delta_{\mbBn}+\lambda$ at infinity; the smaller root $\beta_-:=\frac{n-1}{2}-\sqrt{\lambda_1-\lambda}$ belongs to the slowly decaying radial solutions, which satisfy $2\beta_-<n-1$ and are not square integrable. Taking $0<\ve<\sqrt{\lambda_1-\lambda}$ in \eqref{equ:ground-state-tail} gives $2(\beta_+-\ve)>n-1$, hence $\mcU\in L^2(\mbBn)$.

The following lemma shows that the orbit map $\xi\mapsto\mcU_\xi$ is of class $C^2$ into $H^1(\mbBn)$.
\begin{lemma}\label{lem:pairing-center-derivative}
  Let $r=d(x,o)$ and write $\mcU(r)$ for the radial profile of the ground state.
  \begin{enumerate}[label=\textnormal{(\alph*)}]
    \item\label{item:ground-state-smooth} $\mcU\in C^\infty(\mbBn)$, $\mcU'(0)=0$, and for $k=0,1,2,3$ and every $\ve>0$,
    \begin{equation}\label{equ:ground-state-higher-tail}
      |\mcU^{(k)}(r)|\leqslant C\mre^{-(\beta_+-\ve)r}\qquad(r\geqslant0),
    \end{equation}
    with $C$ depending only on $n$, $\lambda$, $k$ and $\ve$.
    \item\label{item:orbit-map-C2} The orbit map $\xi\mapsto\mcU_\xi$ is of class $C^2$ from $\mbBn$ into $H^1(\mbBn)$; its derivatives up to order two are the pointwise derivatives of $(x,\xi)\mapsto\mcU(d(x,\xi))$.
    \item\label{item:pairing-center-partial} For every $f\in H^1(\mbBn)$ the function $\xi\mapsto\lag f,\mcU_\xi\rag_\lambda$ is of class $C^2$ on $\mbBn$, and $\partial_{\xi_i}\lag f,\mcU_\xi\rag_\lambda=\lag f,\partial_{\xi_i}\mcU_\xi\rag_\lambda$.
  \end{enumerate}
\end{lemma}
\begin{proof}
  Throughout, $0<\ve<\sqrt{\lambda_1-\lambda}$; larger values of $\ve$ only weaken the bounds. Since $\lambda_1-\lambda<\lambda_1-\lambda_0=\frac{1}{4}$, we have
  \begin{equation*}
    \frac{n-1}{2}<\beta_+-\ve<\beta_+<\frac{n}{2}\leqslant n-1.
  \end{equation*}
  We first prove \ref{item:ground-state-smooth}, then use it to prove \ref{item:orbit-map-C2}, and finally deduce \ref{item:pairing-center-partial} from \ref{item:orbit-map-C2}.

  \emph{Proof of \ref{item:ground-state-smooth}.} We first show that $\mcU\in C^\infty(\mbBn)$. By the formula for $\Delta_{\mbBn}$ above, \eqref{equ:EL} reads in Euclidean coordinates
  \begin{equation*}
    -\Delta\mcU-(n-2)\rho^{-1}\nabla\rho\cdot\nabla\mcU=\rho^2\bigl(\lambda\mcU+\mcU^q\bigr)\qquad\text{in }\mbBn .
  \end{equation*}
  The coefficients on the left-hand side are smooth in $\mbBn$, and, since $\mcU>0$ and $t\mapsto\lambda t+t^q$ is smooth on $(0,\infty)$, the right-hand side is a smooth function of $(x,\mcU(x))$. Hence, if $\mcU\in C^{k,\alpha}_{\mrloc}(\mbBn)$ for some $k\geqslant1$ and $0<\alpha<1$, then the right-hand side lies in $C^{k,\alpha}_{\mrloc}(\mbBn)$, and interior Schauder estimates \cite[Theorem~6.17]{GilbargTrudinger2001} give $\mcU\in C^{k+2,\alpha}_{\mrloc}(\mbBn)$. Starting from $\mcU\in C^2(\mbBn)$ (Lemma~\ref{lem:ground-state}) and iterating, we obtain $\mcU\in C^\infty(\mbBn)$.

  Next we pass to the profile. Since $d(o,\tanh(\frac{r}{2})e_1)=|r|$ for $r\in\mbR$ and $\mcU$ is radial, the profile $r\mapsto\mcU(r)=\mcU(\tanh(\frac{r}{2})e_1)$ extends to a smooth even function on $\mbR$. Hence $\mcU'(0)=0$ and every $\mcU^{(k)}$ is bounded on $[0,1]$. Since $\mre^{-(\beta_+-\ve)r}\geqslant\mre^{-\beta_+}$ on $[0,1]$, it suffices to prove \eqref{equ:ground-state-higher-tail} for $r\geqslant1$. The case $k=0$ is \eqref{equ:ground-state-tail}.

  For $k=1$, by \eqref{equ:polar-volume} and \eqref{equ:polar-laplacian}, in geodesic polar coordinates about $o$ we have $\diopv=\sinh^{n-1}r\diopr\dioptheta$ and $\Delta_{\mbBn}\mcU=\mcU''+(n-1)\coth r\,\mcU'$ for radial $\mcU$. So, for $r>0$, \eqref{equ:EL} reads
  \begin{equation*}
    \bigl(\sinh^{n-1}r\,\mcU'(r)\bigr)'=-\sinh^{n-1}r\,\bigl(\lambda\mcU(r)+\mcU(r)^q\bigr).
  \end{equation*}
  Integrating over $[0,r]$ and using $\mcU'(0)=0$, we find
  \begin{align*}
    \bigl|\sinh^{n-1}r\,\mcU'(r)\bigr|
    &=\int_0^r\sinh^{n-1}s\,\bigl(\lambda\mcU+\mcU^q\bigr)\diop s\\
    &\leqslant\bigl(\lambda+\norm{\mcU}_{L^\infty(\mbBn)}^{q-1}\bigr)C_\ve\int_0^r\mre^{(n-1-\beta_++\ve)s}\diop s
    \leqslant C\mre^{(n-1-\beta_++\ve)r},
  \end{align*}
  where in the first inequality we have used $\sinh s\leqslant\mre^{s}$, $\mcU^q\leqslant\norm{\mcU}_{L^\infty(\mbBn)}^{q-1}\mcU$ and \eqref{equ:ground-state-tail}, and in the last inequality we have used $n-1-\beta_++\ve>0$, which follows from the first display of this proof. Since $\sinh r\geqslant\frac{1-\mre^{-2}}{2}\mre^{r}$ for $r\geqslant1$, dividing by $\sinh^{n-1}r$ gives
  \begin{equation*}
    |\mcU'(r)|\leqslant C\mre^{-(\beta_+-\ve)r}\qquad(r\geqslant1).
  \end{equation*}

  For $k=2,3$, the radial equation and its derivative give, for $r\geqslant1$,
  \begin{align*}
    \mcU''&=-(n-1)\coth r\,\mcU'-\lambda\mcU-\mcU^q,\\
    \mcU'''&=(n-1)\sinh^{-2}r\,\mcU'-(n-1)\coth r\,\mcU''-\lambda\mcU'-q\mcU^{q-1}\mcU'.
  \end{align*}
  Since $\coth r\leqslant\coth1$, $\sinh^{-2}r\leqslant\sinh^{-2}1$ and $\mcU^{q-1}\leqslant\norm{\mcU}_{L^\infty(\mbBn)}^{q-1}$, we obtain $|\mcU''|\leqslant C(|\mcU'|+\mcU)$ and $|\mcU'''|\leqslant C(|\mcU'|+|\mcU''|)$ on $[1,\infty)$. Hence the cases $k=2,3$ of \eqref{equ:ground-state-higher-tail} follow from the cases $k=0,1$.

  \emph{Proof of \ref{item:orbit-map-C2}.} Choose $\gamma_\xi$ to be the M\"obius translation of $\mbBn$ taking $o$ to $\xi$. Its inverse is real-analytic in $(x,\xi)\in\mbBn\times\mbBn$, and, by \ref{item:ground-state-smooth}, the function
  \begin{equation*}
    F(x,\xi):=\mcU_\xi(x)=\mcU(\gamma_\xi^{-1}x)
  \end{equation*}
  is $C^\infty$ on $\mbBn\times\mbBn$. Since $\gamma_\xi^{-1}$ is an isometry with $\gamma_\xi^{-1}(\xi)=o$ and $\mcU$ is radial, we also have $F(x,\xi)=\mcU(d(x,\xi))$ in terms of the profile. We prove \ref{item:orbit-map-C2} in three steps: a square-integrable majorant of the derivatives of $F$, uniform for $\xi$ in a compact set; membership of these derivatives in $H^1(\mbBn)$; and differentiability in $H^1(\mbBn)$ by dominated convergence.

  \emph{Step 1: a uniform majorant.} Fix a compact set $K_0\subset\mbBn$ and let $R_0:=\max_{\xi\in K_0}d(\xi,o)$. We claim that, for $\xi\in K_0$, $x\in\mbBn$ and $|\alpha|\leqslant2$,
  \begin{equation}\label{equ:orbit-map-majorant}
    |\partial_\xi^\alpha F(x,\xi)|+\bigl|\nabla_{\mbBn,x}\partial_\xi^\alpha F(x,\xi)\bigr|\leqslant G(x):=C\mre^{-(\beta_+-\ve)d(x,o)},
  \end{equation}
  with $C=C(n,\lambda,K_0,\ve)$. We distinguish two cases.

  Where $d(x,\xi)\leqslant1$, the pair $(x,\xi)$ ranges over the compact set $\{(x,\xi):\xi\in K_0,\ d(x,\xi)\leqslant1\}$, on which the smooth function $F$ and its derivatives are bounded; moreover $|\nabla_{\mbBn}w|=\rho^{-1}|\nabla w|\leqslant\frac{1}{2}|\nabla w|$. Since also $d(x,o)\leqslant1+R_0$ there, so that $\mre^{-(\beta_+-\ve)d(x,o)}\geqslant\mre^{-\beta_+(1+R_0)}$, \eqref{equ:orbit-map-majorant} holds in this case.

  Where $d(x,\xi)\geqslant1$, we first note that $|x-\xi|\geqslant c(K_0)>0$. Indeed, let $2\eta$ be the Euclidean distance from $K_0$ to $\partial\mbBn$. If $|x-\xi|\leqslant\eta$, then $1-|x|^2\geqslant\eta$ and $1-|\xi|^2\geqslant2\eta$, so the distance formula of this subsection gives $\cosh d(x,\xi)-1\leqslant\eta^{-2}|x-\xi|^2$, which is less than $\cosh1-1$ if $|x-\xi|<\eta\sqrt{\cosh1-1}$. Hence $|x-\xi|\geqslant\eta\sqrt{\cosh1-1}$ whenever $d(x,\xi)\geqslant1$. Next, by the distance formula, $d(x,\xi)=\phi(t)$ with $\phi(t):=\operatorname{arcosh}(1+\mre^{t})$ and
  \begin{equation*}
    t=t(x,\xi):=\log2+\log|x-\xi|^2-\log(1-|x|^2)-\log(1-|\xi|^2)\geqslant t_0:=\log(\cosh1-1).
  \end{equation*}
  We collect three bounds, with constants depending only on $n$ and $K_0$.
  \begin{enumerate}[label=\textnormal{(\arabic*)}]
    \item Since $\phi'(t)=(1+2\mre^{-t})^{-\frac{1}{2}}$, the derivatives $\phi',\phi'',\phi'''$ are bounded on $[t_0,\infty)$.
    \item For $1\leqslant|\alpha'|\leqslant2$, $\partial_\xi^{\alpha'} t$ is a sum of $\partial_\xi^{\alpha'}\log|x-\xi|^2$, which is bounded by $C|x-\xi|^{-|\alpha'|}\leqslant C$, and $-\partial_\xi^{\alpha'}\log(1-|\xi|^2)$, which is bounded on $K_0$.
    \item $\nabla_xt=\frac{2(x-\xi)}{|x-\xi|^2}+\frac{2x}{1-|x|^2}$, so $|\nabla_xt|\leqslant C+\rho(x)\leqslant C\rho(x)$, as $\rho\geqslant2$; and for $1\leqslant|\alpha'|\leqslant2$ only the first term depends on $\xi$, so $|\nabla_x\partial_\xi^{\alpha'} t|\leqslant C|x-\xi|^{-1-|\alpha'|}\leqslant C\rho(x)$.
  \end{enumerate}
  By the chain rule and (1)--(3), the derivatives $\partial_\xi^{\alpha'} d(x,\xi)$ with $|\alpha'|\leqslant2$ are bounded, and $|\nabla_x\partial_\xi^{\alpha'} d(x,\xi)|\leqslant C\rho(x)$. Applying the chain rule once more to $F=\mcU(d(x,\xi))$, $\partial_\xi^\alpha F$ is a sum of terms $\mcU^{(j)}(d(x,\xi))$, $j\leqslant2$, times products $\partial_\xi^{\alpha_1}d\cdots\partial_\xi^{\alpha_j}d$ of bounded factors with $\alpha_1+\dots+\alpha_j=\alpha$, and $\nabla_x\partial_\xi^\alpha F$ is a sum of such terms with $j\leqslant3$ and with one factor $\nabla_x\partial_\xi^{\alpha'} d$. Since $|\nabla_{\mbBn}w|=\rho^{-1}|\nabla w|$, we find
  \begin{align*}
    |\partial_\xi^\alpha F(x,\xi)|+\bigl|\nabla_{\mbBn,x}\partial_\xi^\alpha F(x,\xi)\bigr|
    &\leqslant C\max_{j\leqslant3}|\mcU^{(j)}(d(x,\xi))|
    \leqslant C\mre^{-(\beta_+-\ve)d(x,\xi)}\\
    &\leqslant C\mre^{(\beta_+-\ve)R_0}\mre^{-(\beta_+-\ve)d(x,o)},
  \end{align*}
  where in the second inequality we have used \eqref{equ:ground-state-higher-tail}, and in the last inequality $d(x,o)\leqslant d(x,\xi)+R_0$. This proves \eqref{equ:orbit-map-majorant}.

  Finally, $G\in L^2(\mbBn)$: in polar coordinates about $o$, by \eqref{equ:polar-volume},
  \begin{equation*}
    \int_{\mbBn}G^2\diopv=C^2|\mbS^{n-1}|\int_0^\infty\mre^{-2(\beta_+-\ve)r}\sinh^{n-1}r\diopr<\infty,
  \end{equation*}
  since $2(\beta_+-\ve)>n-1$ by the first display of this proof.

  \emph{Step 2: the derivatives lie in $H^1(\mbBn)$.} Fix $\xi\in K_0$ and $|\alpha|\leqslant2$, and let $w:=\partial_\xi^\alpha F(\cdot,\xi)$, a smooth function on $\mbBn$ with $|w|+|\nabla_{\mbBn}w|\leqslant G$ by Step~1. Let $\chi\in C^\infty(\mbR)$ satisfy $\chi=1$ on $(-\infty,1]$ and $\chi=0$ on $[2,\infty)$, and set $\chi_R(x):=\chi(d(x,o)/R)$ for $R\geqslant1$. Then $\chi_Rw\in\Ccit(\mbBn)$ and $|\nabla_{\mbBn}\chi_R|\leqslant\frac{C}{R}$, so
  \begin{equation*}
    |w-\chi_Rw|+\bigl|\nabla_{\mbBn}(w-\chi_Rw)\bigr|\leqslant\Bigl(2+\frac{C}{R}\Bigr)G\,\mathbf{1}_{\{d(\cdot,o)\geqslant R\}}.
  \end{equation*}
  Since $G\in L^2(\mbBn)$, the right-hand side tends to $0$ in $L^2(\mbBn)$ as $R\to\infty$. Hence $\chi_Rw\to w$ in $H^1(\mbBn)$, so $w\in H^1(\mbBn)$.

  \emph{Step 3: differentiability in $H^1(\mbBn)$.} For $|\alpha|\leqslant2$ write $w_\alpha(\xi):=\partial_\xi^\alpha F(\cdot,\xi)\in H^1(\mbBn)$. We first show that $\xi\mapsto w_\alpha(\xi)$ is continuous on $K_0$. If $\xi_k\to\xi$ in $K_0$, then $w_\alpha(\xi_k)-w_\alpha(\xi)$ and its gradient tend to $0$ pointwise, since $F$ is smooth, and are dominated by $2G\in L^2(\mbBn)$ by Step~1; dominated convergence gives $\norm{w_\alpha(\xi_k)-w_\alpha(\xi)}_{H^1(\mbBn)}\to0$. Next, let $|\alpha|\leqslant1$, let $\xi$ be an interior point of $K_0$, and let $h\neq0$ be so small that the segment from $\xi$ to $\xi+he_i$ lies in $K_0$. By the fundamental theorem of calculus, pointwise in $x$,
  \begin{equation*}
    \frac{w_\alpha(\xi+he_i)-w_\alpha(\xi)}{h}-w_{\alpha+e_i}(\xi)=\int_0^1\bigl(w_{\alpha+e_i}(\xi+she_i)-w_{\alpha+e_i}(\xi)\bigr)\diop s,
  \end{equation*}
  and the same identity holds for the gradients in $x$. The integrand tends to $0$ pointwise as $h\to0$ and, since $|\alpha+e_i|\leqslant2$, it is dominated by $2G$ by Step~1. Dominated convergence again shows that the left-hand side tends to $0$ in $H^1(\mbBn)$. Hence $\xi\mapsto w_\alpha(\xi)$, $|\alpha|\leqslant1$, has the partial derivatives $w_{\alpha+e_i}(\xi)$ in $H^1(\mbBn)$, and these are continuous. A map into a Banach space with continuous partial derivatives is continuously differentiable; applying this to $w_0$ and to $w_{e_i}$ shows that $\xi\mapsto\mcU_\xi=w_0(\xi)$ is $C^2$ into $H^1(\mbBn)$ on the interior of $K_0$, with derivatives $w_\alpha(\xi)$. Since $K_0$ is arbitrary, the same holds on $\mbBn$.

  \emph{Proof of \ref{item:pairing-center-partial}.} By \ref{item:orbit-map-C2}, $\xi\mapsto\lag f,\mcU_\xi\rag_\lambda$ is the composition of the $C^2$ map $\xi\mapsto\mcU_\xi$ with the bounded linear functional $\lag f,\cdot\rag_\lambda$. So it is of class $C^2$ on $\mbBn$, and its partial derivatives are $\lag f,\partial_{\xi_i}\mcU_\xi\rag_\lambda$. This finishes the proof of Lemma~\ref{lem:pairing-center-derivative}.
\end{proof}

\subsection{Elementary identities and decay at infinity}\label{subsec:elementary}
\begin{lemma}\label{lem:ground-state-normalization}
  For every $\xi\in\mbBn$,
  \begin{gather}
    (-\Delta_{\mbBn}-\lambda)\mcU_\xi=\mcU_\xi^{q},\label{equ:ground-state-equation}\\
    \norm{\mcU_\xi}_\lambda^2
    =\norm{\mcU_\xi}_{L^{2^\ast}(\mbBn)}^{2^\ast}=S_\lambda^{\frac{n}{2}},\label{equ:ground-state-energy-normalization}
  \end{gather}
  and
  \begin{equation}
    \norm{\mcU_\xi^{q}}_{L^{(2^{\ast})'}(\mbBn)}
    =S_\lambda^{\frac{n+2}{4}}.\label{equ:ground-state-dual-norm}
  \end{equation}
\end{lemma}
\begin{proof}
  We first prove \eqref{equ:ground-state-equation}. Let $\xi\in\mbBn$. 
  Since $\gamma_\xi^{-1}$ is an isometry of $\mbBn$, it commutes with $-\Delta_{\mbBn}$ and preserves $\diopv$. Hence, by \eqref{equ:EL},
  \begin{equation*}
    (-\Delta_{\mbBn}-\lambda)\mcU_\xi
    =\bigl((-\Delta_{\mbBn}-\lambda)\mcU\bigr)\circ\gamma_\xi^{-1}
    =\mcU^{q}\circ\gamma_\xi^{-1}
    =\mcU_\xi^{q},
  \end{equation*}
  which is \eqref{equ:ground-state-equation}.

  Next we prove \eqref{equ:ground-state-energy-normalization}.
  On the one hand, by the definition of $\norm{\cdot}_\lambda$ and integration by parts, we find
  \begin{equation}
    \norm{\mcU_\xi}_\lambda^2
    =\int_{\mbBn}\bigl(|\nabla_{\mbBn} \mcU_\xi|^2-\lambda\, \mcU_\xi^2\bigr)\diopv
     =\int_{\mbBn}\bigl(-\Delta_{\mbBn}\mcU_\xi-\lambda\,\mcU_\xi\bigr)\mcU_\xi\diopv
     =\int_{\mbBn}\mcU_\xi^{2^{\ast}}\diopv>0,\label{equ:normmcUxiLambda}
  \end{equation}
  where in the last equality we have used \eqref{equ:ground-state-equation}.
  On the other hand, since $\mcU_\xi\in\mcZ$ is an extremal of \eqref{equ:PS} by \eqref{equ:Zo}, equality in \eqref{equ:PS} at $\mcU_\xi$ gives
  \begin{align*}
    \norm{\mcU_\xi}_\lambda^2
    =S_\lambda\Bigl(\int_{\mbBn}\mcU_\xi^{2^{\ast}}\diopv\Bigr)^{\frac{2}{2^{\ast}}}.
  \end{align*}
  Hence, the above two equalities and the fact that $1-\frac{2}{2^{\ast}}=\frac{2}{n}$ yield
  \begin{equation*}
    \int_{\mbBn}\mcU_\xi^{2^{\ast}}\diopv=S_\lambda^{\frac{n}{2}},
  \end{equation*}
  which, using \eqref{equ:normmcUxiLambda} again, gives \eqref{equ:ground-state-energy-normalization}.

  Finally, we prove \eqref{equ:ground-state-dual-norm}. Since $(2^{\ast})'=\frac{2n}{n+2}$ is the conjugate of $2^\ast$, we have $q(2^{\ast})'=\frac{n+2}{n-2}\cdot\frac{2n}{n+2}=2^{\ast}$. Hence \eqref{equ:ground-state-energy-normalization} gives
  \begin{equation*}
    \norm{\mcU_\xi^{q}}_{L^{(2^{\ast})'}(\mbBn)}^{(2^{\ast})'}
    =\int_{\mbBn}\mcU_\xi^{q(2^{\ast})'}\diopv
    =\int_{\mbBn}\mcU_\xi^{2^{\ast}}\diopv
    =S_\lambda^{\frac{n}{2}},
  \end{equation*}
  and thus \eqref{equ:ground-state-dual-norm} is true, which finishes the proof of Lemma~\ref{lem:ground-state-normalization}.
\end{proof}
\begin{lemma}\label{lem:duality}
  For every $f\in H^1(\mbBn)$ and $\xi\in\mbBn$,
  \begin{equation}\label{equ:duality}
    \lag f,\mcU_\xi\rag_\lambda=\int_{\mbBn}f\mcU_\xi^q\diopv.
  \end{equation}
\end{lemma}
\begin{proof}
  Notice that both sides of \eqref{equ:duality} are continuous in $f$ with respect to $\norm{\cdot}_\lambda$, because the left-hand side is controlled by the Cauchy--Schwarz inequality, and the right-hand side by H\"older's inequality, \eqref{equ:PS} and \eqref{equ:ground-state-dual-norm}. Hence, by density, it suffices to consider $f\in\Ccit(\mbBn)$. For such $f$, by the definition of $\lag\cdot,\cdot\rag_\lambda$ and integration by parts, we find
  \begin{equation*}
    \lag f,\mcU_\xi\rag_\lambda
    =\int_{\mbBn}\bigl(\nabla_{\mbBn}f\cdot\nabla_{\mbBn}\mcU_\xi-\lambda f\mcU_\xi\bigr)\diopv
    =\int_{\mbBn}f\bigl(-\Delta_{\mbBn}\mcU_\xi-\lambda\,\mcU_\xi\bigr)\diopv
    =\int_{\mbBn}f\mcU_\xi^{q}\diopv,
  \end{equation*}
  where in the last equality we have used \eqref{equ:ground-state-equation}. This finishes the proof of Lemma~\ref{lem:duality}.
\end{proof}
\begin{lemma}\label{lem:dist-formula}
  For $f\in H^1(\mbBn)$,
  \begin{equation}\label{equ:dist-formula}
    \dist(f,\mcZ)^2=\norm{f}_\lambda^2-S_\lambda^{-\frac{n}{2}}\sup_{\xi\in\mbBn}\lag f,\mcU_\xi\rag_\lambda^2
    =\norm{f}_\lambda^2-S_\lambda^{-\frac{n}{2}}\sup_{\xi\in\mbBn}\left(\int_{\mbBn}f\mcU_\xi^q\diopv\right)^2.
  \end{equation}
\end{lemma}
\begin{proof}
  Fix $\xi\in\mbBn$. By \eqref{equ:ground-state-energy-normalization}, for every $c\in\mbR$,
  \begin{equation*}
    \norm{f-c\,\mcU_\xi}_\lambda^2=\norm{f}_\lambda^2-2c\lag f,\mcU_\xi\rag_\lambda+c^2S_\lambda^{\frac{n}{2}}.
  \end{equation*}
  Since the right-hand side is a strictly convex quadratic polynomial in $c$, it is minimized at $c=S_\lambda^{-\frac{n}{2}}\lag f,\mcU_\xi\rag_\lambda$. Hence
  \begin{equation*}
    \min_{c\in\mbR}\norm{f-c\,\mcU_\xi}_\lambda^2=\norm{f}_\lambda^2-S_\lambda^{-\frac{n}{2}}\lag f,\mcU_\xi\rag_\lambda^2,
  \end{equation*}
  and, by continuity in $c$, the infimum over $c\neq0$ equals this minimum. Taking the infimum over $\xi\in\mbBn$ and recalling the definition \eqref{equ:Zo} of $\mcZ$, we obtain
  \begin{equation*}
    \dist(f,\mcZ)^2
    =\inf_{\xi\in\mbBn}\inf_{c\neq0}\norm{f-c\,\mcU_\xi}_\lambda^2
    =\norm{f}_\lambda^2-S_\lambda^{-\frac{n}{2}}\sup_{\xi\in\mbBn}\lag f,\mcU_\xi\rag_\lambda^2,
  \end{equation*}
  which is the first equality in \eqref{equ:dist-formula}. The second equality of \eqref{equ:dist-formula} then follows directly from the first equality and \eqref{equ:duality}.

  Finally, \eqref{equ:dist-formula} obviously yields $\dist(f,\mcZ)^2\leqslant\norm{f}_\lambda^2$, which gives $\sigma(f)\in[0,1]$ in \eqref{equ:deficit-sigma}. This finishes the proof of Lemma~\ref{lem:dist-formula}.
\end{proof}
The pairing decays as the center escapes to infinity.
\begin{lemma}\label{lem:overlap}
  Let $f\in L^{2^{\ast}}(\mbBn)$. Then
  \begin{equation*}
    \int_{\mbBn}f\,\mcU_\xi^{q}\diopv\longrightarrow0
    \qquad\text{as } d(o,\xi)\to\infty .
  \end{equation*}
  In particular $\lag\mcU,\mcU_\xi\rag_\lambda\to0$ as $d(o,\xi)\to\infty$.
\end{lemma}
\begin{proof}
  Let $\ve>0$. Since compactly supported functions are dense in $L^{2^{\ast}}(\mbBn)$, there is a compactly supported $g$ with $\norm{f-g}_{L^{2^{\ast}}(\mbBn)}<\ve$, and we split
  \begin{equation}\label{equ:intmbRnfmcUdslj}
    \int_{\mbBn}f\,\mcU_\xi^{q}\diopv
    =\int_{\mbBn}(f-g)\,\mcU_\xi^{q}\diopv+\int_{\mbBn}g\,\mcU_\xi^{q}\diopv.
  \end{equation}
  On the one hand, by H\"older's inequality and \eqref{equ:ground-state-dual-norm}, we find
  \begin{equation}\label{equ:intmbRnfgalsjfso}
    \Bigl|\int_{\mbBn}(f-g)\,\mcU_\xi^{q}\diopv\Bigr|
    \leqslant\norm{f-g}_{L^{2^{\ast}}(\mbBn)}\norm{\mcU_\xi^{q}}_{L^{(2^{\ast})'}(\mbBn)}
    \leqslant\ve S_\lambda^{\frac{n+2}{4}}
  \end{equation}
  uniformly in $\xi$.
  On the other hand, let $K$ be a compact set containing $\supp g$. By H\"older's inequality and $q(2^{\ast})'=2^{\ast}$, we find
  \begin{equation*}
    \Bigl|\int_{\mbBn}g\,\mcU_\xi^{q}\diopv\Bigr|
    \leqslant\norm{g}_{L^{2^{\ast}}(\mbBn)}\Bigl(\int_{K}\mcU_\xi^{2^{\ast}}\diopv\Bigr)^{\frac{n+2}{2n}}
    =\norm{g}_{L^{2^{\ast}}(\mbBn)}\Bigl(\int_{\gamma_\xi^{-1}K}\mcU^{2^{\ast}}\diopv\Bigr)^{\frac{n+2}{2n}},
  \end{equation*}
  where in the last equality we have used that $\gamma_\xi$ preserves $\diopv$. For $x\in K$, since $\gamma_\xi$ is an isometry with $\gamma_\xi(o)=\xi$, we have
  \begin{equation*}
    d(o,\gamma_\xi^{-1}x)=d(\xi,x)\geqslant d(o,\xi)-\max_{K}d(o,\cdot)=:R_\xi.
  \end{equation*}
  Thus $\gamma_\xi^{-1}K\subset\{y\in\mbBn:d(o,y)\geqslant R_\xi\}$, and $R_\xi\to\infty$ as $d(o,\xi)\to\infty$. As a result, by noticing $\mcU^{2^{\ast}}\in L^1(\mbBn)$, we obtain
  \begin{equation*}
    \int_{\gamma_\xi^{-1}K}\mcU^{2^{\ast}}\diopv
    \leqslant\int_{\{d(o,\cdot)\geqslant R_\xi\}}\mcU^{2^{\ast}}\diopv\longrightarrow0
    \qquad\text{as } d(o,\xi)\to\infty .
  \end{equation*}
  Combining this estimate with \eqref{equ:intmbRnfmcUdslj} and \eqref{equ:intmbRnfgalsjfso}, we find
  \begin{equation*}
    \limsup_{d(o,\xi)\to\infty}\Bigl|\int_{\mbBn}f\,\mcU_\xi^{q}\diopv\Bigr|\leqslant\ve S_\lambda^{\frac{n+2}{4}},
  \end{equation*}
  and letting $\ve\downarrow0$ proves the claim.

  Finally, since $\mcU\in L^{2^{\ast}}(\mbBn)$ by \eqref{equ:ground-state-energy-normalization}, the last assertion follows from the case $f=\mcU$ and \eqref{equ:duality}. This finishes the proof of Lemma~\ref{lem:overlap}.
\end{proof}

\begin{lemma}\label{lem:dist-positive}
  Every $h\in H^1(\mbBn)\setminus(\mcZ\cup\{0\})$ satisfies $\dist(h,\mcZ)>0$.
\end{lemma}
\begin{proof}
  We argue by contradiction by supposing $\dist(h,\mcZ)=0$. Then there are $\{z_k\}\subset\mcZ$ with $\norm{h-z_k}_\lambda\to0$, and hence, by \eqref{equ:PS}, 
  \begin{align*}
    \norm{h-z_k}_{L^{2^{\ast}}(\mbBn)}^2\leqslant S_\lambda^{-1}\norm{h-z_k}_\lambda^2\to 0.
  \end{align*}
  Since $\delta(z_k)=0$ and $\delta$ is continuous in these two norms, we find $\delta(h)=0$. The equality case of \eqref{equ:PS} then forces $h\in\mcZ$, which causes a contradiction and thus finishes the proof of Lemma~\ref{lem:dist-positive}.
\end{proof}

\subsection[The conformal change of variables]{The conformal change of variables}\label{subsec:conformal-dictionary}
With $\rho$ as in Section~\ref{sec:intro}, the map
\begin{equation*}
  f\longmapsto v:=\rho^{\frac{n-2}{2}}f
\end{equation*}
is a bijection of $H^1(\mbBn)$ onto $H_0^1(B_1)$ under which
\begin{equation}\label{equ:dictionary}
  \begin{aligned}
    &\int_{\mbBn}(|\nabla_{\mbBn}f|^2-\lambda_0f^2)\diopv=\int_{B_1}|\nabla v|^2\diopx,
    \quad\int_{\mbBn}f^2\diopv=\int_{B_1}\rho^2v^2\diopx,\\
    &\int_{\mbBn}|f|^{2^{\ast}}\diopv=\int_{B_1}|v|^{2^{\ast}}\diopx
  \end{aligned}
\end{equation}
\cite[(2.9)--(2.10)]{BhaktaSandeep2012}. The first identity is the conformal invariance of the conformal Laplacian, whose zeroth-order term is exactly $-\lambda_0=-\frac{n(n-2)}{4}$; the other two follow from $n-(n-2)=2$ and $2^{\ast}\cdot\frac{n-2}{2}=n$.
With $\Lambda:=\lambda-\lambda_0$ one has $\norm{\cdot}_\lambda^2=\norm{\cdot}_{\lambda_0}^2-\Lambda\int_{\mbBn}(\cdot)^2\diopv$, so the change of variables turns the quadratic form into $\int_{B_1}|\nabla\cdot|^2\diopx-\Lambda\int_{B_1}\rho^2(\cdot)^2\diopx$.

The Euclidean model is the Aubin--Talenti bubble
\begin{equation}\label{equ:AT-bubble}
  \alpha_n:=[n(n-2)]^{\frac{n-2}{4}},\qquad V(x):=\alpha_n(1+|x|^2)^{-\frac{n-2}{2}},
\end{equation}
normalized by $-\Delta V=V^{q}$ on $\mbRn$; it attains the Euclidean Sobolev constant $S_n$ \cite[Th\'eor\`eme~7 and Lemme~3]{Aubin1976}, \cite{Talenti1976}.

\subsection{The linearized spectrum}\label{subsec:linearized-spectrum}
The weighted eigenvalue problem \eqref{equ:const-linearized} is understood in the weak sense, that is,
\begin{equation*}
  \lag\Psi,\vp\rag_\lambda=\mu\int_{\mbBn}\mcU^{q-1}\Psi\vp\diopv\quad\text{for every }\vp\in H^1(\mbBn).
\end{equation*}
By \cite[Lemmas~3.1 and~3.2]{BGKM2025}, the eigenvalues of \eqref{equ:const-linearized} are positive, have finite multiplicity and tend to $\pit$. We arrange the eigenvalues in increasing order without counting multiplicity, and let $\mu_i$ denote the $i$-th eigenvalue. By \cite[Proposition~3.1]{BGKM2025} and \cite[Theorem~1.1]{GangulySandeep2015}, $\mu_1=1$ and $\mu_2=q$. In particular, $\mu_3>\mu_2=q$.

For $1\leqslant i\leqslant n$ let $\Phi_i:=\partial_{\xi_i}\mcU_\xi|_{\xi=o}$, which belongs to $H^1(\mbBn)$ by Lemma~\ref{lem:pairing-center-derivative}\ref{item:orbit-map-C2}.

\begin{lemma}\label{lem:linearized-eigenspaces}
  \begin{enumerate}[label=\textnormal{(\alph*)}]
    \item\label{item:eigen-one} The eigenspace of \eqref{equ:const-linearized} for the eigenvalue $1$ is spanned by $\mcU$.
    \item\label{item:eigen-translation} For $1\leqslant i\leqslant n$, $\Phi_i(x)=-2\mcU'(d(x,o))\frac{x_i}{|x|}$ for $x\neq0$, and for every $f\in H^1(\mbBn)$,
    \begin{equation*}
      \partial_{\xi_i}\lag f,\mcU_\xi\rag_\lambda\big|_{\xi=o}=\lag f,\Phi_i\rag_\lambda=q\int_{\mbBn}f\mcU^{q-1}\Phi_i\diopv.
    \end{equation*}
    In particular, $\Phi_i$ is a weak solution of the linearized equation $(-\Delta_{\mbBn}-\lambda)\Phi=q\mcU^{q-1}\Phi$, that is, an eigenfunction of \eqref{equ:const-linearized} with eigenvalue $q$.
    \item\label{item:eigen-q} The ground state $\mcU$ is a \emph{nondegenerate} solution of \eqref{equ:EL}, that is, the eigenspace of \eqref{equ:const-linearized} for the eigenvalue $q$, which consists of the weak solutions of the linearized equation, is spanned by $\Phi_1,\ldots,\Phi_n$.
    \item\label{item:eigen-tangent} The functions $\mcU,\Phi_1,\ldots,\Phi_n$ are linearly independent, and their span is the tangent space of $\mcZ$ at $\mcU$. Hence this tangent space is the sum of the eigenspaces for the eigenvalues $1$ and $q$.
  \end{enumerate}
\end{lemma}
\begin{proof}
  Assertion \ref{item:eigen-one} is part of \cite[Proposition~3.1]{BGKM2025}.

  For \ref{item:eigen-translation}, by Lemma~\ref{lem:pairing-center-derivative}\ref{item:orbit-map-C2}, $\Phi_i$ is the pointwise derivative of $\mcU(d(x,\xi))$ in $\xi_i$ at $\xi=o$. Differentiating the distance formula of Section~\ref{subsec:geometry-analysis} in $\xi_i$ at $\xi=o$ gives $\partial_{\xi_i}d(x,\xi)|_{\xi=o}=-\frac{2x_i}{|x|}$ for $x\neq0$, and the chain rule gives the formula for $\Phi_i$. By Lemma~\ref{lem:pairing-center-derivative}\ref{item:pairing-center-partial}, $\partial_{\xi_i}\lag f,\mcU_\xi\rag_\lambda|_{\xi=o}=\lag f,\Phi_i\rag_\lambda$. The second equality is the weak form of the linearized equation. By \cite[Theorem~1.3]{ManciniSandeep2008}, $\mcU$ is the unique positive radial solution of \eqref{equ:EL} in $H^1(\mbBn)$, and the present range of $n$ and $\lambda$ is the critical case of the hypothesis (1.2) of \cite{GangulySandeep2015}. So \cite[Theorem~1.1]{GangulySandeep2015} applies: the kernel elements of \cite{GangulySandeep2015} lie in $H^1(\mbBn)$ and solve $(-\Delta_{\mbBn}-\lambda)\Phi=q\mcU^{q-1}\Phi$ weakly in $H^1(\mbBn)$ (their operator $L$ is our $-\Delta_{\mbBn}-\lambda-q\mcU^{q-1}$). By \cite[(4.2)]{GangulySandeep2015} these kernel elements are $\frac{\mrd}{\mrd t}|_{t=0}\,\mcU\circ\tau_{te_i}$, where $\tau_b$ is the hyperbolic translation \cite[(2.2)]{GangulySandeep2015}, an isometry with $\tau_b(o)=b$ and, by its formula, $\tau_b(-b)=o$. Hence $d(\tau_{te_i}x,o)=d(x,-te_i)$, that is, $\mcU\circ\tau_{te_i}=\mcU_{-te_i}$, and their kernel elements are $-\Phi_i$. So $\Phi_i$ weakly solves the linearized equation, and testing with $f$ gives the second equality.

  For \ref{item:eigen-q}, by \cite[Theorem~1.1]{GangulySandeep2015}, the weak solutions of the linearized equation are spanned by the kernel elements of \cite{GangulySandeep2015}, which are $-\Phi_1,\ldots,-\Phi_n$ by the proof of \ref{item:eigen-translation}.

  For \ref{item:eigen-tangent}, $\mcU$ is radial, and by \ref{item:eigen-translation}, $\Phi_i=-2\mcU'(r)\theta_i$, where $\mcU'$ is not identically zero by Lemma~\ref{lem:ground-state}. Since $\theta_1,\ldots,\theta_n$ are linearly independent spherical harmonics of degree $1$ and are orthogonal to the constants, the functions $\mcU,\Phi_1,\ldots,\Phi_n$ are linearly independent. By Lemma~\ref{lem:pairing-center-derivative}\ref{item:orbit-map-C2}, the map $(c,\xi)\mapsto c\mcU_\xi$ is of class $C^2$ into $H^1(\mbBn)$, with partial derivatives $\mcU$ and $\Phi_1,\ldots,\Phi_n$ at $(1,o)$. Hence the tangent space of $\mcZ$ at $\mcU$ is their span, and the last assertion follows from \ref{item:eigen-one}, \ref{item:eigen-translation} and \ref{item:eigen-q}. This finishes the proof of Lemma~\ref{lem:linearized-eigenspaces}.
\end{proof}

\begin{lemma}\label{lem:linearized-coercivity}
  Let $\Psi\in H^1(\mbBn)$ be $\lag\cdot,\cdot\rag_\lambda$-orthogonal to $\mcU$ and to $\Phi_1,\ldots,\Phi_n$. Then
  \begin{equation}\label{equ:linearized-coercivity}
    \norm{\Psi}_\lambda^2\geqslant\mu_3\int_{\mbBn}\mcU^{q-1}\Psi^2\diopv,
  \end{equation}
  with equality if $\Psi$ is an eigenfunction of \eqref{equ:const-linearized} with eigenvalue $\mu_3$. In particular, if $\Psi\neq0$, then
  \begin{equation*}
    \norm{\Psi}_\lambda^2-q\int_{\mbBn}\mcU^{q-1}\Psi^2\diopv\geqslant(\mu_3-q)\int_{\mbBn}\mcU^{q-1}\Psi^2\diopv>0.
  \end{equation*}
\end{lemma}
\begin{proof}
  We follow the argument in the proof of \cite[Lemma~4.1]{BGKM2025}.
  By \cite[Lemma~3.2]{BGKM2025} and the spectral theorem, $H^1(\mbBn)$ has a $\lag\cdot,\cdot\rag_\lambda$-orthonormal basis $\{\Psi_k\}_{k\geqslant1}$ of eigenfunctions of \eqref{equ:const-linearized}; denote by $\nu_k$ the eigenvalue of $\Psi_k$. Testing the weak form of \eqref{equ:const-linearized} for $\Psi_k$ with $\Psi_m$, we find
  \begin{equation*}
    \delta_{km}=\lag\Psi_k,\Psi_m\rag_\lambda=\nu_k\int_{\mbBn}\mcU^{q-1}\Psi_k\Psi_m\diopv,
    \quad\text{that is,}\quad
    \int_{\mbBn}\mcU^{q-1}\Psi_k\Psi_m\diopv=\nu_k^{-1}\delta_{km}.
  \end{equation*}
  Expanding $\Psi=\sum_kc_k\Psi_k$ with $c_k:=\lag\Psi,\Psi_k\rag_\lambda$, we therefore have
  \begin{equation*}
    \norm{\Psi}_\lambda^2=\sum_kc_k^2,
    \qquad
    \int_{\mbBn}\mcU^{q-1}\Psi^2\diopv=\sum_k\nu_k^{-1}c_k^2.
  \end{equation*}
  By Lemma~\ref{lem:linearized-eigenspaces}\ref{item:eigen-one} and~\ref{item:eigen-q}, $\mcU$ and $\Phi_1,\ldots,\Phi_n$ span the eigenspaces of the eigenvalues $\mu_1=1$ and $\mu_2=q$. Since $\Psi$ is orthogonal to them, $c_k=0$ whenever $\nu_k\in\{1,q\}$, and every remaining $\nu_k$ satisfies $\nu_k\geqslant\mu_3$ by the definition of $\mu_3$. Hence
  \begin{equation*}
    \norm{\Psi}_\lambda^2=\sum_{\nu_k\geqslant\mu_3}c_k^2\geqslant\mu_3\sum_{\nu_k\geqslant\mu_3}\nu_k^{-1}c_k^2=\mu_3\int_{\mbBn}\mcU^{q-1}\Psi^2\diopv,
  \end{equation*}
  which is \eqref{equ:linearized-coercivity}. If $\Psi$ is an eigenfunction with eigenvalue $\mu_3$, then $c_k=0$ unless $\nu_k=\mu_3$, and equality holds.

  Finally, subtracting $q\int_{\mbBn}\mcU^{q-1}\Psi^2\diopv$ from both sides of \eqref{equ:linearized-coercivity} gives the last inequality with $\geqslant$, and its right-hand side is positive since $\mu_3>q$ (see the beginning of this subsection) and $\int_{\mbBn}\mcU^{q-1}\Psi^2\diopv>0$ for $\Psi\neq0$, as $\mcU>0$. This finishes the proof of Lemma~\ref{lem:linearized-coercivity}.
\end{proof}
\section{Explicit local stability}\label{sec:crit-local}
This section proves Theorem~\ref{thm:crit-local}. We decompose $f$ at a nearest point of the extremal cone, normalize $\norm{f}_\lambda^2=S_\lambda^{\frac{n}{2}}$, and expand $\norm{f}_{L^{2^{\ast}}(\mbBn)}^{2^{\ast}}$, the only non-quadratic term of the deficit, to bound the orthogonal remainder.
\subsection{Modulation at the nearest point}\label{subsec:local-modulation}
\begin{lemma}\label{lem:const-projection}
  Let $f\in H^1(\mbBn)$ satisfy $0<\sigma(f)<1$ and let $d_f:=\dist(f,\mcZ)$.
  Then the distance is attained, and after an isometry and possibly replacing $f$ by $-f$
  the nearest point is $a\mcU$ with $a>0$, so that
  \begin{equation}\label{equ:const-modulation}
    f=a\mcU+d_f\psi,\qquad\norm{\psi}_\lambda=1,\qquad a^2S_\lambda^{\frac{n}{2}}=\norm{f}_\lambda^2-d_f^2.
  \end{equation}
  Moreover
  \begin{equation}\label{equ:const-orthogonality}
    \lag\psi,\mcU\rag_\lambda=\lag\psi,\Phi_i\rag_\lambda=0\quad(1\leqslant i\leqslant n),
  \end{equation}
  equivalently $\int\mcU^{q-1}\mcU\psi\diopv=\int\mcU^{q-1}\Phi_i\psi\diopv=0$, and
  \begin{equation}\label{equ:const-spectral-moment}
    \int_{\mbBn}\mcU^{q}\psi\diopv=0,\qquad\int_{\mbBn}\mcU^{q-1}\psi^2\diopv\leqslant\frac{1}{\mu_3}.
  \end{equation}
\end{lemma}
\begin{proof}
  We first show that the distance is attained at a point of $\mcZ$, then derive \eqref{equ:const-modulation} and \eqref{equ:const-orthogonality} from the stationarity of the distance at this point, and finally deduce \eqref{equ:const-spectral-moment}.

  \emph{Step 1: the distance is attained.} By \eqref{equ:dist-formula}, $\sigma(f)<1$ forces $\sup_\xi|\lag f,\mcU_\xi\rag_\lambda|>0$, and by Lemma~\ref{lem:pairing-center-derivative}\ref{item:pairing-center-partial} the pairing $\lag f,\mcU_\xi\rag_\lambda$ is continuous in $\xi$. By Lemma~\ref{lem:overlap}, the pairing tends to $0$ as $d(o,\xi)\to\infty$. Hence a maximizing sequence of centers stays in a bounded, hence relatively compact, set, and a subsequence converges to a center $\xi_\ast$ at which the positive supremum is attained. At $\xi_\ast$, the map $c\mapsto\norm{f-c\mcU_{\xi_\ast}}_\lambda^2$ is a strictly convex quadratic minimized at $c_\ast:=S_\lambda^{-\frac{n}{2}}\lag f,\mcU_{\xi_\ast}\rag_\lambda\neq0$ by Lemma~\ref{lem:ground-state-normalization}, and by \eqref{equ:dist-formula} its minimum value is $d_f^2$. So the distance is attained at $c_\ast\mcU_{\xi_\ast}$, which lies in $\mcZ$ and not merely in $\mcZ\cup\{0\}$.

  \emph{Step 2: \eqref{equ:const-modulation} and \eqref{equ:const-orthogonality}.}
  The norms $\norm{\cdot}_\lambda$ and $\norm{\cdot}_{L^{2^{\ast}}(\mbBn)}$ and the cone $\mcZ$ are invariant under composition with isometries and under $f\mapsto-f$. Hence we may assume $\xi_\ast=o$ and $c_\ast>0$. Set $a:=c_\ast=S_\lambda^{-\frac{n}{2}}\lag f,\mcU\rag_\lambda>0$ and $\psi:=\frac{f-a\mcU}{d_f}$, which is well defined since $\sigma(f)>0$ gives $d_f>0$; then $\norm{\psi}_\lambda=1$. Stationarity in the amplitude gives $\lag f-a\mcU,\mcU\rag_\lambda=0$, that is, $\lag\psi,\mcU\rag_\lambda=0$. By the Pythagorean theorem and \eqref{equ:ground-state-energy-normalization}, $\norm{f}_\lambda^2=a^2S_\lambda^{\frac{n}{2}}+d_f^2$, which completes \eqref{equ:const-modulation}.

  Stationarity in the center gives $0=2\lag f,\mcU\rag_\lambda\lag f,\Phi_i\rag_\lambda$: here we differentiate $\xi\mapsto\lag f,\mcU_\xi\rag_\lambda^2$, which attains its maximum at $o$, along $\xi=te_i$ at $t=0$, which is legitimate by Lemma~\ref{lem:pairing-center-derivative}\ref{item:pairing-center-partial}, and the derivative is given by Lemma~\ref{lem:linearized-eigenspaces}\ref{item:eigen-translation}. The first factor is nonzero, so $\lag f,\Phi_i\rag_\lambda=0$. Since $\mcU$ and $\Phi_i$ are eigenfunctions of \eqref{equ:const-linearized} with different eigenvalues, $\lag\mcU,\Phi_i\rag_\lambda=0$. Therefore $\lag\psi,\Phi_i\rag_\lambda=d_f^{-1}\bigl(\lag f,\Phi_i\rag_\lambda-a\lag\mcU,\Phi_i\rag_\lambda\bigr)=0$, which together with $\lag\psi,\mcU\rag_\lambda=0$ proves \eqref{equ:const-orthogonality}. Pairing \eqref{equ:const-linearized} for $\mcU$ and $\Phi_i$ with $\psi$ converts this into the equivalent weighted form.

  \emph{Step 3: \eqref{equ:const-spectral-moment}.} The first identity is the weighted form of $\lag\psi,\mcU\rag_\lambda=0$, since $\mcU^{q}=\mcU^{q-1}\mcU$. For the second, \eqref{equ:const-orthogonality} shows that $\psi$ satisfies the hypotheses of Lemma~\ref{lem:linearized-coercivity}, and \eqref{equ:linearized-coercivity} gives $\mu_3\int_{\mbBn}\mcU^{q-1}\psi^2\diopv\leqslant\norm{\psi}_\lambda^2=1$. This finishes the proof of Lemma~\ref{lem:const-projection}.
\end{proof}
Since $\mu_3>q$, the bound $\frac{1}{\mu_3}$ in \eqref{equ:const-spectral-moment} is strictly smaller than $\frac{1}{q}$. The distance to the cone is measured by $\frac{d_f}{\norm{f}_\lambda}=\sqrt{\sigma(f)}$. For a tube parameter $0<s<1$, let $u:=\frac{s}{1-s}$; equivalently
\begin{equation*}
  s=\frac{u}{1+u},
\end{equation*}
and $s\mapsto u$ is an increasing bijection of $(0,1)$ onto $(0,\infty)$. In this section $u$ always denotes this ratio.

\subsection{An elementary inequality}\label{subsec:local-elementary}
\begin{lemma}\label{lem:taylor-remainder}
  Let $p>2$ and $R_p(t):=|1+t|^p-1-pt-\frac{p(p-1)}{2}t^2$ for $t\in\mbR$.
  \begin{enumerate}[label=\textnormal{(\alph*)}]
    \item\label{item:remainder-small} If $2<p\leqslant3$, then $R_p(t)\leqslant|t|^p$ for every $t\in\mbR$.
    \item\label{item:remainder-large} If $p>3$, then $R_p(t)\leqslant\kappa_p\bigl(C_p|t|^3+|t|^p\bigr)$ for every $t\in\mbR$, where $C_p:=\frac{p(p-1)(p-2)}{6}$ and $\kappa_p:=\max\{1,2^{p-4}\}$.
  \end{enumerate}
  The constant in \ref{item:remainder-small} cannot be decreased. For $3<p\leqslant4$, neither coefficient in \ref{item:remainder-large} can be decreased.
\end{lemma}
\begin{proof}
  We first record an elementary inequality, then prove both bounds on the range $t\leqslant-1$, where no Taylor expansion is needed, then prove \ref{item:remainder-small} and \ref{item:remainder-large} on the range $t\geqslant-1$ by Taylor's formula, and finally prove the optimality statements.

  \emph{An elementary inequality.} For $\theta>0$ and $r\geqslant0$ we use $(1+r)^\theta\leqslant\max\{1,2^{\theta-1}\}\bigl(1+r^\theta\bigr)$. If $0<\theta\leqslant1$, the difference $1+r^\theta-(1+r)^\theta$ equals $0$ at $r=0$ and has derivative $\theta\bigl[r^{\theta-1}-(1+r)^{\theta-1}\bigr]\geqslant0$, so it is nonnegative. If $\theta>1$, the convexity of $r\mapsto r^\theta$ gives $\bigl(\frac{1+r}{2}\bigr)^\theta\leqslant\frac{1+r^\theta}{2}$.

  \emph{The range $t\leqslant-1$.} Write $t=-x$ with $x\geqslant1$. Then $R_p(-x)=(x-1)^p+N(x)$ with $N(x):=-1+px-\frac{p(p-1)}{2}x^2$. The quadratic $-N$ has positive leading coefficient and discriminant $p^2-2p(p-1)=p(2-p)<0$, so $N<0$, and hence $R_p(-x)<(x-1)^p\leqslant x^p=|t|^p$. This gives \ref{item:remainder-small} on this range, and also \ref{item:remainder-large}, because $\kappa_p\geqslant1$.

  \emph{Part \ref{item:remainder-small} for $t\geqslant-1$.} Since $0<p-2\leqslant1$, Taylor's formula with integral remainder gives
  \begin{equation*}
    R_p(t)=p(p-1)\int_0^t(t-r)\bigl[(1+r)^{p-2}-1\bigr]\diopr\qquad(t\geqslant-1).
  \end{equation*}
  For $-1\leqslant t\leqslant0$ the integrand is nonnegative on $[t,0]$, so $R_p(t)\leqslant0$. For $t\geqslant0$, by the elementary inequality with $\theta=p-2\leqslant1$, we find
  \begin{equation*}
    R_p(t)\leqslant p(p-1)\int_0^t(t-r)r^{p-2}\diopr=p(p-1)t^p\left(\frac{1}{p-1}-\frac{1}{p}\right)=t^p.
  \end{equation*}

  \emph{Part \ref{item:remainder-large} for $t\geqslant-1$.} Since $p>3$, Taylor's formula with one more term gives
  \begin{equation*}
    R_p(t)=\frac{p(p-1)(p-2)}{2}\int_0^{t}(t-r)^2(1+r)^{p-3}\diopr\qquad(t\geqslant-1).
  \end{equation*}
  For $-1\leqslant t\leqslant0$ this is nonpositive. For $t\geqslant0$, the elementary inequality with $\theta=p-3$ and
  \begin{equation*}
    \int_0^{t}(t-r)^2r^{p-3}\diopr
    =t^{p}\left(\frac{1}{p-2}-\frac{2}{p-1}+\frac{1}{p}\right)
    =\frac{2\,t^{p}}{p(p-1)(p-2)}.
  \end{equation*}
  give $R_p(t)\leqslant\kappa_p\bigl(C_pt^3+t^{p}\bigr)$.

  \emph{Optimality.} Since $p>2$, $\frac{R_p(t)}{t^p}\to1$ as $t\to\pit$. Hence the constant in \ref{item:remainder-small} cannot be decreased, and for $3<p\leqslant4$ neither can the coefficient $\kappa_p=1$ of $|t|^p$ in \ref{item:remainder-large}. For $p>3$, Taylor's formula gives $\frac{R_p(t)}{t^3}\to C_p$ as $t\to0^+$. Since $t^p=o(t^3)$ there, the coefficient $\kappa_pC_p=C_p$ of $|t|^3$ cannot be decreased when $3<p\leqslant4$. This finishes the proof of Lemma~\ref{lem:taylor-remainder}.
\end{proof}

\subsection{Proof of Theorem~\ref{thm:crit-local}}\label{subsec:local-theorem}
\begin{proposition}\label{prop:local-lower-function}
  The function $\nu$ of \eqref{equ:const-nu-explicit} is continuous and strictly decreasing on $(0,1)$, and has a unique zero $s_0=s_0(n,\lambda)\in(0,1)$. For every $s\in(0,1)$ and every $f\in H^1(\mbBn)\setminus\{0\}$ with $0<\sigma(f)\leqslant s$,
  \begin{equation}\label{equ:local-prop-inequality}
    \delta(f)^2\geqslant\nu(s)\dist(f,\mcZ)^2.
  \end{equation}
  In particular, $\delta(f)^2\geqslant\nu(\sigma(f))\dist(f,\mcZ)^2$ for every $f\in H^1(\mbBn)$ with $0<\sigma(f)<1$.
\end{proposition}
\begin{proof}
  We first study the function $\nu$ and its zero $s_0$, and then prove \eqref{equ:local-prop-inequality} by expanding $\norm{f}_{L^{2^{\ast}}(\mbBn)}^{2^{\ast}}$ around the nearest point given by Lemma~\ref{lem:const-projection}; the expansion is done separately for $n\geqslant6$, $n=5$ and $n=4$.

  \emph{Step 1: the function $\nu$ and the radius $s_0$.} Since $\frac{2}{2^{\ast}}=\frac{n-2}{n}$ and $\frac{q-1}{2}=\frac{2}{n-2}$, the third branch is $1-\frac{q}{\mu_3}-\frac{n-2}{n}u^{\frac{2}{n-2}}$; at $n=6$ it is $1-\frac{q}{\mu_3}-\frac{2}{3}\sqrt{u}$. Since $s\mapsto u=\frac{s}{1-s}$ is an increasing bijection of $(0,1)$ onto $(0,\infty)$ and every penalty term in \eqref{equ:const-nu-explicit} is a positive multiple of a positive power of $u$, each branch of $\nu$ is continuous and strictly decreasing, and
  \begin{equation}\label{equ:const-nu-limits}
    \lim_{s\downarrow0}\nu(s)=1-\frac{q}{\mu_3}>0,\qquad
    \lim_{s\uparrow1}\nu(s)=\nit,
  \end{equation}
  where the first limit is positive since $\mu_3>q$. Hence, by the intermediate value theorem,
  \begin{equation*}
    s_0:=\sup\{s\in(0,1):\nu(s)>0\}
  \end{equation*}
  is a well-defined number in $(0,1)$; by strict monotonicity it is the unique zero of $\nu$, and $\nu>0$ on $(0,s_0)$. Writing $u_0:=\frac{s_0}{1-s_0}$, one has
  \begin{equation*}
    \begin{aligned}
      n=4:\quad&u_0=\Bigl(\sqrt{4+2\bigl(1-\frac{q}{\mu_3}\bigr)}-2\Bigr)^2,\qquad s_0=\frac{u_0}{1+u_0};\\
      n=5:\quad&\frac{K^2}{4+K^2}\leqslant s_0\leqslant\frac{K^2}{1+K^2},\qquad K:=\frac{27}{28}\bigl(1-\frac{q}{\mu_3}\bigr);\\
      n\geqslant6:\quad&u_0=\left(\frac{2^{\ast}}{2}\bigl(1-\frac{q}{\mu_3}\bigr)\right)^{\frac{2}{q-1}},\qquad s_0=\frac{u_0}{1+u_0}.
    \end{aligned}
  \end{equation*}
  These formulas follow by solving $\nu=0$ as an equation in $u$. For $n=4$, $\nu=0$ is the quadratic equation $u+4\sqrt u-2\bigl(1-\frac{q}{\mu_3}\bigr)=0$ in $\sqrt u$, whose positive root gives the formula. For $n\geqslant6$ the penalty is a single power $u^{\frac{q-1}{2}}$, which can be solved directly. For $n=5$, the zero solves $\sqrt{u_0}+u_0^{\frac{2}{3}}=K$, and $K\in(0,1)$ since $0<1-\frac{q}{\mu_3}<1$. Thus $u_0<1$, and $u^{\frac{2}{3}}\leqslant\sqrt u$ on $(0,1)$ gives $\sqrt{u_0}\leqslant K\leqslant2\sqrt{u_0}$, that is, $\frac{K^2}{4}\leqslant u_0\leqslant K^2$; this is the displayed bracket for $s_0=\frac{u_0}{1+u_0}$.

  \emph{Step 2: reduction to the nearest point.} Fix $s\in(0,1)$ and $f\in H^1(\mbBn)\setminus\{0\}$ with $0<\sigma(f)\leqslant s$.
  Both $\delta(\cdot)^2$ and $\dist(\cdot,\mcZ)^2$ are 2-homogeneous and $\sigma$ is 0-homogeneous, so we may normalize $\norm{f}_\lambda^2=S_\lambda^{\frac{n}{2}}$, which by Lemma~\ref{lem:ground-state-normalization} equals $\norm{\mcU}_\lambda^2=\int_{\mbBn}\mcU^{2^{\ast}}\diopv$. Since $0<\sigma(f)<1$, Lemma~\ref{lem:const-projection} gives, after an isometry and possibly replacing $f$ by $-f$, the decomposition $f=a\mcU+d_f\psi$, and \eqref{equ:const-modulation} gives
  \begin{equation*}
    a^2=1-\sigma(f)\in(0,1),\quad d_f^2=S_\lambda^{\frac{n}{2}}\sigma(f),\quad\text{so}\quad
    a\geqslant\sqrt{1-s},\quad d_f\leqslant S_\lambda^{\frac{n}{4}}\sqrt{s}\quad\text{when }\sigma(f)\leqslant s.
  \end{equation*}
  In particular $\frac{d_f}{aS_\lambda^{\frac{n}{4}}}=\bigl(\frac{\sigma(f)}{1-\sigma(f)}\bigr)^{\frac{1}{2}}\leqslant\sqrt u$ when $\sigma(f)\leqslant s$. We shall also use two facts on $\psi$: the first identity of \eqref{equ:const-spectral-moment} removes the linear term in the expansions below, and \eqref{equ:PS} applied to $\psi$ gives $\norm{\psi}_{L^{2^{\ast}}(\mbBn)}\leqslant S_\lambda^{-\frac{1}{2}}$.

  \emph{Step 3: the scheme of the expansion.} Since $\norm{f}_\lambda^2=a^2S_\lambda^{\frac{n}{2}}+d_f^2$ by \eqref{equ:const-modulation}, a lower bound for $\delta(f)^2$ amounts to an upper bound for $S_\lambda\norm{f}_{L^{2^{\ast}}(\mbBn)}^2$, the only non-quadratic term of the deficit. We obtain it in three moves. First, we bound $\norm{f}_{L^{2^{\ast}}(\mbBn)}^{2^{\ast}}=\int|a\mcU+d_f\psi|^{2^{\ast}}\diopv$ from above by expanding in the summands $a\mcU$ and $d_f\psi$, which are of very different size; the quality of the result depends on the regularity of $t\mapsto|1+t|^{2^{\ast}}$ at $0$, that is, on $n$. Since $\mcU>0$, Lemma~\ref{lem:taylor-remainder} with $p=2^{\ast}$ and $t=\frac{d_f\psi}{a\mcU}$, multiplied by $(a\mcU)^{2^{\ast}}$, gives pointwise bounds for $|f|^{2^{\ast}}$: case~\ref{item:remainder-small} applies when $2<2^{\ast}\leqslant3$, that is, $n\geqslant6$, and case~\ref{item:remainder-large} applies when $n=5$; for $n=4$ the expansion is an identity. Second, we pass from $\norm{f}_{L^{2^{\ast}}(\mbBn)}^{2^{\ast}}$ to $S_\lambda\norm{f}_{L^{2^{\ast}}(\mbBn)}^2$ by the concavity of $x\mapsto(1+x)^{\frac{2}{2^{\ast}}}$. Third, we subtract the result from $\norm{f}_\lambda^2$.

  \emph{Dimensions $n\geqslant6$} ($1<q\leqslant2$). Apply Lemma~\ref{lem:taylor-remainder}\ref{item:remainder-small} with $t=\frac{d_f\psi}{a\mcU}$, multiply by $(a\mcU)^{2^{\ast}}$, and integrate. The zeroth-order term gives $a^{2^{\ast}}S_\lambda^{\frac{n}{2}}$, the linear term vanishes, the quadratic term is bounded by the second part of \eqref{equ:const-spectral-moment}, and the remainder by $\norm{\psi}_{L^{2^{\ast}}(\mbBn)}\leqslant S_\lambda^{-\frac{1}{2}}$. We find
  \begin{equation*}
    \norm{f}_{L^{2^{\ast}}(\mbBn)}^{2^{\ast}}\leqslant a^{2^{\ast}}S_\lambda^{\frac{n}{2}}+\frac{2^{\ast}q}{2\mu_3}a^{q-1}d_f^2
    +S_\lambda^{-\frac{2^{\ast}}{2}}d_f^{2^{\ast}}.
  \end{equation*}
  Write the right-hand side as $a^{2^{\ast}}S_\lambda^{\frac{n}{2}}(1+x)$ with $x\geqslant0$. Since $0<\frac{2}{2^{\ast}}<1$, concavity gives $(1+x)^{\frac{2}{2^{\ast}}}\leqslant1+\frac{2}{2^{\ast}}x$. Multiplying by $S_\lambda$ and using $S_\lambda\bigl(a^{2^{\ast}}S_\lambda^{\frac{n}{2}}\bigr)^{\frac{2}{2^{\ast}}}=a^2S_\lambda^{\frac{n}{2}}$, we find
  \begin{equation*}
    S_\lambda\norm{f}_{L^{2^{\ast}}(\mbBn)}^2\leqslant a^2S_\lambda^{\frac{n}{2}}+\frac{q}{\mu_3}d_f^2
    +\frac{2}{2^{\ast}}a^{1-q}S_\lambda^{-\frac{2^{\ast}}{2}}d_f^{2^{\ast}}.
  \end{equation*}
  Subtracting this from $\norm{f}_\lambda^2=a^2S_\lambda^{\frac{n}{2}}+d_f^2$ gives $\delta(f)^2\geqslant[1-\frac{q}{\mu_3}-\frac{2}{2^{\ast}}a^{1-q}S_\lambda^{-\frac{2^{\ast}}{2}}d_f^{q-1}]d_f^2$. Since $a^2=1-\sigma$ and $d_f^2=S_\lambda^{\frac{n}{2}}\sigma$, the powers of $S_\lambda$ cancel and the penalty is $\frac{2}{2^{\ast}}\bigl(\frac{\sigma}{1-\sigma}\bigr)^{\frac{q-1}{2}}$. So, when $\sigma(f)\leqslant s$,
  \begin{equation*}
    \delta(f)^2\geqslant\left[1-\frac{q}{\mu_3}-\frac{2}{2^{\ast}}\left(\frac{\sigma(f)}{1-\sigma(f)}\right)^{\frac{q-1}{2}}\right]d_f^2
    \geqslant\left[1-\frac{q}{\mu_3}-\frac{2}{2^{\ast}}u^{\frac{q-1}{2}}\right]d_f^2,
  \end{equation*}
  which is the third branch of \eqref{equ:const-nu-explicit}.

  \emph{Dimensions $n=4,5$} ($2<q\leqslant3$). Here the expansion produces a cubic term, which we first bound. The exponents $\frac{2^{\ast}}{q-2}$ and $\frac{2^{\ast}}{3}$ are conjugate; they are $4$ and $\frac{4}{3}$ for $n=4$, and $10$ and $\frac{10}{9}$ for $n=5$. H\"older's inequality with these exponents and $\norm{\psi}_{L^{2^{\ast}}(\mbBn)}\leqslant S_\lambda^{-\frac{1}{2}}$ give, with Lemma~\ref{lem:ground-state-normalization}, $\int_{\mbBn}\mcU^{q-2}|\psi|^3\diopv\leqslant\bigl(\int_{\mbBn}\mcU^{2^{\ast}}\diopv\bigr)^{\frac{q-2}{2^{\ast}}}\norm{\psi}_{L^{2^{\ast}}(\mbBn)}^3\leqslant S_\lambda^{\frac{6-n}{4}}S_\lambda^{-\frac{3}{2}}=S_\lambda^{-\frac{n}{4}}$.

  \emph{Dimension $5$.} Apply Lemma~\ref{lem:taylor-remainder}\ref{item:remainder-large} in the same way, with $\kappa_{2^{\ast}}=1$, and bound the coefficient $1$ of $|t|^{2^{\ast}}$ by $C_{2^{\ast}}>1$. Integrating and using the cubic bound just proved, with $C_q:=C_{2^{\ast}}=\frac{(q+1)q(q-1)}{6}$, we find
  \begin{equation*}
    \norm{f}_{L^{2^{\ast}}(\mbBn)}^{2^{\ast}}\leqslant a^{2^{\ast}}S_\lambda^{\frac{n}{2}}+\frac{2^{\ast}q}{2\mu_3}a^{q-1}d_f^2
    +C_q\bigl(a^{q-2}S_\lambda^{-\frac{n}{4}}d_f^3+S_\lambda^{-\frac{2^{\ast}}{2}}d_f^{2^{\ast}}\bigr).
  \end{equation*}
  The concavity and normalization steps of the case $n\geqslant6$ then give, when $\sigma(f)\leqslant s$,
  \begin{equation*}
    \delta(f)^2\geqslant\left\{1-\frac{q}{\mu_3}-\frac{2C_q}{2^{\ast}}
    \left[u^{\frac{1}{2}}+u^{\frac{q-1}{2}}\right]\right\}d_f^2,
  \end{equation*}
  the penalty $u^{\frac{1}{2}}$ coming from the cubic term and $u^{\frac{q-1}{2}}$ from the $2^{\ast}$-power. For $n=5$, $q=\frac{7}{3}$ and $C_q=\frac{140}{81}$, so $\frac{2C_q}{2^{\ast}}=\frac{28}{27}$, and $\frac{q-1}{2}=\frac{2}{3}$. This is the second branch of \eqref{equ:const-nu-explicit}.

  \emph{Dimension $4$.}
  Here $q=3$, $2^{\ast}=4$ and $|f|^4=f^4$, so the expansion is an identity. Expanding $f=a\mcU+d_f\psi$, the linear term vanishes by \eqref{equ:const-spectral-moment}, and we bound the other terms by the second part of \eqref{equ:const-spectral-moment}, the cubic bound $\int_{\mbBn}\mcU|\psi|^3\diopv\leqslant S_\lambda^{-1}$ and $\int_{\mbBn}\psi^4\diopv\leqslant S_\lambda^{-2}$. We find
  \begin{equation*}
    \begin{aligned}
      \norm{f}_{L^{2^{\ast}}(\mbBn)}^{2^{\ast}}&=a^4S_\lambda^2+6a^2d_f^2\int_{\mbBn}\mcU^2\psi^2\diopv+4ad_f^3\int_{\mbBn}\mcU\psi^3\diopv+d_f^4\int_{\mbBn}\psi^4\diopv\\
      &\leqslant a^4S_\lambda^2+\frac{6}{\mu_3}a^2d_f^2+4aS_\lambda^{-1}d_f^3+S_\lambda^{-2}d_f^4.
    \end{aligned}
  \end{equation*}
  Write the right-hand side as $a^4S_\lambda^2(1+x)$ with $x\geqslant0$. Concavity now reads $\sqrt{1+x}\leqslant1+\frac{x}{2}$, which gives
  \begin{equation*}
    S_\lambda\norm{f}_{L^4(\mbBn)}^2\leqslant a^2S_\lambda^2+\frac{3}{\mu_3}d_f^2+2a^{-1}S_\lambda^{-1}d_f^3+\frac{1}{2}a^{-2}S_\lambda^{-2}d_f^4 .
  \end{equation*}
  Subtracting this from $\norm{f}_\lambda^2=a^2S_\lambda^2+d_f^2$ and using $\frac{d_f}{aS_\lambda}\leqslant\sqrt u$, we obtain $\delta(f)^2\geqslant[1-\frac{q}{\mu_3}-2\sqrt u-\frac{1}{2}u]d_f^2$, the first branch of \eqref{equ:const-nu-explicit}.

  In all dimensions we have thus proved \eqref{equ:local-prop-inequality}. Finally, taking $s=\sigma(f)$ in \eqref{equ:local-prop-inequality} gives the last assertion. This finishes the proof of Proposition~\ref{prop:local-lower-function}.
\end{proof}

\begin{proof}[Proof of Theorem~\ref{thm:crit-local}]
  We first deduce the lower bound $c_{\mrloc}\geqslant1-\frac{q}{\mu_3}$ from Proposition~\ref{prop:local-lower-function}, and then prove the matching upper bound by testing along an eigendirection of $\mu_3$, for which we compute the distance to $\mcZ$ exactly.

  \emph{Step 1: the lower bound.} Since $\mcT_s$ increases with $s$, the infimum in the definition of $c_{\mrloc}$ is nonincreasing in $s$, so the limit defining $c_{\mrloc}$ exists in $[0,\infty]$. By Proposition~\ref{prop:local-lower-function}, this infimum is at least $\nu(s)$ for every $s\in(0,1)$. Letting $s\downarrow0$ and using \eqref{equ:const-nu-limits}, we obtain $c_{\mrloc}\geqslant1-\frac{q}{\mu_3}$.

  \emph{Step 2: the test family.} Let $\Psi$ be a real eigenfunction of \eqref{equ:const-linearized} with eigenvalue $\mu_3$, normalized by $\norm{\Psi}_\lambda=1$. Since $\mu_3\notin\{1,q\}$ and eigenfunctions with different eigenvalues are $\lag\cdot,\cdot\rag_\lambda$-orthogonal, $\Psi$ is orthogonal to $\MaOpspan\{\mcU,\Phi_1,\ldots,\Phi_n\}$. Testing the weak form of \eqref{equ:const-linearized} for $\Psi$ with $\mcU$ and with $\Psi$, we find
  \begin{equation}\label{equ:const-eigen-moments}
    \int_{\mbBn}\mcU^{q}\Psi\diopv=0,\qquad\int_{\mbBn}\mcU^{q-1}\Psi^2\diopv=\frac{1}{\mu_3}.
  \end{equation}
  Set $f_\ve:=\mcU+\ve\Psi$. We claim that for all sufficiently small $|\ve|\neq0$,
  \begin{equation}\label{equ:const-exact-distance}
    \dist(f_\ve,\mcZ)^2=\ve^2.
  \end{equation}

  \emph{Step 3: proof of \eqref{equ:const-exact-distance}.} We first show that $(1,o)$ is a strict local minimum of the distance function, then that the nearest point of $f_\ve$ is close to $\mcU$, and finally combine the two. By Lemma~\ref{lem:pairing-center-derivative}\ref{item:orbit-map-C2} and Lemma~\ref{lem:linearized-eigenspaces}\ref{item:eigen-tangent}, $(c,\xi)\mapsto c\mcU_\xi$ is a $C^2$ local parametrization of $\mcZ$ near $\mcU$ with tangent space $\MaOpspan\{\mcU,\Phi_1,\ldots,\Phi_n\}$ at $(1,o)$. The function $E(\ve,c,\xi):=\norm{f_\ve-c\mcU_\xi}_\lambda^2$ has the critical point $(c,\xi)=(1,o)$ for every $\ve$, by the orthogonality of $\Psi$. At $\ve=0$ the Hessian of $E$ in $(c,\xi)$ at $(1,o)$ is twice the Gram matrix of $\mcU,\Phi_1,\ldots,\Phi_n$, which is positive definite since these functions are linearly independent by Lemma~\ref{lem:linearized-eigenspaces}\ref{item:eigen-tangent}. The gradient of $E$ in $(c,\xi)$ is $C^1$ in $(\ve,c,\xi)$. Hence, by the implicit function theorem, there are a neighborhood $N$ of $(1,o)$ and $\ve_0>0$ such that, for $|\ve|<\ve_0$, $(1,o)$ is the unique critical point of $E(\ve,\cdot,\cdot)$ in $N$ and a strict local minimum there.

  Next, the competitor $\mcU$ gives $\dist(f_\ve,\mcZ)\leqslant\norm{\ve\Psi}_\lambda=|\ve|$, so $\sigma(f_\ve)<1$ for small $\ve$, and Lemma~\ref{lem:const-projection} provides a nearest point $c_\ve\mcU_{\xi_\ve}$ with $\norm{f_\ve-c_\ve\mcU_{\xi_\ve}}_\lambda\leqslant|\ve|$; hence $c_\ve\mcU_{\xi_\ve}\to\mcU$ in $H^1(\mbBn)$ as $\ve\to0$. The centers $\xi_\ve$ cannot escape, since translates converge weakly to $0$ at infinity (Lemma~\ref{lem:overlap}). So every limit point $(c,\xi)$ of $(c_\ve,\xi_\ve)$ as $\ve\to0$ satisfies $c\mcU_\xi=\mcU$, and the positivity of $\mcU$, its strict radial monotonicity (Lemma~\ref{lem:ground-state}) and the uniqueness of its maximum point force $(c,\xi)=(1,o)$. Hence $(c_\ve,\xi_\ve)\in N$ for small $|\ve|$. Being a minimizer of $E(\ve,\cdot,\cdot)$, it is a critical point, so $(c_\ve,\xi_\ve)=(1,o)$, and $\dist(f_\ve,\mcZ)^2=\norm{f_\ve-\mcU}_\lambda^2=\ve^2$, which is \eqref{equ:const-exact-distance}.

  \emph{Step 4: expansion of the deficit.} Since $2^{\ast}>2$, the integral form of the Taylor remainder gives
  \begin{equation*}
    |\mcU+\ve\Psi|^{2^{\ast}}=\mcU^{2^{\ast}}+2^{\ast}\mcU^{q}\ve\Psi+2^{\ast}q\,\ve^2\Psi^2\int_0^1(1-t)|\mcU+t\ve\Psi|^{q-1}\diopt.
  \end{equation*}
  For $|\ve|\leqslant1$, the integrand of the last term is dominated by $C(\mcU^{q-1}\Psi^2+|\Psi|^{2^{\ast}})\in L^1(\mbBn)$, by H\"older's inequality and $H^1\hookrightarrow L^{2^{\ast}}$, and it converges pointwise to $\frac{1}{2}\mcU^{q-1}\Psi^2$ as $\ve\to0$. By dominated convergence, Lemma~\ref{lem:ground-state-normalization} and \eqref{equ:const-eigen-moments}, we find
  \begin{equation}\label{equ:local-eigenfunction-expansion}
    \norm{f_\ve}_{L^{2^{\ast}}(\mbBn)}^{2^{\ast}}=S_\lambda^{\frac{n}{2}}+\frac{2^{\ast}q}{2\mu_3}\ve^2+o(\ve^2),\qquad
    \norm{f_\ve}_\lambda^2=S_\lambda^{\frac{n}{2}}+\ve^2,
  \end{equation}
  where the second identity follows from the orthogonality of $\Psi$ and $\mcU$. Expanding the outer power $x\mapsto x^{\frac{2}{2^{\ast}}}$ at $S_\lambda^{\frac{n}{2}}$, and using $\frac{2^{\ast}q}{2\mu_3}\cdot\frac{2}{2^{\ast}}=\frac{q}{\mu_3}$, we obtain
  \begin{equation*}
    S_\lambda\norm{f_\ve}_{L^{2^{\ast}}(\mbBn)}^2=S_\lambda^{\frac{n}{2}}+\frac{q}{\mu_3}\ve^2+o(\ve^2).
  \end{equation*}
  Subtracting this from $\norm{f_\ve}_\lambda^2$ and dividing by \eqref{equ:const-exact-distance}, we find
  \begin{equation}\label{equ:local-eigenfunction-stability-quotient}
    \frac{\delta(f_\ve)^2}{\dist(f_\ve,\mcZ)^2}
    =\left(1-\frac{q}{\mu_3}\right)+o(1)\quad(\ve\to0).
  \end{equation}

  \emph{Step 5: the upper bound.} By \eqref{equ:const-exact-distance} and \eqref{equ:local-eigenfunction-expansion}, $\sigma(f_\ve)=\frac{\ve^2}{S_\lambda^{\frac{n}{2}}+\ve^2}$. Fix $s\in(0,1)$. If $\ve\neq0$ is small and $\ve^2\leqslant S_\lambda^{\frac{n}{2}}s$, then $\sigma(f_\ve)\leqslant s$, so $f_\ve\in\mcT_s$. Letting $\ve\to0$ in \eqref{equ:local-eigenfunction-stability-quotient} gives
  \begin{equation*}
    \inf_{f\in\mcT_s}\frac{\delta(f)^2}{\dist(f,\mcZ)^2}\leqslant1-\frac{q}{\mu_3}.
  \end{equation*}
  Letting $s\downarrow0$ gives $c_{\mrloc}\leqslant1-\frac{q}{\mu_3}$. Together with Step~1, this proves $c_{\mrloc}=1-\frac{q}{\mu_3}$. This finishes the proof of Theorem~\ref{thm:crit-local}.
\end{proof}
\section{Global lower bound}\label{sec:crit-global-bounds-proof}
\begingroup
Proposition~\ref{prop:local-lower-function} controls the stability quotient on the tubes $\mcT_s$ with $s<s_0$. This section extends the control to all of $H^1(\mbBn)\setminus(\mcZ\cup\{0\})$ and proves the following weak form of the lower bound in Theorem~\ref{thm:crit-global-bounds}, which still contains the concentration constant $c_{\mathrm{conc}}$.
\begin{proposition}\label{prop:weak-global-lower-bound}
  Every $f\in H^1(\mbBn)\setminus(\mcZ\cup\{0\})$ satisfies
  \begin{equation}\label{equ:const-global-quotient}
    \delta(f)^2\geqslant\min\Bigl\{\sup_{0<s<1}s\nu(s),\,c_{\mathrm{conc}},\,\frac{1}{2}c_{\mathrm{dich}}\Bigr\}\dist(f,\mcZ)^2,
  \end{equation}
  and the minimum in \eqref{equ:const-global-quotient} is positive. In particular $c_{\mathrm{BE}}$ is at least this minimum.
\end{proposition}
The proof is completed at the end of Section~\ref{subsec:crit-const-global}. In Section~\ref{subsec:proof-global} we combine Proposition~\ref{prop:weak-global-lower-bound} with the strict upper bounds of Section~\ref{sec:strict-upper} and remove $c_{\mathrm{conc}}$ from the minimum.

\emph{Normalization and the elementary bound.} For $t\neq0$ we have $\delta(tf)^2=t^2\delta(f)^2$ and $\dist(tf,\mcZ)=|t|\dist(f,\mcZ)$, the latter because $z\mapsto tz$ is a bijection of $\mcZ$. Hence the stability quotient and $\sigma$ do not change under $f\mapsto tf$, and it suffices to consider functions in the \emph{unit-mass constraint}
\begin{equation}\label{equ:const-unit-mass}
  \mcM:=\{f\in H^1(\mbBn):\norm{f}_{L^{2^{\ast}}(\mbBn)}=1\}.
\end{equation}
On $\mcM$ the deficit, and $\sigma$ times the stability quotient, are determined by the energy $\norm{f}_\lambda^2$.
\begin{lemma}\label{lem:const-trivial-bound}
  Let $f\in\mcM\setminus\mcZ$ and $E:=\norm{f}_\lambda^2$. Then $\dist(f,\mcZ)>0$, $\delta(f)^2=E-S_\lambda$, and
  \begin{equation}\label{equ:const-trivial-bound}
    \sigma(f)\frac{\delta(f)^2}{\dist(f,\mcZ)^2}=1-\frac{S_\lambda}{E},\qquad
    \frac{\delta(f)^2}{\dist(f,\mcZ)^2}\geqslant1-\frac{S_\lambda}{E}.
  \end{equation}
\end{lemma}
\begin{proof}
  The distance is positive by Lemma~\ref{lem:dist-positive}, and $\delta(f)^2=E-S_\lambda$ by the normalization. By the definition \eqref{equ:deficit-sigma} of $\sigma$, $\sigma(f)\frac{\delta(f)^2}{\dist(f,\mcZ)^2}=\frac{\delta(f)^2}{E}=1-\frac{S_\lambda}{E}$. The inequality follows, since $\sigma(f)\leqslant1$ by Lemma~\ref{lem:dist-formula}. This finishes the proof of Lemma~\ref{lem:const-trivial-bound}.
\end{proof}\endgroup\par
\begingroup
We set the \emph{threshold}
\begin{equation*}
  \wtilde{S}:=\min\{2^{\frac{2}{n}}S_\lambda,S_n\}>S_\lambda,
\end{equation*}
where the inequality holds because $2^{\frac{2}{n}}>1$ and $S_\lambda<S_n$. We split $\mcM\setminus\mcZ$ into the \emph{high-energy region} $\{\norm{f}_\lambda^2\geqslant\wtilde{S}\}$ and the \emph{low-energy region} $\{\norm{f}_\lambda^2<\wtilde{S}\}$.

The threshold $\wtilde{S}$ is the optimal choice. In the high-energy region a higher threshold is better, because Lemma~\ref{lem:const-trivial-bound} bounds the stability quotient from below by $1-\frac{S_\lambda}{\norm{f}_\lambda^2}$, which increases with the energy $\norm{f}_\lambda^2$, so the higher the energy at which the region starts, the better the bound on it. In the low-energy region, however, a lower threshold is better, because there we need Palais--Smale sequences with energy below the threshold to converge, up to isometries, to a single ground state (Lemma~\ref{lem:psdec-subthreshold}). A sequence can lose compactness by concentrating at a point, which costs energy at least $S_n$, or by splitting into two pieces, which costs energy at least $2^{\frac{2}{n}}S_\lambda$. So compactness holds below $\wtilde{S}$. It fails at $\wtilde{S}$, as the concentration and dichotomy families of Section~\ref{sec:strict-upper} suggest, so the threshold cannot be raised, although we do not need this fact in this section. Balancing the two requirements gives the threshold $\wtilde{S}$. At this threshold Proposition~\ref{prop:const-high-energy} gives $1-\frac{S_\lambda}{\wtilde{S}}=\min\{c_{\mathrm{conc}},\frac{1}{2}c_{\mathrm{dich}}\}$, which is where these two constants in Proposition~\ref{prop:weak-global-lower-bound} come from. How they compare with the upper bounds of Section~\ref{sec:strict-upper} is explained in Remark~\ref{rem:channels}.

\emph{Plan of the section.} Section~\ref{subsec:crit-const-threshold} treats the high-energy region, where Lemma~\ref{lem:const-trivial-bound} gives the bound at once, and Section~\ref{subsec:low-energy} treats the low-energy region; Proposition~\ref{prop:weak-global-lower-bound} is proved at the end of Section~\ref{subsec:low-energy}.\endgroup\par

\subsection{The high-energy region}\label{subsec:crit-const-threshold}
\begingroup
\begin{proposition}\label{prop:const-high-energy}
  Every $f\in\mcM\setminus\mcZ$ with $\norm{f}_\lambda^2\geqslant\wtilde{S}$ satisfies
  \begin{equation*}
    \frac{\delta(f)^2}{\dist(f,\mcZ)^2}\geqslant1-\frac{S_\lambda}{\norm{f}_\lambda^2}\geqslant1-\frac{S_\lambda}{\wtilde{S}}=\min\Bigl\{c_{\mathrm{conc}},\frac{1}{2}c_{\mathrm{dich}}\Bigr\}>0.
  \end{equation*}
\end{proposition}
\begin{proof}
  The first inequality is \eqref{equ:const-trivial-bound}, and the second holds because $t\mapsto1-\frac{S_\lambda}{t}$ is increasing. For positive $a,b$, $1-\frac{S_\lambda}{\min\{a,b\}}=\min\{1-\frac{S_\lambda}{a},1-\frac{S_\lambda}{b}\}$. With $a=2^{\frac{2}{n}}S_\lambda$ and $b=S_n$ this gives $1-\frac{S_\lambda}{\wtilde{S}}=\min\{1-2^{-\frac{2}{n}},c_{\mathrm{conc}}\}$, and $1-2^{-\frac{2}{n}}=\frac{1}{2}(2-2^{\frac{n-2}{n}})=\frac{1}{2}c_{\mathrm{dich}}$. Both entries are positive, because $2^{\frac{2}{n}}>1$ and $S_\lambda<S_n$. This finishes the proof of Proposition~\ref{prop:const-high-energy}.
\end{proof}\endgroup\par

\subsection{The low-energy region}\label{subsec:low-energy}
\begingroup
In this subsection we prove that every $f\in\mcM\setminus\mcZ$ with $\norm{f}_\lambda^2<\wtilde{S}$ satisfies $\frac{\delta(f)^2}{\dist(f,\mcZ)^2}\geqslant\sup_{0<s<1}s\nu(s)$ (Proposition~\ref{prop:const-transfer}). Such an $f$ need not be close to $\mcZ$, so the local estimate does not apply to it directly. Instead we follow the constrained energy flow from $f$. The energy does not increase along the flow, so by Lemma~\ref{lem:const-trivial-bound} neither does the deficit or $\sigma$ times the stability quotient. The argument has three steps. Section~\ref{subsec:flow} constructs the flow. Section~\ref{subsec:psdec} proves compactness below the threshold and deduces that the trajectory from $f$ eventually enters every tube $\mcT_s$. Section~\ref{subsec:crit-const-global} applies Proposition~\ref{prop:local-lower-function} at the first entrance into $\mcT_s$ and transfers the bound back to $f$.\endgroup\par

\subsubsection{An energy gradient flow on the constraint}\label{subsec:flow}
We use the negative gradient flow of the energy on $\mcM$. It keeps the mass equal to $1$ and does not increase the energy, which is all that the transfer in Section~\ref{subsec:crit-const-global} needs, and it acts on sign-changing data.
Equip the unit-mass constraint $\mcM$ of \eqref{equ:const-unit-mass} with the metric induced by $\lag\cdot,\cdot\rag_\lambda$.
Write $\norm{h}_{H^{-1}(\mbBn)}^2:=\lag h,(-\Delta_{\mbBn}-\lambda)^{-1}h\rag$ for the squared dual norm, where $H^{-1}(\mbBn)$ is the dual of $H^1(\mbBn)$ and the brackets denote the $H^{-1}$--$H^1$ pairing.
Thus $\norm{(-\Delta_{\mbBn}-\lambda)u}_{H^{-1}(\mbBn)}=\norm{u}_\lambda$. The dual of \eqref{equ:PS} gives
$\norm{h}_{H^{-1}(\mbBn)}\leqslant S_\lambda^{-\frac{1}{2}}\norm{h}_{L^{(2^\ast)'}(\mbBn)}$.

The \emph{constrained energy flow} is
\begin{equation}\label{equ:intro-projected-flow}
  \partial_\tau f_\tau=-2f_\tau+
  \frac{2\norm{f_\tau}_{L^{2^{\ast}}(\mbBn)}^{2^{\ast}}}{\bignorm{|f_\tau|^{q-1}f_\tau}_{H^{-1}(\mbBn)}^2}
  (-\Delta_{\mbBn}-\lambda)^{-1}(|f_\tau|^{q-1}f_\tau),
\end{equation}
defined as long as $f_\tau\neq0$.
The gradient of $f\mapsto\norm{f}_\lambda^2$ in this Hilbert space is $2f$, and the gradient of
$f\mapsto\norm{f}_{L^{2^{\ast}}(\mbBn)}^{2^{\ast}}$ is $2^\ast(-\Delta_{\mbBn}-\lambda)^{-1}(|f|^{q-1}f)$.
Their inner product is $2\cdot2^\ast\norm{f}_{L^{2^{\ast}}(\mbBn)}^{2^{\ast}}$.
Consequently the right-hand side of \eqref{equ:intro-projected-flow} is the
negative orthogonal projection of $2f$ onto the tangent space $\{h:\int_{\mbBn}|f|^{q-1}fh\diopv=0\}$.

\begin{proposition}\label{prop:flow-properties}
  For every $f_0\in\mcM$, equation
  \eqref{equ:intro-projected-flow} has a unique solution
  $f\in C^1([0,\infty);H^1(\mbBn))$ with $f|_{\tau=0}=f_0$.
  For every $\tau\geqslant0$,
  \begin{equation}\label{equ:flow-deficit-mass-one}
    \norm{f_\tau}_{L^{2^{\ast}}(\mbBn)}=1,\qquad
    \delta(f_\tau)^2=\norm{f_\tau}_\lambda^2-S_\lambda.
  \end{equation}
  For $0\leqslant a\leqslant b$,
  \begin{equation}\label{equ:flow-dissipation-normalized}
    \norm{f_b}_\lambda^2+\int_a^b\norm{\partial_\tau f_\tau}_\lambda^2\dioptau=\norm{f_a}_\lambda^2.
  \end{equation}
  Moreover,
  \begin{equation}\label{equ:flow-full-residual}
    \norm{\partial_\tau f_\tau}_\lambda\longrightarrow0,\qquad
    \bignorm{(-\Delta_{\mbBn}-\lambda)f_\tau-\norm{f_\tau}_\lambda^2|f_\tau|^{q-1}f_\tau}_{H^{-1}(\mbBn)}
    \longrightarrow0
    \quad(\tau\to\infty).
  \end{equation}
\end{proposition}
\begin{proof}
  We first construct a unique local solution, then show that the flow preserves the mass and dissipates the energy, then use the resulting bounds to show that the solution is global, and finally prove the two limits in \eqref{equ:flow-full-residual}.

  \emph{Step 1: local existence.} For real $a,b$, the mean-value formula gives
  \begin{equation*}
    \bigl||a|^{q-1}a-|b|^{q-1}b\bigr|
    \leqslant q(|a|+|b|)^{q-1}|a-b|.
  \end{equation*}
  H\"older's inequality therefore shows that
  $f\mapsto |f|^{q-1}f$ is Lipschitz from bounded subsets of
  $L^{2^\ast}$ into $L^{(2^\ast)'}$.
  Composing with the Sobolev embedding and $(-\Delta_{\mbBn}-\lambda)^{-1}$ gives a locally
  Lipschitz map on $H^1(\mbBn)$.
  The denominator in \eqref{equ:intro-projected-flow} is positive when $f\neq0$,
  because
  $\lag |f|^{q-1}f,f\rag=\norm{f}_{L^{2^{\ast}}(\mbBn)}^{2^{\ast}}>0$.
  All terms on the right-hand side are thus locally Lipschitz on
  $H^1(\mbBn)\setminus\{0\}$. The integral equation is a contraction on a sufficiently
  short time interval in a closed ball about $f_0$. It gives a unique
  local $C^1$ solution.

  \emph{Step 2: mass and energy.} Pairing \eqref{equ:intro-projected-flow} with the gradient of $f\mapsto\norm{f}_{L^{2^{\ast}}(\mbBn)}^{2^{\ast}}$ gives
  $\frac{\mrd}{\mrd\tau}\norm{f_\tau}_{L^{2^{\ast}}(\mbBn)}^{2^{\ast}}=0$.
  Pairing it with $2f_\tau$ gives
  \begin{equation*}
    \frac{\mrd}{\mrd\tau}\norm{f_\tau}_\lambda^2
    =-4\norm{f_\tau}_\lambda^2
    +\frac{4\norm{f_\tau}_{L^{2^{\ast}}(\mbBn)}^{2\cdot2^{\ast}}}{\bignorm{|f_\tau|^{q-1}f_\tau}_{H^{-1}(\mbBn)}^2}
    =-\norm{\partial_\tau f_\tau}_\lambda^2.
  \end{equation*}
  Here the last equality also follows by taking the norm of the
  orthogonal projection.
  This proves \eqref{equ:flow-deficit-mass-one} and
  \eqref{equ:flow-dissipation-normalized} on the existence interval.
  In particular,
  $S_\lambda\leqslant \norm{f_\tau}_\lambda^2\leqslant \norm{f_0}_\lambda^2$.

  \emph{Step 3: global existence.} By Step~2 the trajectory stays in the set of $f\in\mcM$ with $\norm{f}_\lambda^2\leqslant \norm{f_0}_\lambda^2$. For such $f$, testing the
  nonlinearity against $f$ and using the dual Sobolev bound yields
  \begin{equation}\label{equ:flow-normal-bounds}
    \frac{1}{\norm{f_0}_\lambda}
    \leqslant\frac{1}{\norm{f}_\lambda}
    \leqslant\bignorm{|f|^{q-1}f}_{H^{-1}(\mbBn)}
    \leqslant\frac{1}{\sqrt{S_\lambda}}.
  \end{equation}
  The nonlinearity has a uniform Lipschitz constant on this bounded set.
  The difference of its squared dual norms is bounded by its norm
  difference times the sum of the norms.
  By \eqref{equ:flow-normal-bounds}, the denominator is bounded
  below by $\frac{1}{\norm{f_0}_\lambda^2}$.
  These bounds give a single Lipschitz constant for the right-hand side
  of \eqref{equ:intro-projected-flow} on this set.
  Also $\norm{\partial_\tau f_\tau}_\lambda\leqslant2\norm{f_0}_\lambda$.
  If the maximal existence time were finite, this speed bound would
  make $f_\tau$ Cauchy at that time.
  Its limit would lie in $\mcM$ by continuity of the $L^{2^{\ast}}$ norm on $H^1(\mbBn)$.
  The local construction at that limit would extend the solution,
  a contradiction. The solution is global.

  \emph{Step 4: the first limit.} The same Lipschitz constant and speed bound show that
  $\tau\mapsto\partial_\tau f_\tau$ is uniformly Lipschitz in $H^1$.
  Its squared norm is integrable by
  \eqref{equ:flow-dissipation-normalized}.
  If its norm did not tend to zero, there would be $\ve>0$
  and times tending to infinity at which that norm is at least
  $\ve$. Uniform continuity would give disjoint intervals of
  one fixed positive length on which the norm is at least
  $\frac{\ve}{2}$.
  Their contributions to the integral would have infinite sum.
  This proves the first limit in \eqref{equ:flow-full-residual}.

  \emph{Step 5: the second limit.} On $\mcM$, the projection identity is
  \begin{equation*}
    \norm{f_\tau}_\lambda^2
    -\frac{1}{\bignorm{|f_\tau|^{q-1}f_\tau}_{H^{-1}(\mbBn)}^2}
    =\frac{1}{4}\norm{\partial_\tau f_\tau}_\lambda^2.
  \end{equation*}
  Applying $-\Delta_{\mbBn}-\lambda$ to \eqref{equ:intro-projected-flow} and using this identity gives
  \begin{equation*}
    (-\Delta_{\mbBn}-\lambda)f_\tau-\norm{f_\tau}_\lambda^2|f_\tau|^{q-1}f_\tau
    =-\frac{1}{2}(-\Delta_{\mbBn}-\lambda)\partial_\tau f_\tau
    -\frac{1}{4}\norm{\partial_\tau f_\tau}_\lambda^2|f_\tau|^{q-1}f_\tau.
  \end{equation*}
  By \eqref{equ:flow-normal-bounds}, its dual norm is bounded by
  $\frac{1}{2}\norm{\partial_\tau f_\tau}_\lambda
  +\frac{1}{4\sqrt{S_\lambda}}\norm{\partial_\tau f_\tau}_\lambda^2$,
  which tends to $0$ by Step~4. This finishes the proof of Proposition~\ref{prop:flow-properties}.
\end{proof}

\subsubsection{Subthreshold compactness and low-energy orbits}\label{subsec:psdec}
We show that every trajectory starting in the low-energy region approaches $\mcZ$ (Proposition~\ref{prop:low-energy-orbit}).
For $u\in H^1(\mbBn)$, let
\begin{equation*}
  \mcI(u):=\frac{1}{2}\norm{u}_\lambda^2-\frac{1}{2^\ast}\norm{u}_{L^{2^{\ast}}(\mbBn)}^{2^{\ast}}.
\end{equation*}
A \emph{Palais--Smale sequence} at level $c$ satisfies
$\mcI(u_k)\to c$ and
$\mcI'(u_k)=(-\Delta_{\mbBn}-\lambda)u_k-|u_k|^{q-1}u_k\to0$ in $H^{-1}(\mbBn)$.
For an isometry $\varphi$ of $\mbBn$, the map $u\mapsto u\circ\varphi$ is an isometry for $\norm{\cdot}_\lambda$ and preserves the $L^{2^{\ast}}$ norm and
$\mcI$. Its dual action is also isometric, so it preserves the
dual norm of the Palais--Smale residual.

\begin{lemma}\label{lem:psdec-first}
  A Palais--Smale sequence $\{u_k\}$ at level $c$ is bounded in $H^1(\mbBn)$,
  and
  \begin{equation}\label{equ:psdec-ps-mass}
    \norm{u_k}_\lambda^2\longrightarrow nc,\qquad \norm{u_k}_{L^{2^{\ast}}(\mbBn)}^{2^{\ast}}\longrightarrow nc.
  \end{equation}
  After passage to a subsequence, $u_k\rightharpoonup u_\infty$,
  $\mcI'(u_\infty)=0$, and $\{u_k-u_\infty\}$ is a
  Palais--Smale sequence at level $c-\mcI(u_\infty)$.
  Both $\norm{\cdot}_\lambda^2$ and $\norm{\cdot}_{L^{2^{\ast}}(\mbBn)}^{2^{\ast}}$ split between $u_\infty$ and $u_k-u_\infty$, up to $o(1)$.
  Every Palais--Smale level is either $0$ or at least
  $\frac{1}{n}S_\lambda^{\frac{n}{2}}$, and a sequence at level $0$
  converges strongly to zero.
\end{lemma}
\begin{proof}
  We first prove boundedness and \eqref{equ:psdec-ps-mass}, then show that the weak limit is a critical point, then show that the remainder is a Palais--Smale sequence with the stated splittings, and finally prove the alternative for the levels.

  \emph{Step 1: boundedness.} The identity
  \begin{equation*}
    \frac{1}{n}\norm{u_k}_\lambda^2
    =\mcI(u_k)-\frac{1}{2^\ast}
    \lag\mcI'(u_k),u_k\rag
  \end{equation*}
  bounds $\frac{1}{n}\norm{u_k}_\lambda^2$ by $|\mcI(u_k)|+\frac{1}{2^\ast}\norm{\mcI'(u_k)}_{H^{-1}(\mbBn)}\norm{u_k}_\lambda$.
  The action is bounded and the residual tends to zero, so $\norm{u_k}_\lambda$ is bounded.
  The two $H^1$ norms are equivalent, which gives boundedness in $H^1(\mbBn)$.
  Testing $\mcI'(u_k)$ against $u_k$ then gives
  $\norm{u_k}_\lambda^2-\norm{u_k}_{L^{2^{\ast}}(\mbBn)}^{2^{\ast}}\to0$ and proves \eqref{equ:psdec-ps-mass}.

  \emph{Step 2: the weak limit.} Local Rellich compactness and a diagonal subsequence give a weak
  limit $u_\infty$, almost everywhere convergence, and local strong
  $L^s$ convergence for $s<2^\ast$.
  The nonlinear terms are bounded in $L^{(2^\ast)'}$ and converge
  almost everywhere. Their local weak limit is
  $|u_\infty|^{q-1}u_\infty$.
  Testing against compactly supported smooth functions therefore gives
  $\mcI'(u_\infty)=0$, first on these tests and then by density.

  \emph{Step 3: the remainder.} Hilbert space orthogonality gives the splitting of $\norm{\cdot}_\lambda^2$.
  The Brezis--Lieb lemma \cite[Theorem~1]{BrezisLieb1983} gives the splitting
  of $\norm{\cdot}_{L^{2^{\ast}}(\mbBn)}^{2^{\ast}}$.
  For the residual, let $F(t)=|t|^{q-1}t$. The mean-value and Young inequalities give, for each $\ve>0$,
  \begin{equation*}
    |F(a+b)-F(a)-F(b)|
    \leqslant\ve|a|^{q}+C_\ve|b|^{q}.
  \end{equation*}
  Set $D_k:=F(u_k)-F(u_k-u_\infty)-F(u_\infty)$, so that $\mcI'(u_k-u_\infty)=\mcI'(u_k)-\mcI'(u_\infty)+D_k$.
  Then $D_k\to0$ almost everywhere and
  \begin{equation*}
    0\leqslant
    (|D_k|-\ve|u_k-u_\infty|^{q})_+
    \leqslant C_\ve|u_\infty|^{q}.
  \end{equation*}
  Since $q(2^\ast)'=2^\ast$,
  dominated convergence sends the middle expression to zero in $L^{(2^\ast)'}$.
  Hence
  $\limsup_k\norm{D_k}_{L^{(2^\ast)'}(\mbBn)}
  \leqslant\ve\sup_k\norm{u_k-u_\infty}_{L^{2^\ast}(\mbBn)}^{q}$.
  Letting $\ve\downarrow0$ gives $D_k\to0$ in $L^{(2^\ast)'}$.
  Since $\mcI'(u_k)\to0$ and $\mcI'(u_\infty)=0$, the dual Sobolev embedding
  gives $\mcI'(u_k-u_\infty)\to0$.
  The action consequently splits as well.

  \emph{Step 4: the levels.} By \eqref{equ:psdec-ps-mass}, $c\geqslant0$.
  If $c>0$, apply \eqref{equ:PS} to $u_k$ and pass to the limit to obtain
  $nc\geqslant S_\lambda(nc)^{\frac{2}{2^\ast}}$.
  Thus $nc\geqslant S_\lambda^{\frac{n}{2}}$.
  If $c=0$, \eqref{equ:psdec-ps-mass} gives $\norm{u_k}_\lambda\to0$. This finishes the proof of Lemma~\ref{lem:psdec-first}.
\end{proof}

\begingroup
Below the threshold, a Palais--Smale sequence could a priori fail to converge to a single ground state in three ways: its mass could concentrate at a point, its weak limit could change sign, or a second nonzero piece could split off. In the proof of Lemma~\ref{lem:psdec-subthreshold}, Step~2 excludes the first by means of the entry $S_n$ of $\wtilde{S}$, and Steps~3 and~5 exclude the other two by means of the entry $2^{\frac{2}{n}}S_\lambda$.\endgroup\par

\begin{lemma}\label{lem:psdec-subthreshold}
  Let $\{u_k\}$ be a Palais--Smale sequence for
  $\mcI$ at a level
  \begin{equation}\label{equ:psdec-subthreshold-level}
    c<\frac{1}{n}\min\{2S_\lambda^{\frac{n}{2}},S_n^{\frac{n}{2}}\}
    =\frac{1}{n}\wtilde{S}^{\frac{n}{2}}.
  \end{equation}
  After passage to a subsequence, one of the following holds.
  Either $c=0$ and $u_k\to0$ in $H^1(\mbBn)$, or
  $c=\frac{1}{n}S_\lambda^{\frac{n}{2}}$ and there are isometries $\varphi_k$
  and a function $Z=\mcU_\xi$ or $Z=-\mcU_\xi$ such that
  $u_k\circ\varphi_k\to Z$ strongly in $H^1(\mbBn)$.
\end{lemma}
\begin{proof}
  Lemma~\ref{lem:psdec-first} gives $c\geqslant0$ and handles $c=0$.
  Assume $c>0$. We first locate a fixed amount of mass in a geodesic ball and move its center to the origin, then show that the weak limit $Z$ is nonzero, then show that $Z$ has one sign and identify it, and finally show that the remainder vanishes.

  \emph{Step 1: covering and nonvanishing.} Choose a covering of $\mbBn$ by unit geodesic balls
  $B_{\mbH}(y_\alpha,1)$ whose doubled balls have uniformly bounded
  overlap.
  For example, a maximal separated net gives such a covering.
  The smaller disjoint balls with centers in any one doubled ball lie
  in a fixed larger ball, and hyperbolic volume comparison bounds their
  number.
  The Sobolev inequality on a doubled ball, with a cutoff equal to one
  on the unit ball, has a constant independent of the center, by
  homogeneity.
  Consequently
  \begin{align*}
    \norm{u_k}_{L^{2^{\ast}}(\mbBn)}^{2^{\ast}}
    &\leqslant
    C\left(\sup_\alpha\int_{B_{\mbH}(y_\alpha,2)}
                 |u_k|^{2^\ast}\diopv\right)^{\frac{2}{n}}
    \norm{u_k}_{H^1(\mbBn)}^2. 
  \end{align*}
  To obtain this inequality, write $\int|u_k|^{2^\ast}$ on each unit ball
  as its $\frac{2}{n}$ power times its $\frac{2}{2^\ast}$ power, bound
  the first factor by the displayed supremum, apply the local Sobolev
  inequality to the second factor, and sum using bounded overlap.

  Since $\norm{u_k}_{L^{2^{\ast}}(\mbBn)}^{2^{\ast}}\to nc>0$, there are centers $y_k$ and
  $\eta_{\mathrm{mass}}>0$ such that
  $\int_{B_{\mbH}(y_k,2)}|u_k|^{2^\ast}\diopv\geqslant\eta_{\mathrm{mass}}$.
  Take $\varphi_k(o)=y_k$ and, after a subsequence, let
  $u_k\circ\varphi_k\rightharpoonup Z$.
  The centers $y_k$ may approach the ideal boundary.
  The constants above do not depend on $y_k$, and the pullback by $\varphi_k$ preserves the mass, the energy and the dual residual norm.
  All cutoff and conformal estimates below take place in a fixed neighborhood of the origin.

  \emph{Step 2: $Z\not\equiv0$.}
  Suppose $Z\equiv0$.
  Fix $\Omega\Subset\mbBn$ containing $\ol{B_{\mbH}(o,2)}$ and a
  smooth compactly supported cutoff $\chi$ in $\Omega$ with
  $\chi=1$ on that closed ball.
  This is possible because $B_{\mbH}(o,2)=B_{\tanh1}$, the Euclidean ball
  of radius $\tanh1$ centered at the origin.
  Set $v_k:=\rho^{\frac{n-2}{2}}(u_k\circ\varphi_k)$.
  The conformal change of variables gives boundedness in $H^1(\Omega)$ and
  $v_k\to0$ in $L^2(\Omega)$.
  Multiplication by $\rho^{-\frac{n-2}{2}}$ followed by zero extension maps
  $H_0^1(\Omega)$ boundedly into $H^1(\mbBn)$.
  Testing the conformal identity with this inverse transform therefore
  transfers the residual to $H^{-1}(\Omega)$ and gives
  \begin{equation*}
    -\Delta v_k-\Lambda\rho^2v_k-|v_k|^{q-1}v_k\longrightarrow0.
  \end{equation*}
  The test $\chi^2v_k$ is bounded in $H_0^1(\Omega)$.
  Testing, and adding the cutoff gradient term, gives
  \begin{equation*}
    \int_{\mbRn}|\nabla(\chi v_k)|^2\diopx
    =\int_\Omega\chi^2|v_k|^{2^\ast}\diopx+o(1).
  \end{equation*}
  Here $\Lambda\rho^2$ and the derivatives of $\chi$ are bounded on
  $\Omega$, so their terms tend to zero by the local $L^2$ convergence.
  Extend $\chi v_k$ by zero.
  Recall the sharp Euclidean Sobolev inequality
  \begin{equation}\label{equ:crit-Sn}
    S_n\norm{w}_{L^{2^\ast}(\mbRn)}^2
    \leqslant\int_{\mbRn}|\nabla w|^2\diopx,
    \qquad w\in\dot H^1(\mbRn).
  \end{equation}
  Since $\chi^2|v_k|^{2^\ast}=|v_k|^{q-1}|\chi v_k|^2$,
  this inequality and H\"older's inequality imply
  \begin{align*}
    S_n\norm{\chi v_k}_{L^{2^\ast}(\mbRn)}^2
    &\leqslant\int_{\mbRn}|\nabla(\chi v_k)|^2\diopx\\
    &\leqslant
    \left(\int_\Omega|v_k|^{2^\ast}\diopx\right)^{\frac{2}{n}}
    \norm{\chi v_k}_{L^{2^\ast}(\mbRn)}^2+o(1)\\
    &\leqslant \norm{u_k}_{L^{2^{\ast}}(\mbBn)}^{q-1}\norm{\chi v_k}_{L^{2^\ast}(\mbRn)}^2+o(1).
  \end{align*}
  But $\norm{u_k}_{L^{2^{\ast}}(\mbBn)}^{q-1}\to(nc)^{\frac{2}{n}}<S_n$ by
  \eqref{equ:psdec-subthreshold-level}.
  For all sufficiently large $k$, the difference
  $S_n-\norm{u_k}_{L^{2^{\ast}}(\mbBn)}^{q-1}$ is at least
  $\frac{1}{2}(S_n-(nc)^{\frac{2}{n}})>0$.
  Absorption gives $\norm{\chi v_k}_{L^{2^\ast}(\mbRn)}\to0$.
  Conformal invariance of the $L^{2^\ast}$ norm contradicts the mass
  $\eta_{\mathrm{mass}}$ in $B_{\mbH}(o,2)$.
  Thus $Z\not\equiv0$. This is where the entry $S_n$ of $\wtilde{S}$ is used.

  \emph{Step 3: $Z$ has one sign.} Lemma~\ref{lem:psdec-first}, applied in this moving frame, shows that
  $\mcI'(Z)=0$ and that $u_k\circ\varphi_k-Z$ is a
  Palais--Smale sequence at level $c-\mcI(Z)$.
  Applying Lemma~\ref{lem:psdec-first} to this remainder gives
  $c-\mcI(Z)\geqslant0$.
  If $Z$ had both a nonzero positive part and a nonzero negative part,
  set $Z^+:=\max\{Z,0\}$ and $Z^-:=\max\{-Z,0\}$.
  Testing its equation against $Z^+$ and $-Z^-$ gives
  $\norm{Z^\pm}_\lambda^2=\norm{Z^\pm}_{L^{2^{\ast}}(\mbBn)}^{2^{\ast}}$.
  The Poincar\'e--Sobolev inequality would then give
  $\norm{Z^\pm}_\lambda^2\geqslant S_\lambda^{\frac{n}{2}}$.
  Since the cross terms between the two parts equal zero,
  $\norm{Z}_\lambda^2\geqslant2S_\lambda^{\frac{n}{2}}$.
  This contradicts
  $\frac{1}{n}\norm{Z}_\lambda^2=\mcI(Z)\leqslant c
  <\frac{2}{n}S_\lambda^{\frac{n}{2}}$.
  Hence $Z$ has one sign. Here, as in Step~5, the entry $2^{\frac{2}{n}}S_\lambda$ of $\wtilde{S}$ enters through $c<\frac{2}{n}S_\lambda^{\frac{n}{2}}$.

  \emph{Step 4: identification of $Z$.} Replacing $Z$ by $-Z$ if necessary, we may assume that $Z\geqslant0$. Then $Z$ is not identically zero and satisfies $(-\Delta_{\mbBn}-\lambda)Z=Z^q$ weakly. We show that $Z$ is a positive classical solution, so that the classification applies.
  Use $\rho$ and $\Lambda$ from
  Section~\ref{subsec:conformal-dictionary}.
  The conformal change of variables, polarized against compactly supported tests,
  gives
  \begin{equation}\label{equ:psdec-conformal}
    v:=\rho^{\frac{n-2}{2}}Z,\qquad
    -\Delta v=\Lambda\rho^2v+v^q.
  \end{equation}
  Thus $-\Delta v=(\Lambda\rho^2+v^{q-1})v$, and the coefficient
  lies in $L^{\frac{n}{2}}_{\mrloc}$
  because $(q-1)\frac{n}{2}=2^\ast$.
  The local Brezis--Kato estimate
  \cite{BrezisKato1979}, in the form
  \cite[Lemma~B.3]{StruweVM}, applies to this weak $H^1_{\mrloc}$
  solution and gives $v\in L^s_{\mrloc}$ for every finite $s$.
  Indeed the coefficient has precisely the required local
  $L^{\frac{n}{2}}$ integrability; no boundary condition is needed for
  this local conclusion.
  The right-hand side of \eqref{equ:psdec-conformal} therefore belongs
  to every finite local $L^s$ space.
  Interior elliptic regularity
  \cite[Theorem~8.8 and Lemma~9.16]{GilbargTrudinger2001}
  gives $v\in W^{2,s}_{\mrloc}$.
  Taking $s>n$ and using the Sobolev embedding gives
  $v\in C^{1,\alpha}_{\mrloc}$ for some $\alpha>0$.
  Schauder regularity
  \cite[Theorem~9.19]{GilbargTrudinger2001}
  then gives $v\in C^{2,\alpha}_{\mrloc}$, after reducing $\alpha$
  if needed.
  Since $\Lambda>0$ and $v\geqslant0$,
  $-\Delta v\geqslant0$.
  The strong maximum principle
  \cite[Theorem~3.5]{GilbargTrudinger2001} gives $v>0$.
  On each compact set $v$ is bounded away from zero.
  Further interior regularity gives smoothness of $v$ and $Z$.
  The classification \cite[Theorem~1.3]{ManciniSandeep2008}
  now gives $Z=\mcU_\xi$ for some $\xi\in\mbBn$.
  Its energy and mass follow from
  \eqref{equ:ground-state-energy-normalization}.
  Thus $Z=\pm\mcU_\xi$ and
  $\mcI(Z)=\frac{1}{n}S_\lambda^{\frac{n}{2}}$.

  \emph{Step 5: the remainder vanishes.} The remainder level satisfies
  \begin{equation*}
    0\leqslant c-\mcI(Z)<\frac{1}{n}S_\lambda^{\frac{n}{2}},
  \end{equation*}
  since $c<\frac{2}{n}S_\lambda^{\frac{n}{2}}$ and $\mcI(Z)=\frac{1}{n}S_\lambda^{\frac{n}{2}}$.
  Lemma~\ref{lem:psdec-first} forces that level to be zero and the
  remainder to converge strongly to zero.
  It also gives $c=\mcI(Z)$, as required. This finishes the proof of Lemma~\ref{lem:psdec-subthreshold}.
\end{proof}

\begin{proposition}\label{prop:low-energy-orbit}
  Let $f_0\in\mcM$ satisfy $\norm{f_0}_\lambda^2<\wtilde{S}$.
  The trajectory of Proposition~\ref{prop:flow-properties} satisfies
  \begin{equation}\label{equ:prelim-low-energy-orbit}
    \norm{f_\tau}_\lambda^2\longrightarrow S_\lambda,\qquad
    \dist(f_\tau,\mcZ)\longrightarrow0,\qquad
    \sigma(f_\tau)\longrightarrow0
    \quad(\tau\to\infty).
  \end{equation}
\end{proposition}
\begin{proof}
  We first identify the limit energy as $S_\lambda$ by rescaling the trajectory along any sequence of times into a Palais--Smale sequence below the threshold, and then apply the same argument along a contradicting sequence of times to show that the distance tends to zero.

  \emph{Step 1: the energy limit.} The mass constraint, \eqref{equ:PS} and energy monotonicity give a limit
  $E_\infty\in[S_\lambda,\norm{f_0}_\lambda^2]$.
  Choose any sequence $\tau_k\to\infty$ and set
  $u_k:=E_\infty^{\frac{1}{q-1}}f_{\tau_k}$.
  By \eqref{equ:flow-full-residual},
  \begin{align*}
    \mcI'(u_k)
    &=E_\infty^{\frac{1}{q-1}}
    \bigl((-\Delta_{\mbBn}-\lambda)f_{\tau_k}
      -E_\infty|f_{\tau_k}|^{q-1}f_{\tau_k}\bigr)
    \longrightarrow0 \quad\text{in }H^{-1}(\mbBn).
  \end{align*}
  Indeed replacing $\norm{f_{\tau_k}}_\lambda^2$ by $E_\infty$ contributes at most
  $S_\lambda^{-\frac{1}{2}}|\norm{f_{\tau_k}}_\lambda^2-E_\infty|$ before the fixed scalar
  factor.
  Since $\frac{2}{q-1}+1=\frac{n}{2}$ and
  $\frac{2^\ast}{q-1}=\frac{n}{2}$,
  \begin{equation*}
    \norm{u_k}_\lambda^2\longrightarrow E_\infty^{\frac{n}{2}},\qquad
    \norm{u_k}_{L^{2^{\ast}}(\mbBn)}^{2^{\ast}}=E_\infty^{\frac{n}{2}},\qquad
    \mcI(u_k)\longrightarrow\frac{1}{n}E_\infty^{\frac{n}{2}}.
  \end{equation*}
  This level is positive and strictly less than
  $\frac{1}{n}\wtilde{S}^{\frac{n}{2}}$.
  Lemma~\ref{lem:psdec-subthreshold} therefore gives
  $E_\infty^{\frac{n}{2}}=S_\lambda^{\frac{n}{2}}$ and, after a subsequence,
  isometries $\varphi_k$ and $Z\in\mcZ$ with
  $u_k\circ\varphi_k\to Z$ in $H^1(\mbBn)$.
  Thus $E_\infty=S_\lambda$.

  \emph{Step 2: the distance.} If the distance in \eqref{equ:prelim-low-energy-orbit} did not tend to
  zero, there would be times $s_k\to\infty$ and $\ve>0$ such that
  $\dist(f_{s_k},\mcZ)\geqslant\ve$ for every $k$.
  Apply Step~1 to this time sequence.
  It gives a subsequence for which
  $\dist(S_\lambda^{\frac{1}{q-1}}f_{s_k},\mcZ)\to0$.
  Here we used
  $\dist(h\circ\varphi,\mcZ)=\dist(h,\mcZ)$ for every isometry $\varphi$.
  Since $\mcZ$ is a cone containing both signs,
  $\dist(th,\mcZ)=|t|\dist(h,\mcZ)$ for every nonzero real $t$.
  This contradicts the chosen lower bound.
  Hence $\dist(f_\tau,\mcZ)\to0$ along the whole trajectory.
  Finally $\norm{f_\tau}_\lambda^2\geqslant S_\lambda>0$ gives the limit of $\sigma$. This finishes the proof of Proposition~\ref{prop:low-energy-orbit}.
\end{proof}

\subsubsection{Transfer to the initial datum}\label{subsec:crit-const-global}
\begingroup
The factor $s$ in the lower bound comes from the identity in \eqref{equ:const-trivial-bound}: on $\mcM$, $\sigma$ times the stability quotient equals $1-\frac{S_\lambda}{\norm{f}_\lambda^2}$, so it is nonincreasing along the flow. At the first time when the trajectory reaches $\sigma=s$, the local estimate gives at least $s\nu(s)$ for this product, and the bound passes back to the initial datum, where $\sigma\leqslant1$. The same transfer, at a first entrance time and through a quantity monotone along a flow, is used in \cite[Theorem~3.1]{DEFFL2025} with a continuous rearrangement flow and in \cite[Section~2.1]{CLTW2026} with the Yamabe flow, and there too the local estimate is multiplied by the relative distance at the first entrance. Those flows are applied to nonnegative functions, and the Yamabe flow converges for data of any energy. No such convergence result is available for our flow on $\mbBn$, which is why we need the threshold $\wtilde{S}$ and the high-energy region. In exchange, our flow acts directly on sign-changing functions, so no reduction from sign-changing to nonnegative functions is needed; in \cite[Proposition~3.10]{DEFFL2025} that reduction costs a factor $\frac{1}{2}$.\endgroup\par
\begingroup
\begin{proposition}\label{prop:const-transfer}
  Every $f\in\mcM\setminus\mcZ$ with $\norm{f}_\lambda^2<\wtilde{S}$ satisfies
  \begin{equation}\label{equ:const-interface}
    \frac{\delta(f)^2}{\dist(f,\mcZ)^2}
    \geqslant\sup_{0<s<1}s\nu(s).
  \end{equation}
\end{proposition}
\endgroup
\begin{proof}
  Fix $f_0\in\mcM\setminus\mcZ$ with $\norm{f_0}_\lambda^2<\wtilde{S}$ and $s\in(0,1)$.
  If $\nu(s)\leqslant0$, there is nothing to prove, so let $\nu(s)>0$.
  We first record an elementary lower bound for the quotient of $f_0$, and then distinguish whether $f_0$ already lies in $\mcT_s$ or reaches it only along the flow.

  By Lemma~\ref{lem:const-trivial-bound}, the distance is positive and
  \begin{equation*}
    \frac{\delta(f_0)^2}{\dist(f_0,\mcZ)^2}
    \geqslant1-\frac{S_\lambda}{\norm{f_0}_\lambda^2}.
  \end{equation*}

  \emph{Case $0<\sigma(f_0)\leqslant s$.} Proposition~\ref{prop:local-lower-function} gives
  the desired bound, because $\nu(s)>0$ and $s<1$.

  \emph{Case $\sigma(f_0)>s$.}
  The trajectory is continuous in $H^1$, and $\norm{f_\tau}_\lambda^2\geqslant S_\lambda$,
  so $\sigma(f_\tau)$ is continuous.
  Its limit is zero by Proposition~\ref{prop:low-energy-orbit}.
  The first time $\tau_s$ with $\sigma(f_{\tau_s})\leqslant s$ is
  therefore finite, is positive, and satisfies $\sigma(f_{\tau_s})=s$.
  At this time, Proposition~\ref{prop:local-lower-function} gives
  \begin{equation*}
    s\nu(s)
    \leqslant s\frac{\delta(f_{\tau_s})^2}{\dist(f_{\tau_s},\mcZ)^2}
    =1-\frac{S_\lambda}{\norm{f_{\tau_s}}_\lambda^2}
    \leqslant1-\frac{S_\lambda}{\norm{f_0}_\lambda^2}
    \leqslant\frac{\delta(f_0)^2}{\dist(f_0,\mcZ)^2},
  \end{equation*}
  where the equality is the identity in \eqref{equ:const-trivial-bound} at $f_{\tau_s}$, and in the second inequality we have used that the energy is nonincreasing along the flow.
  Taking the supremum over $s$ proves
  \eqref{equ:const-interface}. This finishes the proof of Proposition~\ref{prop:const-transfer}.
\end{proof}

\begin{proof}[Proof of Proposition~\ref{prop:weak-global-lower-bound}]
  We first show that the minimum in \eqref{equ:const-global-quotient} is positive, and then prove \eqref{equ:const-global-quotient} by splitting according to the energy of the normalized function.

  \emph{Step 1: positivity.} By \eqref{equ:const-nu-explicit}, $\nu(s)\to1-\frac{q}{\mu_3}>0$ as $s\to 0$, so $s\nu(s)>0$ for small $s>0$ and
  $\sup_{0<s<1}s\nu(s)>0$.
  Also $\min\{c_{\mathrm{conc}},\frac{1}{2}c_{\mathrm{dich}}\}>0$ by
  Proposition~\ref{prop:const-high-energy}. Thus the minimum is positive.

  \emph{Step 2: the inequality.} For a nonzero $f\notin\mcZ$, let
  $\what f:=\frac{f}{\norm{f}_{L^{2^{\ast}}(\mbBn)}}$.
  The stability quotient is unchanged by this normalization, and
  $\what f\in\mcM\setminus\mcZ$.
  If $\norm{\what f}_\lambda^2\geqslant \wtilde{S}$,
  Proposition~\ref{prop:const-high-energy} gives the bound
  $\min\{c_{\mathrm{conc}},\frac{1}{2}c_{\mathrm{dich}}\}$.
  If $\norm{\what f}_\lambda^2<\wtilde{S}$, Proposition~\ref{prop:const-transfer}
  gives the bound $\sup_{0<s<1}s\nu(s)$.
  Taking the smaller of the two bounds proves
  \eqref{equ:const-global-quotient}.
  Taking the infimum proves the assertion about $c_{\mathrm{BE}}$. This finishes the proof of Proposition~\ref{prop:weak-global-lower-bound}.
\end{proof}

\section{Strict upper bounds}\label{sec:strict-upper}
This section proves the strict upper bounds in Theorem~\ref{thm:crit-global-bounds}. Each bound comes from a test family whose stability quotient tends to the bound, and an expansion of the same family to the next order shows that the quotient goes strictly below it. The three test families are those of Section~\ref{sec:intro}: the functions $\mcU+\ve\Psi$ with $\Psi$ an eigenfunction associated with $\mu_3$ and $\ve\to0$ (the local channel), Euclidean bubbles concentrating at a point (the concentration channel), and sums of two ground states whose centers move apart (the dichotomy channel). Under the normalization $\norm{f}_{L^{2^{\ast}}(\mbBn)}=1$, along these families the energy $\norm{f}_\lambda^2$ tends to $S_\lambda$, $S_n$ and $2^{\frac{2}{n}}S_\lambda$, the relative distance $\sigma$ tends to $0$, $1$ and $\frac{1}{2}$, and the stability quotient tends to $c_{\mrloc}$, $c_{\mathrm{conc}}$ and $c_{\mathrm{dich}}$, respectively. They are also the only ways in which a minimizing sequence for \eqref{equ:cbe-def} can lose compactness (Section~\ref{sec:attainment}). Sections~\ref{subsec:strict-conc}, \ref{subsec:strict-dich} and~\ref{subsec:strict-local} treat the concentration, dichotomy and local channels, respectively. Section~\ref{subsec:proof-global} completes the proof of Theorem~\ref{thm:crit-global-bounds}, combining these bounds with Proposition~\ref{prop:weak-global-lower-bound} and removing $c_{\mathrm{conc}}$ from the lower bound.

\subsection{The concentration channel}\label{subsec:strict-conc}
The concentration channel is tested by Euclidean bubbles concentrating at the origin, transplanted to $\mbBn$ by the conformal change of variables of Section~\ref{subsec:conformal-dictionary}. These are the localized Aubin--Talenti bubbles of \cite[(3.16)]{BhaktaSandeep2012}. We fix once and for all a cut-off
\begin{equation*}
  \eta\in \Ccit[0,\tfrac{1}{2}),\qquad0\leqslant\eta\leqslant1,\qquad\eta\equiv1\text{ on }[0,\tfrac{1}{4}],
\end{equation*}
and, for $\mu\geqslant1$, we set
\begin{equation*}
  v_\mu(x):=\eta(|x|)\,\mu^{\frac{n-2}{2}}V(\mu x),\qquad x\in B_1,
\end{equation*}
where $V$ is the Aubin--Talenti bubble \eqref{equ:AT-bubble}. Thus $v_\mu\in \Ccit(B_{1/2})$. Items~\ref{it:conc-grad}--\ref{it:conc-L2} of Lemma~\ref{lem:conc-profile} below are the estimates of Brezis and Nirenberg \cite[proof of Lemma~1.1, (1.11)--(1.13)]{BrezisNirenberg1983}, written in our normalization; we include the short proof since we need the explicit constant for $n=4$ and the additional estimates~\ref{it:conc-moment} and~\ref{it:conc-L1}.

\begin{lemma}\label{lem:conc-profile}
  Let $K:=\int_{\mbRn}V^2\diopy$ if $n\geqslant5$. As $\mu\to\infty$, the following hold, with constants depending only on $n$ and $\eta$.
  \begin{enumerate}[label=\textnormal{(\roman*)},leftmargin=3.2em]
    \item\label{it:conc-grad}
    $\int_{B_1}|\nabla v_\mu|^2\diopx=S_n^{\frac{n}{2}}+O(\mu^{-(n-2)})$.
    \item\label{it:conc-crit}
    $\int_{B_1}v_\mu^{2^{\ast}}\diopx=S_n^{\frac{n}{2}}+O(\mu^{-n})$.
    \item\label{it:conc-L2}
    $\int_{B_1}v_\mu^2\diopx=K\mu^{-2}+O(\mu^{-(n-2)})$ if $n\geqslant5$, and $\int_{B_1}v_\mu^2\diopx=16\pi^2\mu^{-2}\log\mu+O(\mu^{-2})$ if $n=4$.
    \item\label{it:conc-moment}
    $\int_{B_1}|x|^2v_\mu^2\diopx$ is $O(\mu^{-2})$ if $n=4$, $O(\mu^{-3})$ if $n=5$, $O(\mu^{-4}\log\mu)$ if $n=6$, and $O(\mu^{-4})$ if $n\geqslant7$.
    \item\label{it:conc-L1}
    $\int_{B_1}v_\mu\diopx\leqslant C\mu^{-\frac{n-2}{2}}$.
  \end{enumerate}
\end{lemma}

\begin{proof}
  All integrals are computed by the substitution $y=\mu x$. Under it,
  \begin{gather*}
    \int_{\mbRn}\bigl|\nabla_x\bigl(\mu^{\frac{n-2}{2}}V(\mu x)\bigr)\bigr|^2\diopx=\int_{\mbRn}|\nabla V|^2\diopy=S_n^{\frac{n}{2}},\\
    \int_{\mbRn}\bigl(\mu^{\frac{n-2}{2}}V(\mu x)\bigr)^{2^{\ast}}\diopx=\int_{\mbRn}V^{2^{\ast}}\diopy=S_n^{\frac{n}{2}}.
  \end{gather*}
  By \eqref{equ:AT-bubble}, $V(x)\leqslant\alpha_n\min\{1,|x|^{2-n}\}$ and $|\nabla V(x)|\leqslant(n-2)\alpha_n|x|^{1-n}$ for $x\in\mbRn$. Since $\eta(|x|)=1$ for $|x|\leqslant\frac{1}{4}$, the function $v_\mu$ differs from $\mu^{\frac{n-2}{2}}V(\mu x)$ only on $\{|x|>\frac{1}{4}\}$, and there $\mu^{\frac{n-2}{2}}V(\mu x)\leqslant4^{n-2}\alpha_n\mu^{-\frac{n-2}{2}}$.

  \emph{Proof of \ref{it:conc-grad}.} On $\{|x|>\frac{1}{4}\}$ we have $|\nabla v_\mu|^2\leqslant2\mu^{n}|\nabla V(\mu x)|^2+2\norm{\eta'}_{L^\infty(0,\frac{1}{2})}^2\mu^{n-2}V(\mu x)^2$, and $\eta'$ vanishes outside $[\frac{1}{4},\frac{1}{2}]$. Hence
  \begin{align*}
    \Bigl|\int_{B_1}|\nabla v_\mu|^2\diopx-S_n^{\frac{n}{2}}\Bigr|
    &\leqslant3\int_{\{|y|>\frac{\mu}{4}\}}|\nabla V|^2\diopy
    +2\norm{\eta'}_{L^\infty(0,\frac{1}{2})}^2\int_{\{\frac{1}{4}<|x|<\frac{1}{2}\}}\mu^{n-2}V(\mu x)^2\diopx\\
    &\leqslant C\mu^{-(n-2)},
  \end{align*}
  since $\int_{\{|y|>\frac{\mu}{4}\}}|\nabla V|^2\diopy\leqslant C\int_{\frac{\mu}{4}}^\infty r^{1-n}\diopr=C\mu^{2-n}$ and the last integrand is at most $C\mu^{-(n-2)}$.

  \emph{Proof of \ref{it:conc-crit}.} Since $0\leqslant v_\mu\leqslant\mu^{\frac{n-2}{2}}V(\mu x)$, with equality for $|x|\leqslant\frac{1}{4}$,
  \begin{equation*}
    0\leqslant S_n^{\frac{n}{2}}-\int_{B_1}v_\mu^{2^{\ast}}\diopx
    \leqslant\int_{\{|y|>\frac{\mu}{4}\}}V^{2^{\ast}}\diopy
    \leqslant C\int_{\frac{\mu}{4}}^\infty r^{-n-1}\diopr
    =C\mu^{-n}.
  \end{equation*}

  \emph{Proof of \ref{it:conc-L2}.} We have
  \begin{equation*}
    \int_{B_1}v_\mu^2\diopx=\mu^{-2}\int_{\mbRn}\eta\bigl(\tfrac{|y|}{\mu}\bigr)^2V(y)^2\diopy,
  \end{equation*}
  and the integrand lies between $V^2\mathbf 1_{\{|y|<\frac{\mu}{4}\}}$ and $V^2\mathbf 1_{\{|y|<\frac{\mu}{2}\}}$. If $n\geqslant5$, then $V\in L^2(\mbRn)$ and $\int_{\{|y|>\frac{\mu}{4}\}}V^2\diopy\leqslant C\int_{\frac{\mu}{4}}^\infty r^{3-n}\diopr=C\mu^{4-n}$, which gives the claim. If $n=4$, then $\alpha_4^2=8$, and as $R\to\infty$,
  \begin{equation*}
    \int_{\{|y|<R\}}V^2\diopy=16\pi^2\int_0^R\frac{r^3}{(1+r^2)^2}\diopr=16\pi^2\log R+O(1).
  \end{equation*}
  Taking $R=\frac{\mu}{4}$ and $R=\frac{\mu}{2}$ gives the claim.

  \emph{Proof of \ref{it:conc-moment}.} In the same way,
  \begin{equation*}
    \int_{B_1}|x|^2v_\mu^2\diopx\leqslant\mu^{-4}\int_{\{|y|<\frac{\mu}{2}\}}|y|^2V^2\diopy,
  \end{equation*}
  and $|y|^2V(y)^2\leqslant C\min\{|y|^2,|y|^{6-2n}\}$. Integration in polar coordinates gives the stated orders.

  \emph{Proof of \ref{it:conc-L1}.} We have
  \begin{align*}
    \int_{B_1}v_\mu\diopx
    &\leqslant\mu^{\frac{n-2}{2}}\int_{\{|x|<\frac{1}{2}\}}V(\mu x)\diopx
    =\mu^{-\frac{n+2}{2}}\int_{\{|y|<\frac{\mu}{2}\}}V\diopy\\
    &\leqslant C\mu^{-\frac{n+2}{2}}\int_{\{|y|<\frac{\mu}{2}\}}\min\{1,|y|^{2-n}\}\diopy
    \leqslant C\mu^{-\frac{n-2}{2}}.
  \end{align*}
  This finishes the proof of Lemma~\ref{lem:conc-profile}.
\end{proof}

\begin{proposition}\label{prop:strict-conc}
  We have $c_{\mathrm{BE}}<c_{\mathrm{conc}}$.
\end{proposition}

\begin{proof}
  We use the functions $W_\mu$ defined below as test functions and show that their stability quotient lies strictly below $c_{\mathrm{conc}}$ for all large $\mu$. We first show that $1-\sigma(W_\mu)$ decays like $\mu^{-(n-2)}$, then compute the first correction to the energy ratio $\mcE_\mu$ below, which has a definite sign and is of larger order, and finally combine the two through an exact identity for the quotient.

  \emph{Setup.} Recall that $\Lambda=\lambda-\lambda_0>0$. All constants in this proof may depend on $n$, $\lambda$ and $\eta$. Since $\supp v_\mu\subset B_{1/2}$ and $\rho(x)=\frac{2}{1-|x|^2}$, we have
  \begin{equation}\label{equ:conc-rho}
    2\leqslant\rho\leqslant\frac{8}{3}
    \quad\text{and}\quad
    0\leqslant\rho^2-4\leqslant\frac{128}{9}|x|^2
    \qquad\text{on }\supp v_\mu.
  \end{equation}
  By the conformal change of variables \eqref{equ:dictionary}, the function $\rho^{-\frac{n-2}{2}}v_\mu$ belongs to $H^1(\mbBn)$, and
  \begin{equation*}
    \bignorm{\rho^{-\frac{n-2}{2}}v_\mu}_\lambda^2=\int_{B_1}|\nabla v_\mu|^2\diopx-\Lambda\int_{B_1}\rho^2v_\mu^2\diopx,
    \qquad
    \bignorm{\rho^{-\frac{n-2}{2}}v_\mu}_{L^{2^{\ast}}(\mbBn)}=\norm{v_\mu}_{L^{2^{\ast}}(B_1)}.
  \end{equation*}
  By Lemma~\ref{lem:conc-profile}, items~\ref{it:conc-grad} and~\ref{it:conc-L2}, and by \eqref{equ:conc-rho}, the first quantity tends to $S_n^{\frac{n}{2}}$, so it is positive for all large $\mu$. For such $\mu$ let
  \begin{equation*}
    W_\mu:=\frac{S_\lambda^{\frac{n}{4}}}{\bignorm{\rho^{-\frac{n-2}{2}}v_\mu}_\lambda}\,\rho^{-\frac{n-2}{2}}v_\mu,
    \qquad
    \mcE_\mu:=\frac{\norm{W_\mu}_\lambda^2}{\norm{W_\mu}_{L^{2^{\ast}}(\mbBn)}^2}
    =\frac{\int_{B_1}|\nabla v_\mu|^2\diopx-\Lambda\int_{B_1}\rho^2v_\mu^2\diopx}{\norm{v_\mu}_{L^{2^{\ast}}(B_1)}^2}.
  \end{equation*}
  Then $\norm{W_\mu}_\lambda^2=S_\lambda^{\frac{n}{2}}$, the normalization factor in the definition of $W_\mu$ stays bounded, and $\mcE_\mu\to S_n$ by Lemma~\ref{lem:conc-profile}\ref{it:conc-crit}.

  \emph{Step 1: the pairing estimate.} We show that $W_\mu$ is almost orthogonal to every $\mcU_\xi$. Since $\diopv=\rho^n\diopx$, the bound $\rho\leqslant\frac{8}{3}$ in \eqref{equ:conc-rho}, the boundedness of the normalization factor and Lemma~\ref{lem:conc-profile}\ref{it:conc-L1} give
  \begin{equation*}
    \norm{W_\mu}_{L^1(\mbBn)}
    \leqslant C\int_{B_1}\rho^{\frac{n+2}{2}}v_\mu\diopx
    \leqslant C\mu^{-\frac{n-2}{2}}.
  \end{equation*}
  Hence Lemmas~\ref{lem:ground-state} and~\ref{lem:duality} give
  \begin{equation*}
    \sup_{\xi\in\mbBn}
    |\lag W_\mu,\mcU_\xi\rag_\lambda|
    \leqslant
    \norm{\mcU}_{L^\infty(\mbBn)}^q\norm{W_\mu}_{L^1(\mbBn)}
    \leqslant C\mu^{-\frac{n-2}{2}}.
  \end{equation*}
  The distance formula therefore yields
  \begin{equation*}
    1-\sigma(W_\mu)
    =
    S_\lambda^{-n}
    \sup_{\xi\in\mbBn}
    \lag W_\mu,\mcU_\xi\rag_\lambda^2
    =O(\mu^{-(n-2)}).
  \end{equation*}

  \emph{Step 2: the first correction to $\mcE_\mu$.}
  Suppose first that $n\geqslant5$, and let $K=\int_{\mbRn}V^2\diopy$. Writing $\rho^2=4+(\rho^2-4)$ and using \eqref{equ:conc-rho} together with Lemma~\ref{lem:conc-profile}, items~\ref{it:conc-L2} and~\ref{it:conc-moment}, we get
  \begin{equation*}
    \int_{B_1}\rho^2v_\mu^2\diopx=4K\mu^{-2}+o(\mu^{-2}).
  \end{equation*}
  The errors in items~\ref{it:conc-grad} and~\ref{it:conc-crit} of Lemma~\ref{lem:conc-profile} are $O(\mu^{-(n-2)})=o(\mu^{-2})$. Expanding the outer power in the denominator now yields
  \begin{equation*}
    \mcE_\mu
    =
    S_n-4S_n^{1-\frac{n}{2}}\Lambda K\mu^{-2}
    +o(\mu^{-2}).
  \end{equation*}

  For $n=4$, the same argument gives
  $\int_{B_1}\rho^2v_\mu^2\diopx=64\pi^2\mu^{-2}\log\mu+O(\mu^{-2})$,
  while the errors in items~\ref{it:conc-grad} and~\ref{it:conc-crit} of Lemma~\ref{lem:conc-profile} are $O(\mu^{-2})$. Hence
  \begin{equation*}
    \mcE_\mu
    =
    S_4-\frac{64\pi^2\Lambda}{S_4}\mu^{-2}\log\mu
    +O(\mu^{-2}).
  \end{equation*}

\emph{Step 3: the sign of the quotient.}
  Since $1-\sigma(W_\mu)\to0$, we have
  $\dist(W_\mu,\mcZ)>0$ for all sufficiently large $\mu$.
  Homogeneity and the distance formula give the exact identity
  \begin{align*}
    &\sigma(W_\mu)
    \left(\frac{\delta(W_\mu)^2}{\dist(W_\mu,\mcZ)^2}
    -c_{\mathrm{conc}}\right)
    =-\frac{S_\lambda(S_n-\mcE_\mu)}{S_n\mcE_\mu}
    +c_{\mathrm{conc}}\bigl(1-\sigma(W_\mu)\bigr).
  \end{align*}
  Since $\mcE_\mu\to S_n$ and $\sigma(W_\mu)\to1$, the right-hand side tends to $0$, so the stability quotient of $W_\mu$ tends to $c_{\mathrm{conc}}$. We now determine the sign of the right-hand side by comparing the orders given by Steps~1 and~2.
  For $n\geqslant5$, the negative term has a positive coefficient
  times $\mu^{-2}$ in absolute value, while
  the other term is $O(\mu^{-(n-2)})=o(\mu^{-2})$.
  For $n=4$, the negative term has a positive coefficient times
  $\mu^{-2}\log\mu$ in absolute value, while
  the other term is $O(\mu^{-2})$.
  Thus the right-hand side is negative for all sufficiently large
  $\mu$.
  Each such $W_\mu$ is an admissible test function, so
  \begin{equation*}
    c_{\mathrm{BE}}
    \leqslant
    \frac{\delta(W_\mu)^2}{\dist(W_\mu,\mcZ)^2}
    <c_{\mathrm{conc}}.
  \end{equation*}
  This finishes the proof of Proposition~\ref{prop:strict-conc}.
\end{proof}

\subsection{The dichotomy channel}\label{subsec:strict-dich}

\begin{lemma}\label{lem:dich-elementary}
  Let $2<p\leqslant4$. There is a constant $C_p$ such that, for all $u,v\geqslant0$,
  \begin{equation*}
    \left|
    (u+v)^p-u^p-v^p-p(u^{p-1}v+uv^{p-1})
    \right|
    \leqslant C_p(uv)^{\frac{p}{2}}.
  \end{equation*}
\end{lemma}

\begin{proof}
  Both sides are symmetric in $(u,v)$, positively homogeneous of degree $p$, and vanish if $uv=0$. It therefore suffices to consider $u=1$ and $0<v\leqslant1$, where we must bound
  \begin{equation*}
    F(v):=(1+v)^p-1-v^p-pv-pv^{p-1}
  \end{equation*}
  by $C_pv^{\frac{p}{2}}$. Taylor's formula gives $|(1+v)^p-1-pv|\leqslant\frac{p(p-1)}{2}2^{p-2}v^2$ for $0<v\leqslant1$. Since also $v^p\leqslant v^{p-1}$, we get $|F(v)|\leqslant C_p'v^{\min\{2,p-1\}}$. Since $2<p\leqslant4$, we have $\min\{2,p-1\}\geqslant\frac{p}{2}$, and hence $v^{\min\{2,p-1\}}\leqslant v^{\frac{p}{2}}$ for $0<v\leqslant1$. This finishes the proof of Lemma~\ref{lem:dich-elementary}.
\end{proof}

\begin{proposition}\label{prop:strict-dich}
  We have $c_{\mathrm{BE}}<c_{\mathrm{dich}}$.
\end{proposition}

\begin{proof}
  Inspired by K\"onig's proof of the corresponding Euclidean bound \cite{Koenig2025}, we test with the sum of two ground states whose centers move apart, and measure everything against their overlap $h(D)$. We first expand the deficit up to $o(h(D))$, then prove that $h$ is strictly decreasing with bounded logarithmic derivative, use this to show that the squared distance to $\mcZ$ equals $S_\lambda^{\frac{n}{2}}+o(h(D))$, and finally divide the two expansions.

  \emph{Setup.} Choose $\xi_D$ on a fixed unit speed geodesic starting at $o$,
  with $d(o,\xi_D)=D$.
  Use the family
  $\vp_D=\mcU+\mcU_{\xi_D}$, and let
  \begin{equation*}
    h(D)=\lag\mcU,\mcU_{\xi_D}\rag_\lambda.
  \end{equation*}
  Lemmas~\ref{lem:ground-state-normalization},
  \ref{lem:duality} and~\ref{lem:overlap} give
  \begin{equation*}
    h(D)=\int_{\mbBn}\mcU^q\mcU_{\xi_D}\diopv
    =\int_{\mbBn}\mcU\mcU_{\xi_D}^q\diopv>0,
    \qquad
    h(D)\to0,
    \qquad
    \norm{\vp_D}_\lambda^2=2S_\lambda^{\frac{n}{2}}+2h(D).
  \end{equation*}

  \emph{Step 1: expansion of the deficit.} We first control the nonlinear remainder.
  The logarithmic decay from
  \cite[Lemma~3.4]{ManciniSandeep2008}, already used in
  Lemma~\ref{lem:ground-state}, gives, for every
  $\ve>0$,
  \begin{equation*}
    c_\ve\mre^{-(\beta_++\ve)r}
    \leqslant\mcU(r)
    \leqslant
    C_\ve\mre^{-(\beta_+-\ve)r}
    \qquad(r\geqslant0).
  \end{equation*}
  Integration over $B_{\mbH}(o,1)$ gives
  \begin{equation*}
    h(D)\geqslant
    c_\ve\mre^{-(\beta_++\ve)D}.
  \end{equation*}
  Choose $\ve>0$ so small that
  \begin{equation*}
    a=\frac{2^{\ast}}{2}(\beta_+-\ve)>\beta_++\ve.
  \end{equation*}

  Let $r=d(x,o)$ and $s=d(x,\xi_D)$.
  In cylindrical coordinates about the complete geodesic through the
  two centers, let $t\in\mbR$ be the axial coordinate and $z\geqslant0$
  the distance to the axis.
  Then
  \begin{equation*}
    \cosh r=\cosh z\cosh t,
    \qquad
    \cosh s=\cosh z\cosh(t-D),
    \qquad
    \diopv=
    \cosh z\,\sinh^{n-2}z\,
    \diopz\diopt\diopomega.
  \end{equation*}
  Since $2a>n-1$, these formulas imply
  \begin{align*}
    \int_{\mbBn}(\mcU\mcU_{\xi_D})^{\frac{2^{\ast}}{2}}\diopv
    &\leqslant
    C\int_0^\infty
    \cosh^{1-2a}z\,\sinh^{n-2}z\diopz
    \int_{\mbR}
    \mre^{-a(|t|+|t-D|)}\diopt\\
    &\leqslant C(1+D)\mre^{-aD}
    =o(h(D)).
  \end{align*}

  By Lemma~\ref{lem:dich-elementary} with $p=2^{\ast}$, which applies since $2<2^{\ast}\leqslant4$ for $n\geqslant4$, the integral of the remainder is bounded by $C\int_{\mbBn}(\mcU\mcU_{\xi_D})^{\frac{2^{\ast}}{2}}\diopv=o(h(D))$.
  It follows that
  \begin{equation*}
    \int_{\mbBn}\vp_D^{2^{\ast}}\diopv
    =2S_\lambda^{\frac{n}{2}}+2\cdot2^{\ast}h(D)+o(h(D)).
  \end{equation*}
  Since $\frac{n}{2}(1-\frac{2}{2^{\ast}})=1$, expansion of the outer power gives
  \begin{equation*}
    \delta(\vp_D)^2
    =
    \left(2-2^{\frac{2}{2^{\ast}}}\right)S_\lambda^{\frac{n}{2}}
    -2\left(2^{\frac{2}{2^{\ast}}}-1\right)h(D)
    +o(h(D)).
  \end{equation*}

  \emph{Step 2: properties of $h$.} We next prove the two properties of $h$ needed to refine the distance estimate.
  The equation
  $-\Delta_{\mbBn}\mcU=(\lambda+\mcU^{q-1})\mcU$ has a bounded
  coefficient.
  The Harnack inequality and interior gradient estimates on unit
  geodesic balls give
  \begin{equation*}
    |\nabla_{\mbBn}\mcU|\leqslant C\mcU
  \end{equation*}
  with a constant independent of the center
  \cite{GilbargTrudinger2001}.
  Here homogeneity identifies every such ball with a fixed ball,
  and Lemma~\ref{lem:ground-state} bounds $\mcU$ globally.
  Differentiation along the geodesic therefore gives
  \begin{equation*}
    |h'(s)|\leqslant Ch(s).
  \end{equation*}

  The function $h$ depends only on the distance between the centers.
  It is $C^2$ by Lemma~\ref{lem:pairing-center-derivative}.
  The ground-state equation in the center variable gives, by \eqref{equ:polar-laplacian}, for $s>0$,
  \begin{equation*}
    -\frac{\mrd}{\mrd s}
    \bigl(\sinh^{n-1}s\,h'(s)\bigr)
    =
    \sinh^{n-1}s
    \left(
    \lambda h(s)+
    \int_{\mbBn}\mcU^q\mcU_{\xi_s}^q\diopv
    \right)>0.
  \end{equation*}
  Since $h'(0)=0$ and $\lambda>0$, the function $h$ is strictly
  decreasing on $(0,\infty)$.
  Also $h(0)=S_\lambda^{\frac{n}{2}}$.

  \emph{Step 3: expansion of the distance.} By Lemma~\ref{lem:dist-formula}, we need the largest pairing of $\vp_D$ with a translate of $\mcU$. Let
  \begin{equation*}
    P_D=
    \max_{\zeta\in\mbBn}
    \lag\vp_D,\mcU_\zeta\rag_\lambda.
  \end{equation*}
  The maximum exists by continuity and Lemma~\ref{lem:overlap}.
  Let $\zeta_D$ be a maximizing center.
  Let $a_D$ and $b_D$ be the smaller and larger of
  $d(o,\zeta_D)$ and $d(\xi_D,\zeta_D)$.
  Then
  \begin{equation*}
    P_D=h(a_D)+h(b_D),
    \qquad
    a_D+b_D\geqslant D,
    \qquad
    b_D\geqslant\frac{D}{2}.
  \end{equation*}
  The competitor $\zeta=o$ gives $P_D\geqslant S_\lambda^{\frac{n}{2}}+h(D)$.
  Since $h(b_D)\leqslant h(\frac{D}{2})\to0$, we obtain
  $h(a_D)\to S_\lambda^{\frac{n}{2}}$.
  Strict decrease of $h$ gives $a_D\to0$.

  For large $D$, monotonicity and the bound on the logarithmic
  derivative yield
  \begin{equation*}
    h(b_D)
    \leqslant h(D-a_D)
    \leqslant \mre^{Ca_D}h(D)
    =h(D)+o(h(D)).
  \end{equation*}
  Together with $h(a_D)\leqslant S_\lambda^{\frac{n}{2}}$ and the preceding lower bound,
  this proves
  \begin{equation*}
    P_D=S_\lambda^{\frac{n}{2}}+h(D)+o(h(D)).
  \end{equation*}
  All pairings are positive.
  Lemma~\ref{lem:dist-formula} therefore gives
  \begin{equation*}
    \dist(\vp_D,\mcZ)^2
    =
    2S_\lambda^{\frac{n}{2}}+2h(D)-\frac{P_D^2}{S_\lambda^{\frac{n}{2}}}
    =S_\lambda^{\frac{n}{2}}+o(h(D)).
  \end{equation*}

  \emph{Step 4: conclusion.} Dividing the deficit expansion of Step~1 by the distance expansion of Step~3 gives
  \begin{equation*}
    \frac{\delta(\vp_D)^2}{\dist(\vp_D,\mcZ)^2}
    =
    2-2^{\frac{n-2}{n}}
    -\frac{2\left(2^{\frac{n-2}{n}}-1\right)}
    {S_\lambda^{\frac{n}{2}}}\,h(D)
    +o(h(D)).
  \end{equation*}
  The coefficient of $h(D)$ is negative, and $h(D)>0$.
  The distance tends to $S_\lambda^{\frac{n}{2}}>0$.
  Thus every sufficiently large finite $D$ gives an admissible
  test function whose quotient is strictly smaller than
  $2-2^{\frac{n-2}{n}}$.
  Taking the infimum proves the claim. This finishes the proof of Proposition~\ref{prop:strict-dich}.
\end{proof}

\subsection{The local channel}
\label{subsec:strict-local}
\begingroup
In this subsection we prove the following proposition.
\begin{proposition}\label{prop:strict-local}
  We have $c_{\mathrm{BE}}<c_{\mrloc}$.
\end{proposition}
The proof is given at the end of the subsection. We test the stability quotient along $\mcU+\ve\Psi$, where $\Psi$ is an eigenfunction of \eqref{equ:const-linearized} associated with $\mu_3$. To second order in $\ve$ the quotient equals $c_{\mrloc}$ (Theorem~\ref{thm:crit-local}), so strictness is decided by the third-order term, which is a multiple of $\ve\int_{\mbBn}\mcU^{q-2}\Psi^3\diopv$. Such a third-order expansion has also been used by K\"onig \cite{Koenig2023a} in the corresponding Euclidean case. Since the sign of $\ve$ is free, it suffices to find such a $\Psi$ with nonzero cubic integral; we find one with positive cubic integral. Section~\ref{subsec:strict-local-level} shows that $\mu_3=\min\{\mu_{0,2},\mu_{2,1}\}$, so that $\Psi$ can be taken radial or of angular degree $2$. Section~\ref{subsec:strict-local-expansion} expands the quotient to third order. Section~\ref{subsec:strict-local-sign} proves that the cubic integral is positive in both cases and then proves Proposition~\ref{prop:strict-local}.\endgroup\par

\subsubsection{The spectral level \texorpdfstring{$\mu_3$}{mu3}}\label{subsec:strict-local-level}
For $\ell\geqslant0$, let $\mu_{\ell,j}$ be the $j$-th radial
eigenvalue of \eqref{equ:const-linearized} with angular degree $\ell$.
Here $j$ counts the simple radial eigenvalues for one fixed spherical
harmonic. It does not count the multiplicity of the spherical
harmonics.

\begin{lemma}\label{lem:strict-local-sectors}
  We have
  \begin{equation*}
    \mu_3=\min\{\mu_{0,2},\mu_{2,1}\}.
  \end{equation*}
\end{lemma}

\begin{proof}
  We first show that the radial eigenvalues are nondecreasing in the angular degree, then that they are simple in each angular sector, and finally we locate the eigenvalues $1$ and $q$ and read off $\mu_3$.

  \emph{Step 1: monotonicity in $\ell$.} Expanding in spherical harmonics decomposes $H^1(\mbBn)$ orthogonally into angular sectors respected by both quadratic forms, so the spectrum of \eqref{equ:const-linearized} is the union of the radial spectra $\{\mu_{\ell,j}\}$. Let $Y$ be a spherical harmonic of degree $\ell$. For a separated function $R(r)Y(\theta)$, by \eqref{equ:polar-volume} and \eqref{equ:polar-gradient}, the radial Rayleigh quotient is
  \begin{equation*}
    \frac{
    \ds\int_0^\infty
    \left(R'^2+\frac{\ell(\ell+n-2)}{\sinh^2r}R^2-\lambda R^2\right)
    (\sinh r)^{n-1}\diopr
    }{
    \ds\int_0^\infty
    \mcU^{q-1}R^2(\sinh r)^{n-1}\diopr
    }.
  \end{equation*}
  The angular factor cancels, the denominator is independent of $\ell$, the numerator increases with $\ell$, and the form domain for degree $\ell+1$ is contained in the form domain for degree $\ell$. The min--max principle therefore gives $\mu_{\ell,j}\leqslant\mu_{\ell+1,j}$.

  \emph{Step 2: the radial problems.} The weighted embedding is compact by \cite[Lemmas~3.1 and~3.2]{BGKM2025}, and therefore so is its restriction to each angular sector. The radial realization uses the regular solution at $r=0$ and the finite-energy solution at infinity. Conjugation by $(\sinh r)^{\frac{n-1}{2}}$ gives the inverse-square coefficient
  \begin{equation*}
    \ell(\ell+n-2)+\frac{(n-1)(n-3)}{4}
    \geqslant\frac{3}{4}
  \end{equation*}
  at the origin, hence the origin is a limit-point endpoint. At infinity, the finite-energy condition selects the fast solution and defines a separated self-adjoint boundary condition. The compact weighted embedding shows that the essential spectrum is empty, and the quadratic form is bounded below. By \cite[Theorem~10.12.1(8)(ii)]{Zettl2005}, the radial eigenvalues are simple, and the $j$-th radial eigenfunction has exactly $j-1$ zeros in $(0,\infty)$.

  \emph{Step 3: identification of $\mu_3$.} By Lemma~\ref{lem:linearized-eigenspaces}\ref{item:eigen-one}--\ref{item:eigen-q}, the eigenspace of $1$ is generated by the radial function $\mcU$, and the eigenspace of $q$ is generated by $\Phi_1,\ldots,\Phi_n$, which have angular degree $1$. So $1$ occurs only in degree $0$ and $q$ only in degree $1$; since $1=\mu_1$ and $q=\mu_2$, this gives $\mu_{0,1}=1$ and $\mu_{1,1}=q$. Every other eigenvalue is at least $\mu_3$, hence strictly above $q$, since $\mu_3>q$ (see Section~\ref{subsec:linearized-spectrum}).

  By Steps~1 and~2, every remaining radial eigenvalue in degree $0$ is at least $\mu_{0,2}$, every remaining radial eigenvalue in degree $1$ is at least the second radial eigenvalue $\mu_{0,2}$ in degree $0$, and every eigenvalue in degree $\ell\geqslant2$ is at least $\mu_{2,1}$. Both $\mu_{0,2}$ and $\mu_{2,1}$ occur in the full spectrum and differ from $1$ and $q$. This proves the asserted identity, which finishes the proof of Lemma~\ref{lem:strict-local-sectors}.
\end{proof}

\subsubsection{Expansion to third order}\label{subsec:strict-local-expansion}
\begin{lemma}\label{lem:strict-local-cubic}
  Let $\Psi=R(r)Y(\theta)$ be a separated eigenfunction of \eqref{equ:const-linearized}, where $Y$ is a spherical harmonic. Then $\frac{\Psi}{\mcU}$ is bounded, so that the cubic integral $\int_{\mbBn}\mcU^{q-2}\Psi^3\diopv$ converges absolutely. If moreover $\Psi$ has eigenvalue $\mu_3$ and $\norm{\Psi}_\lambda=1$, then $f_\ve=\mcU+\ve\Psi$ satisfies $\dist(f_\ve,\mcZ)^2=\ve^2$ for all sufficiently small $|\ve|$, and
  \begin{equation}\label{equ:strict-local-cubic}
    \frac{\delta(f_\ve)^2}{\dist(f_\ve,\mcZ)^2}
    =
    1-\frac{q}{\mu_3}
    -\frac{q(q-1)}{3}\ve
    \int_{\mbBn}\mcU^{q-2}\Psi^3\diopv
    +O(\ve^2)\qquad(\ve\to0).
  \end{equation}
\end{lemma}

\begin{proof}
  The lemma extends \eqref{equ:local-eigenfunction-expansion} by one order. We first show that $\frac{\Psi}{\mcU}$ is bounded by comparing $R$ and $\mcU$ at both endpoints; this bound controls the fourth-order Taylor remainder, which lets us expand $\norm{f_\ve}_{L^{q+1}(\mbBn)}^{q+1}$ to third order and then the stability quotient.

  \emph{Step 1: endpoint behavior.} Let $R$ be a nonzero finite-energy radial factor with angular degree $\ell$ and eigenvalue $\mu$. By \eqref{equ:polar-laplacian}, its equation gives
  \begin{equation*}
    \begin{aligned}
      &\frac{\mrd^2}{\mrd r^2}
      \left[(\sinh r)^{\frac{n-1}{2}}R(r)\right]\\
      &\quad=
      \left[
      \lambda_1-\lambda+
      \frac{\ell(\ell+n-2)+\frac{(n-1)(n-3)}{4}}{\sinh^2r}
      -\mu\mcU(r)^{q-1}
      \right]
      \left[(\sinh r)^{\frac{n-1}{2}}R(r)\right].
    \end{aligned}
  \end{equation*}
  The perturbation of the positive constant $\lambda_1-\lambda$ is integrable at infinity by Lemma~\ref{lem:ground-state}. Variation of constants constructs two solutions with respective leading factors $\mre^{\sqrt{\lambda_1-\lambda}\,r}$ and $\mre^{-\sqrt{\lambda_1-\lambda}\,r}$, that is, $R$ behaves like $\mre^{-\beta_-r}$ or $\mre^{-\beta_+r}$, and their derivatives have the corresponding leading factors. One may construct the decaying solution by a Volterra equation on a sufficiently far tail and obtain the other solution by reduction of order. The finite-energy condition excludes the growing solution, so a nonzero solution has a nonzero coefficient of the decaying solution.

  The ground state obeys the same argument with $\ell=0$ and $\mu=1$, and, since $\mcU>0$, it has a strictly positive fast-tail coefficient. It follows that $R$ and $\mcU$ have the same decay exponent $\beta_+$, with nonzero leading coefficients, so their ratio has a finite limit at infinity. At the origin, regularity gives $R(r)=c r^\ell+O(r^{\ell+2})$ for a nonzero constant $c$, and the ground state is smooth and strictly positive there by Section~\ref{sec:intro} and Lemma~\ref{lem:pairing-center-derivative}. Consequently $\frac{RY}{\mcU}$ is bounded for every spherical harmonic $Y$ of degree $\ell$. This proves the first assertion; since $\mcU^{q+1}$ is integrable by \eqref{equ:ground-state-energy-normalization}, so is $\mcU^{q-2}|\Psi|^3$.

  \emph{Step 2: the expansion.} Now let $\Psi$ have eigenvalue $\mu_3$ and $\norm{\Psi}_\lambda=1$, and let $f_\ve=\mcU+\ve\Psi$. Equation~\eqref{equ:const-exact-distance} gives $\dist(f_\ve,\mcZ)^2=\ve^2$ for all sufficiently small $|\ve|$. For $|\ve|\lrnorm{\frac{\Psi}{\mcU}}_{L^\infty(\mbBn)}\leqslant\frac{1}{2}$, Taylor's formula and the eigenfunction moments in \eqref{equ:const-eigen-moments} give
  \begin{align*}
      \int_{\mbBn}|f_\ve|^{q+1}\diopv
      &=S_\lambda^{\frac{n}{2}} +\frac{q(q+1)}{2\mu_3}\ve^2 +\frac{q(q+1)(q-1)}{6}\ve^3
      \int_{\mbBn}\mcU^{q-2}\Psi^3\diopv +O(\ve^4),
  \end{align*}
  where the pointwise fourth-order remainder is bounded by $C\ve^4\mcU^{q+1}|\frac{\Psi}{\mcU}|^4$, which is integrable by \eqref{equ:ground-state-energy-normalization} and Step~1. The argument also applies when $q<2$, since the expansion is taken around $1+\frac{\ve\Psi}{\mcU}\geqslant\frac{1}{2}$.

  \emph{Step 3: the stability quotient.} The same normalization and orthogonality give $\norm{f_\ve}_\lambda^2=S_\lambda^{\frac{n}{2}}+\ve^2$. Expanding the outer power and using $S_\lambda(S_\lambda^{\frac{n}{2}})^{\frac{2}{q+1}-1}=1$ yields
  \begin{align*}
    S_\lambda\norm{f_\ve}_{L^{q+1}(\mbBn)}^2
    &=S_\lambda^{\frac{n}{2}}+\frac{q}{\mu_3}\ve^2
    +\frac{q(q-1)}{3}\ve^3 \int_{\mbBn}\mcU^{q-2}\Psi^3\diopv
    +O(\ve^4).
  \end{align*}
  Subtracting this expansion from $\norm{f_\ve}_\lambda^2$ and dividing by $\dist(f_\ve,\mcZ)^2=\ve^2$ gives \eqref{equ:strict-local-cubic}. This finishes the proof of Lemma~\ref{lem:strict-local-cubic}.
\end{proof}

\subsubsection{The sign of the cubic integral}\label{subsec:strict-local-sign}
By Lemmas~\ref{lem:strict-local-sectors} and~\ref{lem:strict-local-cubic}, it remains to find a separated eigenfunction with eigenvalue $\mu_3$ and positive cubic integral, of angular degree $2$ if $\mu_3=\mu_{2,1}$ and radial if $\mu_3=\mu_{0,2}$. We first treat degree $2$, where an explicit zonal harmonic works.

\begingroup
\begin{lemma}\label{lem:strict-local-degree2}
  There is a separated eigenfunction $\Psi=R_2(r)Y(\theta)$ of \eqref{equ:const-linearized} with eigenvalue $\mu_{2,1}$, angular factor $Y(\theta)=\theta_1^2-\frac{1}{n}$, positive radial factor $R_2$ and $\norm{\Psi}_\lambda=1$, and every such $\Psi$ satisfies
  \begin{equation*}
    \int_{\mbBn}\mcU^{q-2}\Psi^3\diopv>0.
  \end{equation*}
\end{lemma}
\endgroup

\begin{proof}
  The first radial eigenfunction in degree $2$ has no zero in $(0,\infty)$, by the oscillation theorem used in the proof of Lemma~\ref{lem:strict-local-sectors}, so we may choose its radial factor $R_2$ strictly positive. As angular factor we take
  \begin{equation*}
    Y(\theta)=\theta_1^2-\frac{1}{n},
    \qquad \theta\in\mbS^{n-1},
  \end{equation*}
  which is a spherical harmonic of degree $2$. Multiplying $R_2$ by a positive constant, we obtain $\Psi=R_2Y$ with $\norm{\Psi}_\lambda=1$ and eigenvalue $\mu_{2,1}$; this proves the existence assertion.

  Next we show that $Y$ has a positive cubic average. The spherical moments are
  \begin{equation*}
    \begin{aligned}
      \frac{1}{|\mbS^{n-1}|}
      \int_{\mbS^{n-1}}\theta_1^2\dioptheta
      &=\frac{1}{n},\\
      \frac{1}{|\mbS^{n-1}|}
      \int_{\mbS^{n-1}}\theta_1^4\dioptheta
      &=\frac{3}{n(n+2)},\\
      \frac{1}{|\mbS^{n-1}|}
      \int_{\mbS^{n-1}}\theta_1^6\dioptheta
      &=\frac{15}{n(n+2)(n+4)}.
    \end{aligned}
  \end{equation*}
  Expanding the cube gives
  \begin{equation*}
    \frac{1}{|\mbS^{n-1}|}
    \int_{\mbS^{n-1}}Y^3\dioptheta
    =
    \frac{15}{n(n+2)(n+4)}
    -\frac{9}{n^2(n+2)}
    +\frac{2}{n^3}
    =\frac{8(n-1)(n-2)}{n^3(n+2)(n+4)}>0.
  \end{equation*}
  In polar coordinates about $o$, by \eqref{equ:polar-volume}, the cubic integral of $\Psi$ factors as
  \begin{equation*}
    \int_{\mbBn}\mcU^{q-2}\Psi^3\diopv
    =
    \left[
    \int_0^\infty
    \mcU^{q-2}R_2^3(\sinh r)^{n-1}\diopr
    \right]
    \left[
    \int_{\mbS^{n-1}}Y^3\dioptheta
    \right].
  \end{equation*}
  The first factor is finite by Lemma~\ref{lem:strict-local-cubic} and is strictly positive, since $\mcU>0$ and $R_2>0$. The second factor is positive by the computation above. Hence the cubic integral is positive. This finishes the proof of Lemma~\ref{lem:strict-local-degree2}.
\end{proof}

We now turn to the radial case.

\begingroup
The positivity of the radial cubic integral rests on a one-dimensional comparison principle, which we isolate first. In Lemmas~\ref{lem:strict-local-onedim} and~\ref{lem:strict-local-radial} and their proofs, primes are used only on $w$, $h$ and $\mcV$, and denote derivatives with respect to $t$; derivatives of all other functions, in particular of functions of the height variable $y$, and derivatives of composite expressions are written with $\frac{\mrd}{\mrd y}$, $\frac{\mrd}{\mrd t}$ and $\frac{\mrd}{\mrd r}$.

\begin{lemma}\label{lem:strict-local-onedim}
  Let $\kappa>1$, $\alpha_0>0$, $\gamma>0$ and $\beta_0>0$, and let $\mcV\in C^1((\nit,0))$ be strictly decreasing with a finite limit at $\nit$. Let $w,h\in C^2((\nit,0))$ with $w>0$ satisfy
  \begin{equation}\label{equ:onedim-equations}
    -w''+\mcV w=w^\kappa,
    \qquad
    -\frac{\mrd}{\mrd t}\bigl(w^2h'\bigr)=\gamma w^{\kappa+1}h
    \qquad\text{in }(\nit,0),
  \end{equation}
  and assume the following.
  \begin{enumerate}[label=\textnormal{(\roman*)}]
    \item\label{item:onedim-w} As $t\to\nit$, $w(t)\sim\mre^{\alpha_0t}$ and $\frac{w'(t)}{w(t)}\to\alpha_0$; as $t\uparrow0$, $w(t)\sim(-t)^{\beta_0}$ and $\frac{(-t)\,w'(t)}{w(t)}\to-\beta_0$.
    \item\label{item:onedim-h} The function $h$ has limits $h(\nit)>0$ and $h(0)<0$ at the endpoints, and exactly one zero in $(\nit,0)$.
    \item\label{item:onedim-orth} $w^2h'\to0$ at both endpoints, and $\int_{\nit}^0w^{\kappa+1}h\diopt=0$.
  \end{enumerate}
  Then
  \begin{equation*}
    \int_{\nit}^0w^{\kappa+1}h^3\diopt>0.
  \end{equation*}
\end{lemma}

\begin{proof}
  We first study $w$. It has exactly one, nondegenerate, maximum (Step~1). At equal heights, its left flank is slower than its right flank, and the ratio of the speeds increases (Step~2). We then turn to $h$. It is strictly decreasing (Step~3), its zero lies strictly before the peak (Step~4), and its flux on the left flank dominates the flux on the right flank (Step~5). Finally, in Step~6, two integrations by parts express the integral as an integral over heights whose integrand is positive by Step~5. All comparison constants below may depend on $w$, $h$, $\mcV$, $\gamma$ and $\beta_0$.

  \emph{Step 1: $w$ has exactly one nondegenerate maximum.} By the first equation in \eqref{equ:onedim-equations},
  \begin{equation*}
    \frac{\mrd}{\mrd t}
    \left[
    \frac{1}{2}w'^2
    -\frac{1}{2}\mcV w^2
    +\frac{1}{\kappa+1}w^{\kappa+1}
    \right]
    =-\frac{1}{2}\mcV'w^2>0.
  \end{equation*}
  By \ref{item:onedim-w} and the finiteness of $\mcV(\nit)$, the expression in brackets tends to $0$ at the left endpoint. It is therefore strictly positive at every interior point. At a stationary point with $\mcV\geqslant0$ this gives
  \begin{equation*}
    w^{\kappa-1}>\frac{\kappa+1}{2}\mcV\geqslant\mcV,
  \end{equation*}
  and at a stationary point with $\mcV<0$ the inequality $w^{\kappa-1}>\mcV$ holds trivially. Thus every stationary point satisfies
  \begin{equation*}
    w''=(\mcV-w^{\kappa-1})w<0.
  \end{equation*}
  Since $w$ is positive and tends to $0$ at both endpoints, and two strict local maxima would enclose a local minimum, $w$ has exactly one stationary point $t_\ast$. Let
  \begin{equation*}
    M:=w(t_\ast)=\max_{(\nit,0)}w.
  \end{equation*}

  \emph{Step 2: comparison of the speeds on the two flanks.} By Step~1, for each height $y\in(0,M)$ there are exactly two points $t_-(y)<t_\ast<t_+(y)$ with
  \begin{equation*}
    w(t_-(y))=w(t_+(y))=y,
  \end{equation*}
  and $w'(t_-(y))>0>w'(t_+(y))$. We define the speeds
  \begin{equation*}
    p_\pm(y):=\bigl|w'(t_\pm(y))\bigr|>0.
  \end{equation*}
  Since $\frac{\mrd t_-}{\mrd y}=\frac{1}{p_-}$, $\frac{\mrd t_+}{\mrd y}=-\frac{1}{p_+}$ and $w''=\mcV w-w^\kappa$, we have
  \begin{equation*}
    \frac{\mrd}{\mrd y}\bigl(p_\pm^2\bigr)
    =2y\bigl(\mcV(t_\pm(y))-y^{\kappa-1}\bigr).
  \end{equation*}
  Since $\mcV$ strictly decreases and $t_-(y)<t_+(y)$,
  \begin{equation*}
    \frac{\mrd}{\mrd y}\bigl(p_+^2-p_-^2\bigr)
    =-2y\bigl(\mcV(t_-(y))-\mcV(t_+(y))\bigr)<0.
  \end{equation*}
  Both speeds tend to $0$ at the peak height. It follows that
  \begin{equation*}
    p_+>p_->0\qquad(0<y<M).
  \end{equation*}

  We also need strict monotonicity of the speed ratio. Writing $t_\pm=t_\pm(y)$, direct differentiation gives
  \begin{equation*}
    \begin{aligned}
      \frac{\mrd}{\mrd y}\log\frac{p_-}{p_+}
      =\frac{y}{p_-^2p_+^2}
      \Bigl[
      &\bigl(\mcV(t_-)-\mcV(t_+)\bigr)p_-^2\\
      &-\bigl(y^{\kappa-1}-\mcV(t_-)\bigr)
       \bigl(p_+^2-p_-^2\bigr)
      \Bigr].
    \end{aligned}
  \end{equation*}
  The function $y^{\kappa-1}-\mcV(t_-(y))$ strictly increases, and the positive function $\mcV(t_-(y))-\mcV(t_+(y))$ strictly decreases. Moreover,
  \begin{equation*}
    \begin{aligned}
      p_-(y)^2
      &=2\int_y^M
      s\bigl(s^{\kappa-1}-\mcV(t_-(s))\bigr)\diop s,\\
      p_+(y)^2-p_-(y)^2
      &=2\int_y^M
      s\bigl(\mcV(t_-(s))-\mcV(t_+(s))\bigr)
      \diop s.
    \end{aligned}
  \end{equation*}
  If $y^{\kappa-1}-\mcV(t_-(y))\leqslant0$, the bracket in the derivative formula is strictly positive. If this quantity is positive, we compare the two displayed integrals: at every $s>y$ the first integrand contains the strictly larger factor $s^{\kappa-1}-\mcV(t_-(s))$, and the potential difference in the second integrand is strictly smaller than its value at $y$. The bracket is again strictly positive. Therefore
  \begin{equation*}
    \frac{\mrd}{\mrd y}\log\frac{p_-}{p_+}>0
    \qquad(0<y<M).
  \end{equation*}

  \emph{Step 3: $h$ is strictly decreasing.} The function $h$ is bounded by \ref{item:onedim-h}, and $w^{\kappa+1}$ is integrable by \ref{item:onedim-w}. Integrating the second equation in \eqref{equ:onedim-equations} from $\nit$, and using $w^2h'\to0$ at $\nit$ and the orthogonality in \ref{item:onedim-orth}, we obtain
  \begin{equation*}
    -w(t)^2h'(t)
    =\gamma\int_{\nit}^tw(\tau)^{\kappa+1}h(\tau)\dioptau
    =-\gamma\int_t^0w(\tau)^{\kappa+1}h(\tau)\dioptau.
  \end{equation*}
  By \ref{item:onedim-h}, $h>0$ before its zero and $h<0$ after it. Before the zero the first integral is strictly positive, and after the zero the second one is. Therefore $h'<0$ throughout $(\nit,0)$. The same formula, \ref{item:onedim-w} and \ref{item:onedim-h} give
  \begin{equation}\label{equ:onedim-flux-rates}
    \begin{aligned}
      -h'(t)&\sim w(t)^{\kappa-1}
      &&(t\to\nit),\\
      -h'(t)&\sim(-t)^{\beta_0(\kappa-1)+1}
      &&(t\uparrow0).
    \end{aligned}
  \end{equation}

  \emph{Step 4: the zero of $h$ lies before the peak.} Define
  \begin{equation*}
    h_-(y):=h(t_-(y)),\qquad
    h_+(y):=-h(t_+(y)).
  \end{equation*}
  By Step~3 both functions strictly decrease. The chain rule gives
  \begin{equation*}
    \frac{\mrd h_\pm}{\mrd y}(y)=\frac{h'(t_\pm(y))}{p_\pm(y)},
  \end{equation*}
  and the second equation in \eqref{equ:onedim-equations} becomes
  \begin{equation}\label{equ:onedim-height-equation}
    -\frac{\mrd}{\mrd y}
    \left(y^2p_\pm\frac{\mrd h_\pm}{\mrd y}\right)
    =\gamma\frac{y^{\kappa+1}}{p_\pm}h_\pm.
  \end{equation}
  Suppose, to the contrary, that $h(t_\ast)\geqslant0$. Then $h_->0$ on $(0,M)$. Also $h_+(0+)=-h(0)>0$ and $h_+(M)=-h(t_\ast)\leqslant0$. Let $y_0\in(0,M]$ be the unique zero of $h_+$; then $h_+>0$ on $(0,y_0)$. We shall test the left-flank equation, whose positive solution is $h_-$, with $h_+$ and use $p_-<p_+$; first we justify the comparison at the endpoints.

  As $y\downarrow0$, \ref{item:onedim-w}, \ref{item:onedim-h} and \eqref{equ:onedim-flux-rates} give
  \begin{equation*}
    p_-(y)\sim y,\qquad
    h_-(y)\longrightarrow h(\nit)>0,\qquad
    \frac{\mrd h_-}{\mrd y}(y)=O\bigl(y^{\kappa-2}\bigr),
  \end{equation*}
  and
  \begin{equation*}
    h_+(y)=O(1),\qquad
    \frac{\mrd h_+}{\mrd y}(y)=O\Bigl(y^{\kappa-2+\frac{2}{\beta_0}}\Bigr).
  \end{equation*}
  Thus the transplanted derivative energy density satisfies
  \begin{equation*}
    y^2p_-\left(\frac{\mrd h_+}{\mrd y}\right)^2
    =O\Bigl(y^{2\kappa-1+\frac{4}{\beta_0}}\Bigr),
  \end{equation*}
  and the mass density $\gamma\frac{y^{\kappa+1}}{p_-}h_+^2$ is $O(y^\kappa)$; both are integrable. The Picone boundary term satisfies
  \begin{equation*}
    \frac{y^2p_-}{h_-}\frac{\mrd h_-}{\mrd y}h_+^2=O\bigl(y^{\kappa+1}\bigr)\longrightarrow0,
  \end{equation*}
  and the right-flank boundary term $y^2p_+\frac{\mrd h_+}{\mrd y}h_+=y^2h'(t_+(y))h_+(y)$ also tends to $0$, by \eqref{equ:onedim-flux-rates}.

  If $y_0<M$, then $h_+(y_0)=0$ and $h_-(y_0)>0$; the endpoint is regular, and the comparison boundary terms equal $0$. If $y_0=M$, then $h(t_\ast)=0$, and the nondegenerate peak and $h'(t_\ast)<0$ give
  \begin{equation*}
    p_\pm(y)\sim(M-y)^{\frac{1}{2}},
    \qquad
    h_\pm(y)\sim(M-y)^{\frac{1}{2}}
    \qquad(y\uparrow M).
  \end{equation*}
  The derivative energy densities are $O((M-y)^{-\frac{1}{2}})$, hence integrable, and both comparison boundary terms are $O((M-y)^{\frac{1}{2}})$ and tend to $0$.

  We may therefore integrate \eqref{equ:onedim-height-equation} for $h_+$, multiplied by $h_+$, and the Picone identity on truncated intervals and pass to the endpoints. This gives
  \begin{equation*}
    \int_0^{y_0}
    \left(y^2p_+\left(\frac{\mrd h_+}{\mrd y}\right)^2-\gamma\frac{y^{\kappa+1}}{p_+}h_+^2\right)
    \diopy=0
  \end{equation*}
  and
  \begin{equation*}
    \begin{aligned}
      &\int_0^{y_0}
      \left(y^2p_-\left(\frac{\mrd h_+}{\mrd y}\right)^2-\gamma\frac{y^{\kappa+1}}{p_-}h_+^2\right)
      \diopy\\
      &\qquad=
      \int_0^{y_0}
      y^2p_-h_-^2
      \left(\frac{\mrd}{\mrd y}\frac{h_+}{h_-}\right)^2
      \diopy
      \geqslant0.
    \end{aligned}
  \end{equation*}
  The Picone identity follows by expanding the square and using \eqref{equ:onedim-height-equation} for $h_-$; its boundary term is the one estimated above. Since $p_-<p_+$, the integral on the left of the second identity is strictly smaller than the integral in the first identity; the mass term alone makes the comparison strict. This is a contradiction. We have proved
  \begin{equation*}
    h(t_\ast)<0.
  \end{equation*}
  Thus the unique zero of $h$ lies strictly before the peak, and $h_+>0$ on $(0,M]$.

  \emph{Step 5: comparison of the fluxes.} Define the height fluxes
  \begin{equation*}
    F_\pm(y):=-y^2h'(t_\pm(y))
    =-y^2p_\pm(y)\frac{\mrd h_\pm}{\mrd y}(y),
  \end{equation*}
  which are positive by Step~3. By \eqref{equ:onedim-height-equation} and Step~4,
  \begin{equation*}
    \frac{\mrd F_\pm}{\mrd y}=\gamma\frac{y^{\kappa+1}}{p_\pm}h_\pm,
    \qquad
    \frac{\mrd F_+}{\mrd y}>0,
  \end{equation*}
  and eliminating $h_\pm$ by means of $\frac{\mrd h_\pm}{\mrd y}=-\frac{F_\pm}{y^2p_\pm}$ gives
  \begin{equation*}
    \frac{\mrd^2F_\pm}{\mrd y^2}
    +\left(\frac{1}{p_\pm}\frac{\mrd p_\pm}{\mrd y}-\frac{\kappa+1}{y}\right)\frac{\mrd F_\pm}{\mrd y}
    +\gamma\frac{y^{\kappa-1}}{p_\pm^2}F_\pm=0.
  \end{equation*}
  By \ref{item:onedim-w} and \eqref{equ:onedim-flux-rates},
  \begin{equation*}
    F_-(y)\sim y^{\kappa+1},\qquad
    F_+(y)\sim y^{\kappa+1+\frac{1}{\beta_0}}
    \qquad(y\downarrow0),
  \end{equation*}
  and at the peak height
  \begin{equation*}
    F_\pm(y)\longrightarrow
    -M^2h'(t_\ast)>0
    \qquad(y\uparrow M).
  \end{equation*}
  Hence
  \begin{equation*}
    \frac{F_-}{F_+}\longrightarrow\pit
    \quad(y\downarrow0),\qquad
    \frac{F_-}{F_+}\longrightarrow1
    \quad(y\uparrow M).
  \end{equation*}

  At an interior stationary point of this ratio, the two flux equations give
  \begin{equation*}
    \frac{\mrd^2}{\mrd y^2}\frac{F_-}{F_+}
    =-\frac{F_-}{F_+}
    \left[
    \left(\frac{\mrd}{\mrd y}\log\frac{p_-}{p_+}\right)
    \frac{1}{F_+}\frac{\mrd F_+}{\mrd y}
    +\gamma y^{\kappa-1}
    \left(\frac{1}{p_-^2}-\frac{1}{p_+^2}\right)
    \right]<0,
  \end{equation*}
  by Step~2 and $\frac{\mrd F_+}{\mrd y}>0$. The ratio therefore has no interior local minimum, and its endpoint limits imply
  \begin{equation*}
    F_-(y)>F_+(y)\qquad(0<y<M).
  \end{equation*}
  Indeed, an interior value below $1$ would force an interior minimum, and an interior value equal to $1$, with no smaller value, would itself be a local minimum. Consequently
  \begin{equation}\label{equ:onedim-flux-comparison}
    -h'(t_-(y))>-h'(t_+(y))>0
    \qquad(0<y<M).
  \end{equation}

  \emph{Step 6: the cubic identity.} The boundary term $w^2h'h^2$ tends to $0$ at both endpoints by \ref{item:onedim-h} and \ref{item:onedim-orth}. By \ref{item:onedim-w} and \eqref{equ:onedim-flux-rates}, the boundary term $w^{3-\kappa}(-h')^3$ also tends to $0$:
  \begin{equation*}
    w^{3-\kappa}(-h')^3
    =
    \begin{cases}
      O\bigl(w^{2\kappa}\bigr),&t\to\nit,\\
      O\bigl((-t)^{2\beta_0\kappa+3}\bigr),&t\uparrow0.
    \end{cases}
  \end{equation*}
  We multiply the second equation in \eqref{equ:onedim-equations} by $h^2$ and integrate by parts, then use $\frac{\mrd}{\mrd t}\bigl(w^2(-h')\bigr)=\gamma w^{\kappa+1}h$ and integrate by parts once more. This yields
  \begin{equation*}
    \begin{aligned}
      \gamma\int_{\nit}^0w^{\kappa+1}h^3\diopt
      &=2\int_{\nit}^0
      w^2h(h')^2\diopt\\
      &=\frac{2}{3\gamma}\int_{\nit}^0
      w^{-\kappa-3}\frac{\mrd}{\mrd t}\bigl[(-w^2h')^3\bigr]
      \diopt\\
      &=\frac{2(\kappa+3)}{3\gamma}\int_{\nit}^0
      w^{2-\kappa}w'(-h')^3\diopt.
    \end{aligned}
  \end{equation*}
  We change to the height variable $y=w(t)$ on each flank; $y$ runs from $0$ to $M$ on the left flank and from $M$ to $0$ on the right flank. This gives
  \begin{equation*}
    \begin{aligned}
      \int_{\nit}^0w^{\kappa+1}h^3\diopt
      =\frac{2(\kappa+3)}{3\gamma^2}
      \int_0^My^{2-\kappa}
      \Bigl[
      &\bigl(-h'(t_-(y))\bigr)^3\\
      &-\bigl(-h'(t_+(y))\bigr)^3
      \Bigr]\diopy.
    \end{aligned}
  \end{equation*}
  By \eqref{equ:onedim-flux-comparison}, the integrand is strictly positive at every interior height. It is integrable at height $0$, where the left contribution is $O(y^{2\kappa-1})$ and the right contribution has a higher power, and it is bounded at the peak height. Therefore the integral is strictly positive. This finishes the proof of Lemma~\ref{lem:strict-local-onedim}.
\end{proof}

\begin{lemma}\label{lem:strict-local-radial}
  Let $\Psi_{0,2}$ be a radial eigenfunction of \eqref{equ:const-linearized} with eigenvalue $\mu_{0,2}$, normalized by $\norm{\Psi_{0,2}}_\lambda=1$ and $\Psi_{0,2}(0)>0$. Then
  \begin{equation*}
    \int_{\mbBn}\mcU^{q-2}\Psi_{0,2}^3\diopv>0.
  \end{equation*}
\end{lemma}

\begin{proof}
  We apply the conformal change of variables and pass to the Emden--Fowler variable, in which the ground state and the ratio $\frac{\Psi_{0,2}}{\mcU}$ satisfy the equations \eqref{equ:onedim-equations} with a strictly decreasing potential (Step~1), and then verify the hypotheses of Lemma~\ref{lem:strict-local-onedim} (Step~2).

  \emph{Step 1: reduction to one dimension.} Let $r=d(x,o)$, so that $|x|=\tanh\frac{r}{2}$, and let
  \begin{equation*}
    t:=\log|x|=\log\tanh\frac{r}{2}\in(\nit,0).
  \end{equation*}
  Then $\rho(x)|x|=\sinh r$, and
  \begin{equation*}
    \frac{\mrd t}{\mrd r}=\frac{1}{\sinh r},
    \qquad
    \sinh r=-\frac{1}{\sinh t}.
  \end{equation*}
  By the conformal change of variables of Section~\ref{subsec:conformal-dictionary}, as in \eqref{equ:psdec-conformal}, the radial functions $v:=\rho^{\frac{n-2}{2}}\mcU$ and $\vp:=\rho^{\frac{n-2}{2}}\Psi_{0,2}$ satisfy
  \begin{equation*}
    -\Delta v=\Lambda\rho^2v+v^q,
    \qquad
    -\Delta\vp=\Lambda\rho^2\vp+\mu_{0,2}v^{q-1}\vp
    \qquad\text{in }B_1\setminus\{0\},
  \end{equation*}
  where $\Lambda=\lambda-\lambda_0>0$. Writing radial functions as functions of $|x|$ or of $r$, we define, within this proof, the Emden--Fowler transforms
  \begin{equation*}
    w(t):=\mre^{\frac{n-2}{2}t}v(\mre^t)=(\sinh r)^{\frac{n-2}{2}}\mcU(r),
    \qquad
    h(t):=\frac{\vp(\mre^t)}{v(\mre^t)}=\frac{\Psi_{0,2}(r)}{\mcU(r)}.
  \end{equation*}
  For a radial function $g(|x|)=|x|^{-\frac{n-2}{2}}G(\log|x|)$ one has
  \begin{equation*}
    -\Delta g=|x|^{-\frac{n+2}{2}}
    \left(-\frac{\mrd^2G}{\mrd t^2}+\frac{(n-2)^2}{4}G\right).
  \end{equation*}
  Since $\rho^2|x|^2=\frac{1}{\sinh^2t}$ and $\frac{n-2}{2}q=\frac{n+2}{2}$, the equation for $v$ becomes the first equation in \eqref{equ:onedim-equations} with $\kappa=q$ and
  \begin{equation*}
    \mcV(t):=\frac{(n-2)^2}{4}-\frac{\Lambda}{\sinh^2t}.
  \end{equation*}
  The constant $\frac{(n-2)^2}{4}$ comes from the Emden--Fowler substitution and the term $-\frac{\Lambda}{\sinh^2t}$ from $\Lambda\rho^2$. The equation for $\vp=|x|^{-\frac{n-2}{2}}wh$ becomes $-\frac{\mrd^2}{\mrd t^2}(wh)+\mcV wh=\mu_{0,2}w^qh$; subtracting $h$ times the equation for $w$ and multiplying by $w$ gives the second equation in \eqref{equ:onedim-equations} with $\gamma=\mu_{0,2}-1$. The potential satisfies
  \begin{equation*}
    \mcV'(t)=2\Lambda\frac{\cosh t}{\sinh^3t}<0,
    \qquad
    \mcV(t)\longrightarrow\frac{(n-2)^2}{4}\quad(t\to\nit),
  \end{equation*}
  so it is strictly decreasing with a finite limit at $\nit$. Finally, by \eqref{equ:polar-volume}, $\mcU^{q-2}\Psi_{0,2}^3=\mcU^{q+1}h^3$, $(\sinh r)^{n-1}\diopr=(\sinh r)^n\diopt$ and $\frac{n-2}{2}(q+1)=n$, we have
  \begin{equation}\label{equ:radial-cubic-onedim}
    \int_{\mbBn}\mcU^{q-2}\Psi_{0,2}^3\diopv
    =|\mbS^{n-1}|\int_{\nit}^0w^{q+1}h^3\diopt.
  \end{equation}

  \emph{Step 2: the hypotheses of Lemma~\ref{lem:strict-local-onedim}.} By Lemma~\ref{lem:strict-local-sectors} and Section~\ref{subsec:linearized-spectrum}, $\mu_{0,2}\geqslant\mu_3>q>1$, so that $\gamma=\mu_{0,2}-1>0$ and $\kappa=q>1$. Let
  \begin{equation*}
    \beta_0:=\beta_+-\frac{n-2}{2}=\frac{1}{2}+\sqrt{\lambda_1-\lambda}>0.
  \end{equation*}
  As $t\to\nit$, that is, as $r\to0$, the ground state is smooth and strictly positive at the origin, and $\sinh r\sim\mre^t$. As $t\uparrow0$, that is, as $r\to\infty$, the fast-tail argument in the proof of Lemma~\ref{lem:strict-local-cubic} gives $\mcU(r)\sim\mre^{-\beta_+r}$ and $\frac{1}{\mcU}\frac{\mrd\mcU}{\mrd r}\to-\beta_+$, while $-t\sim\mre^{-r}$. Since
  \begin{equation*}
    \frac{w'}{w}=\sinh r\left(\frac{n-2}{2}\coth r+\frac{1}{\mcU}\frac{\mrd\mcU}{\mrd r}\right),
  \end{equation*}
  this gives hypothesis \ref{item:onedim-w} with $\alpha_0=\frac{n-2}{2}$:
  \begin{equation*}
    \begin{aligned}
      w(t)&\sim\mre^{\frac{n-2}{2}t},
      &\frac{w'(t)}{w(t)}&\longrightarrow\frac{n-2}{2}
      &&(t\to\nit),\\
      w(t)&\sim(-t)^{\beta_0},
      &\frac{(-t)\,w'(t)}{w(t)}&\longrightarrow-\beta_0
      &&(t\uparrow0).
    \end{aligned}
  \end{equation*}
  The oscillation theorem used in the proof of Lemma~\ref{lem:strict-local-sectors} shows that $\Psi_{0,2}$ has exactly one zero in $(0,\infty)$, which is simple. By the proof of Lemma~\ref{lem:strict-local-cubic}, $\Psi_{0,2}$ and $\mcU$ have nonzero fast-tail coefficients, so $h$ has a nonzero limit as $r\to\infty$; moreover $h(\nit)=\frac{\Psi_{0,2}(0)}{\mcU(0)}>0$. Hence $h(0)<0$, and hypothesis \ref{item:onedim-h} holds. Next,
  \begin{equation*}
    w^2h'=(\sinh r)^{n-1}
    \left(\mcU\frac{\mrd\Psi_{0,2}}{\mrd r}-\Psi_{0,2}\frac{\mrd\mcU}{\mrd r}\right)
  \end{equation*}
  tends to $0$ as $r\to0$, and also as $r\to\infty$ by the same endpoint asymptotics, since $2\beta_+>n-1$. The ground state is the eigenfunction of \eqref{equ:const-linearized} with eigenvalue $\mu_{0,1}=1\neq\mu_{0,2}$, so $\Psi_{0,2}$ is orthogonal to it in the weighted inner product, and the change of variables of Step~1 gives
  \begin{equation*}
    0=\int_{\mbBn}\mcU^{q-1}\cdot\mcU\,\Psi_{0,2}\diopv
    =|\mbS^{n-1}|\int_{\nit}^0w^{q+1}h\diopt.
  \end{equation*}
  This is hypothesis \ref{item:onedim-orth}.

  Lemma~\ref{lem:strict-local-onedim} and \eqref{equ:radial-cubic-onedim} now give $\int_{\mbBn}\mcU^{q-2}\Psi_{0,2}^3\diopv>0$. This finishes the proof of Lemma~\ref{lem:strict-local-radial}.
\end{proof}
\endgroup

\begin{proof}[Proof of Proposition~\ref{prop:strict-local}]
  By Lemma~\ref{lem:strict-local-cubic}, it suffices to find a separated eigenfunction $\Psi$ with eigenvalue $\mu_3$, normalized by $\norm{\Psi}_\lambda=1$, with $\int_{\mbBn}\mcU^{q-2}\Psi^3\diopv>0$. By Lemma~\ref{lem:strict-local-sectors}, $\mu_3=\min\{\mu_{0,2},\mu_{2,1}\}$, so we distinguish the two cases $\mu_{0,2}\leqslant\mu_{2,1}$ and $\mu_{2,1}<\mu_{0,2}$; in the first case we take a radial eigenfunction, and in the second one of angular degree $2$. Finally we insert $\Psi$ into \eqref{equ:strict-local-cubic}.

  \emph{Case $\mu_{0,2}\leqslant\mu_{2,1}$.} Then Lemma~\ref{lem:strict-local-sectors} gives $\mu_3=\mu_{0,2}$, so we may take $\Psi=\Psi_{0,2}$ in \eqref{equ:strict-local-cubic}. By Lemma~\ref{lem:strict-local-radial}, $\int_{\mbBn}\mcU^{q-2}\Psi_{0,2}^3\diopv>0$.

  \emph{Case $\mu_{2,1}<\mu_{0,2}$.} Then Lemma~\ref{lem:strict-local-sectors} gives $\mu_3=\mu_{2,1}$, and we take for $\Psi$ the eigenfunction of Lemma~\ref{lem:strict-local-degree2}, whose cubic integral is positive.

  \emph{Conclusion.} In either case, $\Psi$ satisfies the hypotheses of Lemma~\ref{lem:strict-local-cubic} and satisfies $\int_{\mbBn}\mcU^{q-2}\Psi^3\diopv>0$. Choosing $\ve>0$ sufficiently small in \eqref{equ:strict-local-cubic}, the resulting quotient is strictly smaller than $1-\frac{q}{\mu_3}$, which equals $c_{\mrloc}$ by Theorem~\ref{thm:crit-local}. Its denominator equals $\ve^2>0$, so the test function is nonzero and lies outside $\mcZ$, and \eqref{equ:cbe-def} gives $c_{\mathrm{BE}}<c_{\mrloc}$. This finishes the proof of Proposition~\ref{prop:strict-local}.
\end{proof}

\subsection{Proof of Theorem~\ref{thm:crit-global-bounds}}\label{subsec:proof-global}

\begin{proof}[Proof of Theorem~\ref{thm:crit-global-bounds}]
  The strict upper bounds are Propositions~\ref{prop:strict-conc}, \ref{prop:strict-dich} and~\ref{prop:strict-local}.
  \begingroup
  For the lower bound, write $A:=\sup_{0<s<1}s\nu(s)$ and $B:=\frac{1}{2}c_{\mathrm{dich}}$. Proposition~\ref{prop:weak-global-lower-bound} gives
  \begin{equation*}
    c_{\mathrm{BE}}\geqslant\min\{A,c_{\mathrm{conc}},B\}>0.
  \end{equation*}
  If $\min\{A,B\}\geqslant c_{\mathrm{conc}}$, then $\min\{A,c_{\mathrm{conc}},B\}=c_{\mathrm{conc}}$, and hence $c_{\mathrm{BE}}\geqslant c_{\mathrm{conc}}$, which contradicts Proposition~\ref{prop:strict-conc}. Therefore $\min\{A,B\}<c_{\mathrm{conc}}$, so $\min\{A,c_{\mathrm{conc}},B\}=\min\{A,B\}$. This is the left-hand side of \eqref{equ:intro-global-two-sided}, which is thus positive and at most $c_{\mathrm{BE}}$.\endgroup\par
  This finishes the proof of Theorem~\ref{thm:crit-global-bounds}.
\end{proof}

\begingroup
\begin{remark}\label{rem:channels}
  The proof shows how the lower and upper bounds in \eqref{equ:intro-global-two-sided} match channel by channel.
  In the local channel, the lower bound $\sup_{0<s<1}s\nu(s)$ is the local estimate transferred along the flow, which costs the factor $\sigma=s$ at the first entrance time, and the upper bound is $c_{\mrloc}$.
  In the concentration channel, Section~\ref{sec:crit-global-bounds-proof} gives $c_{\mathrm{conc}}$, the elementary bound $1-\frac{S_\lambda}{\norm{f}_\lambda^2}$ at the energy $S_n$. It equals the upper bound, because $\sigma\to1$ along the concentration family and the elementary bound is an equality when $\sigma=1$; it can be removed from the lower bound precisely because the upper bound is strict.
  In the dichotomy channel, the lower bound $\frac{1}{2}c_{\mathrm{dich}}$ is the elementary bound at the energy $2^{\frac{2}{n}}S_\lambda$, and the upper bound is $c_{\mathrm{dich}}$. The factor $2$ between them reflects $\sigma\to\frac{1}{2}$ along the dichotomy family. Since the lower bound is only half of the upper bound, the strict upper bound does not remove this term.
\end{remark}\endgroup\par

\section{Existence of extremals}\label{sec:attainment}
This section proves Theorem~\ref{thm:attainment}. Using the strict upper bounds of Theorem~\ref{thm:crit-global-bounds}, we show that every minimizing sequence normalized in $L^{2^{\ast}}$ converges strongly after isometries to a minimizer. The same strategy, a refinement of an argument going back to Brezis and Lieb \cite{BrezisLieb1983} in which the strict upper bounds rule out every loss of compactness, has been used by K\"onig \cite{Koenig2025} in the Euclidean case. There, dilations remove concentration at a point; on $\mbBn$ concentration cannot be removed by isometries and is excluded instead by the strict concentration bound.

For $h\in H^1(\mbBn)$, let
\begin{equation*}
  m(h):=S_\lambda^{-\frac{n}{4}}
  \sup_{\xi\in\mbBn}|\lag h,\mcU_\xi\rag_\lambda|.
\end{equation*}
Lemmas~\ref{lem:dist-formula}, \ref{lem:duality}
and~\ref{lem:ground-state-normalization} give
\begin{equation}\label{equ:minimizing-pairing}
  \dist(h,\mcZ)^2=\norm{h}_\lambda^2-m(h)^2,\quad
  m(h)^2\leqslant S_\lambda\norm{h}_{L^{2^{\ast}}(\mbBn)}^2.
\end{equation}
Lemma~\ref{lem:dist-positive} gives positive distance
for every $h\notin\mcZ\cup\{0\}$.
Since scalar multiples of $\mcU$ approach $0$, this also gives
$\ol{\mcZ}=\mcZ\cup\{0\}$.
The definition of $c_{\mathrm{BE}}$ and \eqref{equ:minimizing-pairing} imply
\begin{equation}\label{equ:minimizing-global}
  (1-c_{\mathrm{BE}})\norm{h}_\lambda^2+c_{\mathrm{BE}} m(h)^2
  \geqslant S_\lambda\norm{h}_{L^{2^{\ast}}(\mbBn)}^2
  \qquad(h\in H^1(\mbBn)).
\end{equation}
For $h\in\mcZ\cup\{0\}$, this inequality is an equality.

\begin{proposition}\label{prop:minimizing-compact}
  Let $\{f_k\}\subset H^1(\mbBn)\setminus(\mcZ\cup\{0\})$ satisfy
  $\norm{f_k}_{L^{2^{\ast}}(\mbBn)}=1$ and
  \begin{equation*}
    \frac{\delta(f_k)^2}{\dist(f_k,\mcZ)^2}\longrightarrow c_{\mathrm{BE}}.
  \end{equation*}
  Then there are isometries $\varphi_k$ and a subsequence such that
  $f_k\circ\varphi_k\to f_\infty$ strongly in $H^1(\mbBn)$.
  Moreover $f_\infty\notin\mcZ\cup\{0\}$, $\norm{f_\infty}_{L^{2^{\ast}}(\mbBn)}=1$, and
  \begin{equation*}
    \delta(f_\infty)^2
    =c_{\mathrm{BE}}\dist(f_\infty,\mcZ)^2.
  \end{equation*}
\end{proposition}

\begin{proof}
  We first show that $\{f_k\}$ is bounded and record the asymptotic identity \eqref{equ:minimizing-equality}. Then we exclude vanishing by the strict concentration bound, which gives a nonzero weak limit after isometries (Step~1); next we exclude a nonzero remainder by the strict dichotomy bound, which gives strong convergence (Step~2); and finally we exclude a limit in $\mcZ$ by the strict local bound (Step~3). Throughout we use that composition with an isometry preserves $\norm{\cdot}_\lambda$, $\norm{\cdot}_{L^{2^{\ast}}(\mbBn)}$ and the family $\{\mcU_\xi\}$, hence also $\delta$, $\dist(\cdot,\mcZ)$ and $m$.

  \emph{Preliminaries.} Theorem~\ref{thm:crit-global-bounds} gives $0\leqslant c_{\mathrm{BE}}<c_{\mathrm{dich}}<1$. Since $\dist(f_k,\mcZ)^2\leqslant\norm{f_k}_\lambda^2$, the minimizing property implies
  \begin{equation*}
    (1-c_{\mathrm{BE}}+o(1))\norm{f_k}_\lambda^2\leqslant S_\lambda,
  \end{equation*}
  and, as $c_{\mathrm{BE}}<1$, $\{f_k\}$ is bounded in $H^1(\mbBn)$.
  The distance is bounded as well, so
  \begin{equation}\label{equ:minimizing-equality}
    (1-c_{\mathrm{BE}})\norm{f_k}_\lambda^2+c_{\mathrm{BE}} m(f_k)^2
    =S_\lambda+o(1).
  \end{equation}

  \emph{Step 1: a nonzero weak limit after isometries.} We argue by contradiction and suppose, along a subsequence, that
  \begin{equation*}
    \sup_{y\in\mbBn}\int_{B_{\mbH}(y,1)}|f_k|^2\diopv
    \longrightarrow0.
  \end{equation*}
  We show that then the stability quotient of $f_k$ is asymptotically at least $c_{\mathrm{conc}}$. This is done in three parts: $f_k\to0$ in $L^{p_0}$, the energy is asymptotically at least $S_n$, and $m(f_k)\to0$.

  First, $f_k\to0$ in $L^{p_0}$ with $p_0:=2+\frac{4}{n}$. Use the covering by unit balls $B_{\mbH}(y_\alpha,1)$ from
  Section~\ref{subsec:psdec}; its doubled balls have uniformly bounded overlap.
  H\"older's inequality and the local Sobolev inequality give
  \begin{equation*}
    \begin{aligned}
      \int_{B_{\mbH}(y_\alpha,1)}|f_k|^{p_0}\diopv
      &\leqslant
      \left(\int_{B_{\mbH}(y_\alpha,1)}|f_k|^2\diopv\right)^{\frac{2}{n}}
      \left(\int_{B_{\mbH}(y_\alpha,1)}|f_k|^{2^{\ast}}\diopv
      \right)^{\frac{2}{2^{\ast}}}\\
      &\leqslant
      C\left(\sup_{y\in\mbBn}
       \int_{B_{\mbH}(y,1)}|f_k|^2\diopv\right)^{\frac{2}{n}}
      \int_{B_{\mbH}(y_\alpha,2)}
       \bigl(|\nabla_{\mbBn}f_k|^2+|f_k|^2\bigr)\diopv ,
    \end{aligned}
  \end{equation*}
  The constant is independent of $\alpha$ by homogeneity.
  Summing over the covering yields
  \begin{equation*}
    \norm{f_k}_{L^{p_0}(\mbBn)}^{p_0}
    \leqslant
    C\left(\sup_{y\in\mbBn}
       \int_{B_{\mbH}(y,1)}|f_k|^2\diopv\right)^{\frac{2}{n}}
    \norm{f_k}_{H^1(\mbBn)}^2
    \longrightarrow0.
  \end{equation*}
  Interpolation with the bounded $L^2$ and $L^{2^{\ast}}$ norms
  gives convergence to $0$ in every $L^s$ with
  $2<s<2^{\ast}$.

  Second, we show that $\norm{f_k}_\lambda^2\geqslant S_n+(\lambda_1-\lambda)\norm{f_k}_{L^2(\mbBn)}^2-o(1)$. To this end we split $|f_k|$ into its large and small parts,
  \begin{equation*}
    h_k:=(|f_k|-1)_+,\qquad
    \ell_k:=\min\{|f_k|,1\},\qquad |f_k|=h_k+\ell_k.
  \end{equation*}
  Sobolev truncation gives $h_k,\ell_k\in H^1(\mbBn)$ and
  \begin{equation*}
    |\nabla_{\mbBn}f_k|^2
    =|\nabla_{\mbBn}h_k|^2+|\nabla_{\mbBn}\ell_k|^2
    \quad\text{almost everywhere}.
  \end{equation*}
  Moreover,
  \begin{equation*}
    \begin{aligned}
      \norm{h_k}_{L^2(\mbBn)}^2&\leqslant\norm{f_k}_{L^{p_0}(\mbBn)}^{p_0}\longrightarrow0,\\
      \norm{\ell_k}_{L^{2^{\ast}}(\mbBn)}^{2^{\ast}}
      &\leqslant\norm{f_k}_{L^{p_0}(\mbBn)}^{p_0}\longrightarrow0,\\
      0\leqslant
      \norm{f_k}_{L^2(\mbBn)}^2-\norm{\ell_k}_{L^2(\mbBn)}^2
      &\leqslant\norm{f_k}_{L^{p_0}(\mbBn)}^{p_0}\longrightarrow0.
    \end{aligned}
  \end{equation*}
  Since $|f_k|=h_k+\ell_k$, we have
  $\norm{h_k}_{L^{2^{\ast}}(\mbBn)}\to1$.
  Applying the conformal change of variables \eqref{equ:dictionary} and the
  Euclidean Sobolev inequality \eqref{equ:crit-Sn} to $h_k$, the spectral lower bound $\lambda_1$ to $\ell_k$, and the gradient identity, we find
  \begin{equation*}
    \begin{aligned}
      \norm{f_k}_\lambda^2
      &\geqslant
      S_n\norm{h_k}_{L^{2^{\ast}}(\mbBn)}^2
      +\lambda_0\norm{h_k}_{L^2(\mbBn)}^2
      +\lambda_1\norm{\ell_k}_{L^2(\mbBn)}^2
      -\lambda\norm{f_k}_{L^2(\mbBn)}^2\\
      &\geqslant
      S_n+(\lambda_1-\lambda)\norm{f_k}_{L^2(\mbBn)}^2-o(1).
    \end{aligned}
  \end{equation*}
  In particular, $\norm{f_k}_\lambda^2\geqslant S_n-o(1)$.

  Third, we show that $m(f_k)\to0$. The local $L^2$ condition is preserved by every isometric pullback, and each fixed compact set is covered by finitely many unit balls, so every weak limit of any such pullback sequence vanishes by local Rellich compactness.
  If $m(f_k)$ did not tend to $0$, we could select centers $\xi_k$ such that $|\lag f_k,\mcU_{\xi_k}\rag_\lambda|$ stays bounded away from $0$ along a subsequence. Moving these centers to $o$ would contradict weak convergence against $\mcU$.
  Hence $m(f_k)\to0$.

  Combining the three parts with \eqref{equ:minimizing-pairing},
  \begin{equation*}
    \frac{\delta(f_k)^2}{\dist(f_k,\mcZ)^2}
    =1-\frac{S_\lambda}{\norm{f_k}_\lambda^2}+o(1)
    \geqslant1-\frac{S_\lambda}{S_n}-o(1).
  \end{equation*}
  Letting $k\to\infty$ gives $c_{\mathrm{BE}}\geqslant c_{\mathrm{conc}}$, which contradicts the strict concentration bound $c_{\mathrm{BE}}<c_{\mathrm{conc}}$ of Theorem~\ref{thm:crit-global-bounds}.

  Consequently, after taking a subsequence, some unit balls carry
  a fixed positive amount of $L^2$ mass.
  Choose their centers $y_k$ and isometries with $\varphi_k(o)=y_k$.
  After these pullbacks, which do not change the hypotheses, we retain the notation $f_k$.
  Weak compactness and local Rellich compactness give
  \begin{equation*}
    f_k\rightharpoonup f_\infty\neq0
    \quad\text{in }H^1(\mbBn).
  \end{equation*}
  A diagonal subsequence also converges almost everywhere.
  Let $r_k:=f_k-f_\infty$, so that $r_k\rightharpoonup0$ in $H^1(\mbBn)$.
  Hilbert space orthogonality and the Brezis--Lieb lemma
  \cite[Theorem~1]{BrezisLieb1983} give
  \begin{equation}\label{equ:minimizing-weak-splitting}
    \begin{aligned}
      \norm{f_k}_\lambda^2
      &=\norm{f_\infty}_\lambda^2+\norm{r_k}_\lambda^2+o(1),\\
      1&=\norm{f_\infty}_{L^{2^{\ast}}(\mbBn)}^{2^{\ast}}
       +\norm{r_k}_{L^{2^{\ast}}(\mbBn)}^{2^{\ast}}+o(1).
    \end{aligned}
  \end{equation}

  \emph{Step 2: the remainder vanishes.} We first show that the pairing term also splits, namely
  \begin{equation}\label{equ:minimizing-max-pairing}
    m(f_k)^2
    =\max\{m(f_\infty)^2,m(r_k)^2\}+o(1).
  \end{equation}
  The point is that the pairings of $f_\infty$ and $r_k$ with $\mcU_\xi$ live on disjoint sets of centers asymptotically.
  On the one hand, by Lemma~\ref{lem:pairing-center-derivative}\ref{item:orbit-map-C2},
  the family $\{\mcU_\xi\}$ is strongly compact in $H^1(\mbBn)$
  when the centers range over a compact set.
  Since $r_k\rightharpoonup0$, a finite net gives uniform convergence
  of $\lag r_k,\mcU_\xi\rag_\lambda$ to $0$ on each such set.
  On the other hand, Lemmas~\ref{lem:duality} and~\ref{lem:overlap} give
  $\lag f_\infty,\mcU_\xi\rag_\lambda\to0$ as
  $d(o,\xi)\to\infty$.
  Since both pairings are uniformly bounded, it follows that
  \begin{equation*}
    \sup_{\xi\in\mbBn}
    \left|
    \lag f_\infty,\mcU_\xi\rag_\lambda
    \lag r_k,\mcU_\xi\rag_\lambda
    \right|\longrightarrow0.
  \end{equation*}
  For two real numbers, the square of their sum differs from the
  larger square by at most three times the absolute value of their
  product. Applying this to the two pairings and taking the supremum
  over $\xi$, the uniform estimate proves \eqref{equ:minimizing-max-pairing}.

  We now compare the masses of $f_\infty$ and $r_k$.
  \begingroup
  By the second line of \eqref{equ:minimizing-weak-splitting}, the limit of $\norm{r_k}_{L^{2^{\ast}}(\mbBn)}$ exists and
  \begin{equation}\label{equ:minimizing-mass-normalization}
    \norm{f_\infty}_{L^{2^{\ast}}(\mbBn)}^{2^{\ast}}
    +\limktoit\norm{r_k}_{L^{2^{\ast}}(\mbBn)}^{2^{\ast}}=1,
  \end{equation}
  and $\norm{f_\infty}_{L^{2^{\ast}}(\mbBn)}>0$.
  \endgroup
  Applying \eqref{equ:minimizing-global} separately to $f_\infty$
  and $r_k$ and adding, then using the first line of \eqref{equ:minimizing-weak-splitting},
  \eqref{equ:minimizing-max-pairing}
  and \eqref{equ:minimizing-equality}, we find
  \begin{equation*}
    \begin{aligned}
      S_\lambda\bigl(\norm{f_\infty}_{L^{2^{\ast}}(\mbBn)}^2+\norm{r_k}_{L^{2^{\ast}}(\mbBn)}^2\bigr)
      &\leqslant
      (1-c_{\mathrm{BE}})\bigl(\norm{f_\infty}_\lambda^2+\norm{r_k}_\lambda^2\bigr)
      +c_{\mathrm{BE}}\bigl(m(f_\infty)^2+m(r_k)^2\bigr)\\
      &=S_\lambda+c_{\mathrm{BE}}\min\{m(f_\infty)^2,m(r_k)^2\}+o(1)\\
      &\leqslant
      S_\lambda+c_{\mathrm{BE}} S_\lambda
      \min\{\norm{f_\infty}_{L^{2^{\ast}}(\mbBn)}^2,\norm{r_k}_{L^{2^{\ast}}(\mbBn)}^2\}+o(1),
    \end{aligned}
  \end{equation*}
  where the last line uses \eqref{equ:minimizing-pairing}.
  Passing to the limit yields
  \begingroup
  \begin{equation}\label{equ:minimizing-mass-comparison}
    \norm{f_\infty}_{L^{2^{\ast}}(\mbBn)}^2+\limktoit\norm{r_k}_{L^{2^{\ast}}(\mbBn)}^2-1
    \leqslant c_{\mathrm{BE}}\min\Bigl\{\norm{f_\infty}_{L^{2^{\ast}}(\mbBn)}^2,\ \limktoit\norm{r_k}_{L^{2^{\ast}}(\mbBn)}^2\Bigr\}.
  \end{equation}

  We compare this with the dichotomy threshold.
  For $s,t\geqslant0$ with $s^{\frac{n}{n-2}}+t^{\frac{n}{n-2}}=1$ we claim
  \begin{equation}\label{equ:minimizing-elementary}
    s+t-1\geqslant\bigl(2-2^{\frac{n-2}{n}}\bigr)\min\{s,t\}.
  \end{equation}
  By symmetry we may assume $t\leqslant s$. Then
  \begin{equation*}
    0\leqslant t\leqslant2^{-\frac{n-2}{n}},\qquad
    s=\bigl(1-t^{\frac{n}{n-2}}\bigr)^{\frac{n-2}{n}}.
  \end{equation*}
  In the interior of this interval,
  \begin{equation*}
    \frac{\mrd^2}{\mrd t^2}
    \bigl(1-t^{\frac{n}{n-2}}\bigr)^{\frac{n-2}{n}}
    =
    -\frac{2}{n-2}\,
    t^{\frac{n}{n-2}-2}
    \bigl(1-t^{\frac{n}{n-2}}\bigr)^{\frac{n-2}{n}-2}<0,
  \end{equation*}
  so the graph of $t\mapsto(1-t^{\frac{n}{n-2}})^{\frac{n-2}{n}}$ lies above the chord joining
  $(0,1)$ and $(2^{-\frac{n-2}{n}},2^{-\frac{n-2}{n}})$.
  Hence
  \begin{equation*}
    s\geqslant1+\bigl(1-2^{\frac{n-2}{n}}\bigr)t,\qquad
    s+t-1\geqslant\bigl(2-2^{\frac{n-2}{n}}\bigr)t,
  \end{equation*}
  which is \eqref{equ:minimizing-elementary}.

  By \eqref{equ:minimizing-mass-normalization}, we may apply \eqref{equ:minimizing-elementary} with $s=\norm{f_\infty}_{L^{2^{\ast}}(\mbBn)}^2$ and $t=\limktoit\norm{r_k}_{L^{2^{\ast}}(\mbBn)}^2$.
  Together with \eqref{equ:minimizing-mass-comparison} and the strict dichotomy bound
  $c_{\mathrm{BE}}<c_{\mathrm{dich}}=2-2^{\frac{n-2}{n}}$ of Theorem~\ref{thm:crit-global-bounds}, this forces the minimum in \eqref{equ:minimizing-mass-comparison} to vanish.
  Since $\norm{f_\infty}_{L^{2^{\ast}}(\mbBn)}>0$, we obtain $\norm{r_k}_{L^{2^{\ast}}(\mbBn)}\to0$ and, by \eqref{equ:minimizing-mass-normalization}, $\norm{f_\infty}_{L^{2^{\ast}}(\mbBn)}=1$.
  \endgroup
  By the second line of \eqref{equ:minimizing-pairing}, $m(r_k)\to0$, and
  thus \eqref{equ:minimizing-max-pairing} gives
  $m(f_k)^2=m(f_\infty)^2+o(1)$.
  Inserting this and the first line of \eqref{equ:minimizing-weak-splitting} into \eqref{equ:minimizing-equality},
  we obtain
  \begin{equation*}
    \begin{aligned}
      (1-c_{\mathrm{BE}})\norm{r_k}_\lambda^2
      &=
      S_\lambda-
      \bigl((1-c_{\mathrm{BE}})\norm{f_\infty}_\lambda^2
          +c_{\mathrm{BE}} m(f_\infty)^2\bigr)+o(1)\\
      &\leqslant o(1),
    \end{aligned}
  \end{equation*}
  where \eqref{equ:minimizing-global} and $\norm{f_\infty}_{L^{2^{\ast}}(\mbBn)}=1$ give the last line.
  Since $c_{\mathrm{BE}}<1$, we have $\norm{r_k}_\lambda\to0$, and
  norm equivalence proves $f_k\to f_\infty$ strongly in $H^1(\mbBn)$.

  \emph{Step 3: the limit is a minimizer outside $\mcZ$.} The limit has $\norm{f_\infty}_{L^{2^{\ast}}(\mbBn)}=1$, so $f_\infty\neq0$.
  Suppose that $f_\infty\in\mcZ$. Then
  \begin{equation*}
    0<\sigma(f_k)
    \leqslant
    \frac{\norm{f_k-f_\infty}_\lambda^2}{\norm{f_k}_\lambda^2}
    \longrightarrow0.
  \end{equation*}
  Hence Proposition~\ref{prop:local-lower-function}, \eqref{equ:const-nu-explicit} and
  \eqref{equ:const-nu-limits} give
  \begin{equation*}
    c_{\mathrm{BE}}
    =\limktoit
    \frac{\delta(f_k)^2}{\dist(f_k,\mcZ)^2}
    \geqslant
    \limktoit\nu(\sigma(f_k))
    =1-\frac{q}{\mu_3}=c_{\mrloc},
  \end{equation*}
  where the last equality is Theorem~\ref{thm:crit-local}.
  This contradicts the strict local bound $c_{\mathrm{BE}}<c_{\mrloc}$ of Theorem~\ref{thm:crit-global-bounds}.
  Hence $f_\infty\notin\mcZ$, and its distance to $\mcZ$ is positive by
  Lemma~\ref{lem:dist-positive}.
  Finally, strong convergence in $H^1(\mbBn)$ and \eqref{equ:PS} make the deficit continuous along $\{f_k\}$, and
  the distance is $1$-Lipschitz in the $\lambda$-norm.
  The stability quotient therefore converges to its value at
  $f_\infty$, which equals $c_{\mathrm{BE}}$. This finishes the proof of Proposition~\ref{prop:minimizing-compact}.
\end{proof}

\begin{proof}[Proof of Theorem~\ref{thm:attainment}]
  By homogeneity, choose a minimizing sequence $\{f_k\}$ for \eqref{equ:cbe-def} with
  $\norm{f_k}_{L^{2^{\ast}}(\mbBn)}=1$.
  Then Proposition~\ref{prop:minimizing-compact} gives $f_\infty\in H^1(\mbBn)\setminus(\mcZ\cup\{0\})$ with
  $\delta(f_\infty)^2=c_{\mathrm{BE}}\dist(f_\infty,\mcZ)^2$. This finishes the proof of Theorem~\ref{thm:attainment}.
\end{proof}

\section*{Acknowledgments}
 The first author is supported by the National Natural
Science Foundation of China (Grant No.~12571127). The second author is supported by the National Natural Science Foundation of China (Grant No.~12526592).

\section*{Conflict of interest}
The authors declare that they have no conflict of interest.

\section*{Declaration on the use of AI tools}
AI-based tools were used for language polishing and to assist in checking the precision of mathematical statements and the correctness of proofs. All the mathematical statements and
their proofs were
reviewed by the authors, who take full responsibility for the content of the
manuscript.


\begin{thebibliography}{99}
  \bibitem{Aubin1976} T.~Aubin, \emph{Probl\`emes isop\'erim\'etriques et espaces de Sobolev}, J. Differential Geometry \textbf{11} (1976), 573--598.

  \bibitem{BGKM2025} M.~Bhakta, D.~Ganguly, D.~Karmakar and S.~Mazumdar, \emph{Sharp quantitative stability of Poincar\'e--Sobolev inequality in the hyperbolic space and applications to fast diffusion flows}, Calc. Var. Partial Differential Equations \textbf{64} (2025), Paper No.~23, 47~pp.

  \bibitem{BGKM2025b} M.~Bhakta, D.~Ganguly, D.~Karmakar and S.~Mazumdar, \emph{Sharp quantitative stability of Struwe's decomposition of the Poincar\'e--Sobolev inequalities on the hyperbolic space: Part~I}, Adv. Math. \textbf{479} (2025), Paper No.~110447, 84~pp.

  \bibitem{BhaktaSandeep2012} M.~Bhakta and K.~Sandeep, \emph{Poincar\'e--Sobolev equations in the hyperbolic space}, Calc. Var. Partial Differential Equations \textbf{44} (2012), 247--269.

  \bibitem{BianchiEgnell1991} G.~Bianchi and H.~Egnell, \emph{A note on the Sobolev inequality}, J. Funct. Anal. \textbf{100} (1991), 18--24.

  \bibitem{BrezisKato1979} H.~Brezis and T.~Kato, \emph{Remarks on the Schr\"odinger operator with singular complex potentials}, J. Math. Pures Appl. (9) \textbf{58} (1979), 137--151.

  \bibitem{BrezisLieb1983} H.~Brezis and E.~H.~Lieb, \emph{A relation between pointwise convergence of functions and convergence of functionals}, Proc. Amer. Math. Soc. \textbf{88} (1983), 486--490.

  \bibitem{BrezisLieb1985} H.~Brezis and E.~H.~Lieb, \emph{Sobolev inequalities with remainder terms}, J. Funct. Anal. \textbf{62} (1985), 73--86.

  \bibitem{BrezisNirenberg1983} H.~Brezis and L.~Nirenberg, \emph{Positive solutions of nonlinear elliptic equations involving critical Sobolev exponents}, Comm. Pure Appl. Math. \textbf{36} (1983), 437--477.

  \bibitem{CKW2026} H.~Chen, S.~Kim and J.~Wei, \emph{Sharp quantitative stability estimates for the Brezis--Nirenberg problem}, J. Funct. Anal. \textbf{291} (2026), Paper No.~111515, 61~pp.

  \bibitem{CLTa2024} L.~Chen, G.~Lu and H.~Tang, \emph{Stability of Hardy--Littlewood--Sobolev inequalities with explicit lower bounds}, Adv. Math. \textbf{450} (2024), Paper No.~109778, 28~pp.

  \bibitem{CLTc2025} L.~Chen, G.~Lu and H.~Tang, \emph{Optimal asymptotic lower bound for stability of fractional Sobolev inequality and the global stability of log-Sobolev inequality on the sphere}, Adv. Math. \textbf{479} (2025), Paper No.~110438, 32~pp.

  \bibitem{CLTb2026} L.~Chen, G.~Lu and H.~Tang, \emph{Optimal stability of Hardy--Littlewood--Sobolev and Sobolev inequalities of arbitrary orders with dimension-dependent constants}, Math. Ann. \textbf{394} (2026), Paper No.~77, 46~pp.

  \bibitem{CLTW2026} L.~Chen, G.~Lu, H.~Tang and B.~Wang, \emph{Asymptotically sharp stability of Sobolev inequalities on the Heisenberg group with dimension-dependent constants}, J. Math. Pures Appl. (9) \textbf{206} (2026), Paper No.~103832, 29~pp.

  \bibitem{CFW2013} S.~Chen, R.~L.~Frank and T.~Weth, \emph{Remainder terms in the fractional Sobolev inequality}, Indiana Univ. Math. J. \textbf{62} (2013), 1381--1397.

  \bibitem{CFMP2009} A.~Cianchi, N.~Fusco, F.~Maggi and A.~Pratelli, \emph{The sharp Sobolev inequality in quantitative form}, J. Eur. Math. Soc. (JEMS) \textbf{11} (2009), 1105--1139.

  \bibitem{CFM2018} G.~Ciraolo, A.~Figalli and F.~Maggi, \emph{A quantitative analysis of metrics on $\mbRn$ with almost constant positive scalar curvature, with applications to fast diffusion flows}, Int. Math. Res. Not. IMRN \textbf{2018}, 6780--6797.

  \bibitem{DSW2025} B.~Deng, L.~Sun and J.~Wei, \emph{Sharp quantitative estimates of Struwe's decomposition}, Duke Math. J. \textbf{174} (2025), 159--228.

  \bibitem{DEFFL2025} J.~Dolbeault, M.~J.~Esteban, A.~Figalli, R.~L.~Frank and M.~Loss, \emph{Sharp stability for Sobolev and log-Sobolev inequalities, with optimal dimensional dependence}, Camb. J. Math. \textbf{13} (2025), 359--430.

  \bibitem{ENS2022} M.~Engelstein, R.~Neumayer and L.~Spolaor, \emph{Quantitative stability for minimizing Yamabe metrics}, Trans. Amer. Math. Soc. Ser. B \textbf{9} (2022), 395--414.

  \bibitem{FigalliGlaudo2020} A.~Figalli and F.~Glaudo, \emph{On the sharp stability of critical points of the Sobolev inequality}, Arch. Ration. Mech. Anal. \textbf{237} (2020), 201--258.

  \bibitem{FigalliNeumayer2019} A.~Figalli and R.~Neumayer, \emph{Gradient stability for the Sobolev inequality: the case $p\geqslant2$}, J. Eur. Math. Soc. (JEMS) \textbf{21} (2019), 319--354.

  \bibitem{FigalliZhang2022} A.~Figalli and Y.~R.--Y.~Zhang, \emph{Sharp gradient stability for the Sobolev inequality}, Duke Math. J. \textbf{171} (2022), 2407--2459.

  \bibitem{Frank2022} R.~L.~Frank, \emph{Degenerate stability of some Sobolev inequalities}, Ann. Inst. H. Poincar\'e C Anal. Non Lin\'eaire \textbf{39} (2022), 1459--1484.

  \bibitem{Frank2024} R.~L.~Frank, \emph{The sharp Sobolev inequality and its stability: an introduction}, in: Geometric and Analytic Aspects of Functional Variational Principles, pp.~1--64, Lecture Notes in Math. 2348, CIME Found. Subser., Springer, Cham, 2024.

  \bibitem{GangulySandeep2015} D.~Ganguly and K.~Sandeep, \emph{Nondegeneracy of positive solutions of semilinear elliptic problems in the hyperbolic space}, Commun. Contemp. Math. \textbf{17} (2015), Paper No.~1450019, 13~pp.

  \bibitem{GilbargTrudinger2001} D.~Gilbarg and N.~S.~Trudinger, \emph{Elliptic Partial Differential Equations of Second Order}, Reprint of the 1998 edition, Classics Math., Springer, Berlin, 2001. xiv+517 pp.

  \bibitem{Koenig2023a} T.~K\"onig, \emph{On the sharp constant in the Bianchi--Egnell stability inequality}, Bull. Lond. Math. Soc. \textbf{55} (2023), 2070--2075.

  \bibitem{Koenig2025} T.~K\"onig, \emph{Stability for the Sobolev inequality: existence of a minimizer}, J. Eur. Math. Soc. (JEMS) (2025), https://doi.org/10.4171/JEMS/1582.

  \bibitem{ManciniSandeep2008} G.~Mancini and K.~Sandeep, \emph{On a semilinear elliptic equation in $\mbH^n$}, Ann. Sc. Norm. Super. Pisa Cl. Sci. (5) \textbf{7} (2008), 635--671.

  \bibitem{McKean1970} H.~P.~McKean, \emph{An upper bound to the spectrum of $\Delta$ on a manifold of negative curvature}, J. Differential Geometry \textbf{4} (1970), 359--366.

  \bibitem{Neumayer2020} R.~Neumayer, \emph{A note on strong-form stability for the Sobolev inequality}, Calc. Var. Partial Differential Equations \textbf{59} (2020), Paper No.~25, 8~pp.

  \bibitem{Nobili2026} F.~Nobili, \emph{An overview of the stability of Sobolev inequalities on Riemannian manifolds with Ricci lower bounds}, Indag. Math. (N.S.) \textbf{37} (2026), 960--992.

  \bibitem{NobiliViolo2024} F.~Nobili and I.~Y.~Violo, \emph{Stability of Sobolev inequalities on Riemannian manifolds with Ricci curvature lower bounds}, Adv. Math. \textbf{440} (2024), Paper No.~109521, 58~pp.

  \bibitem{Struwe1984} M.~Struwe, \emph{A global compactness result for elliptic boundary value problems involving limiting nonlinearities}, Math. Z. \textbf{187} (1984), 511--517.

  \bibitem{StruweVM} M.~Struwe, \emph{Variational Methods: Applications to Nonlinear Partial Differential Equations and Hamiltonian Systems}, Second edition, Ergeb. Math. Grenzgeb. (3) 34, Springer, Berlin, 1996. xvi+272 pp.

  \bibitem{Talenti1976} G.~Talenti, \emph{Best constant in Sobolev inequality}, Ann. Mat. Pura Appl. (4) \textbf{110} (1976), 353--372.

  \bibitem{WeiWu2022} J.~Wei and Y.~Wu, \emph{On the stability of the Caffarelli--Kohn--Nirenberg inequality}, Math. Ann. \textbf{384} (2022), 1509--1546.

  \bibitem{Zettl2005} A.~Zettl, \emph{Sturm--Liouville Theory}, Math. Surveys Monogr. 121, American Mathematical Society, Providence, RI, 2005. xii+328 pp.
  
\end{thebibliography}
\end{document}